\documentclass[12pt]{article}
\usepackage{amsmath,amssymb,amsthm,mathrsfs, esint}
\usepackage{titlesec,hyperref}
\usepackage{color}
\usepackage{fancyhdr}
\usepackage{cases}
\usepackage[margin=2.5cm]{geometry}
\usepackage{enumerate}
\usepackage{tikz}
\usepackage{mathdots}
\usepackage{yhmath}
\usepackage{cancel}
\usepackage{siunitx}
\usepackage{array}
\usepackage{multirow}
\usepackage{amssymb}
\usepackage{gensymb}
\usepackage{tabularx}
\usepackage{extarrows}
\usepackage{booktabs}
\usetikzlibrary{fadings}
\usetikzlibrary{patterns}
\usetikzlibrary{shadows.blur}
\usetikzlibrary{shapes}
\usepackage[integrals]{wasysym}
\usepackage{bm}
\usepackage{graphicx}
\usepackage{caption}
\usepackage{float}
\usepackage{subfigure}
\usepackage{indentfirst}
\usepackage[title]{appendix}
\usepackage[numbers,sort&compress]{natbib}

\newcommand{\dd}{\mathop{}\!\mathrm{d}}
\let\del\partial

\newcommand{\supp}{\text{supp}\,}

\newcommand{\RR}{\mathring R}

\newcommand{\ZZ}{\mathbb Z}

\newcommand{\RRR}{\mathring{\bar R}}
\newcommand{\ootimes}{\mathbin{\mathring{\otimes}}}

\let\div\relax
\DeclareMathOperator{\div}{div}
\DeclareMathOperator{\tr}{Tr}
\newcommand{\vin}[0]{v^{\textup{in}}}

\newcommand{\pex}[0]{\textsf{\textit{p}}}
\newcommand{\vex}[0]{\textsf{\textit{v}}}
\newcommand{\vv}[0]{\bar v}

\newcommand{\ppp}[0]{\bar{\bar p}}
\newcommand{\pp}[0]{\bar p}

\newcommand{\wpq}[1][q+1]{w^{\textup{(p)}}_{#1}}
\newcommand{\wcq}[1][q+1]{w^{\textup{(c)}}_{#1}}

\newcommand{\wtq}[1][q+1]{w^{\textup{(t)}}_{#1}}

\newcommand{\ptq}[1][q+1]{P^{\textup{(t)}}_{#1}}
\newcommand{\Rnash}{R_{q+1}^\textup{Nash}}
\newcommand{\Rtransport}{R_{q+1}^\textup{trans}}
\newcommand{\Rosc}{R_{q+1}^\textup{osc}}

\newcommand{\TTT}[0]{\mathbb{T}}
\newcommand{\Ttwo}[0]{\mathbb{T}^2}

\newcommand{\matd}[1]{D_{t,#1}}

\newcommand{\tRR}{{\mathring{\widetilde R} }}

\newcommand{\tv}{ {\widetilde{{v}}} }

\def\dashint{\,\ThisStyle{\ensurestackMath{%
			\stackinset{c}{.2\LMpt}{c}{.5\LMpt}{\SavedStyle-}{\SavedStyle\phantom{\int}}}%
		\setbox0=\hbox{$\SavedStyle\int\,$}\kern-\wd0}\int}
\def\ddashint{\,\ThisStyle{\ensurestackMath{%
			\stackinset{c}{.2\LMpt}{c}{.5\LMpt+.2\LMex}{\SavedStyle-}{%
				\stackinset{c}{.2\LMpt}{c}{.5\LMpt-.2\LMex}{\SavedStyle-}{%
					\SavedStyle\phantom{\int}}}}\setbox0=\hbox{$\SavedStyle\int\,$}\kern-\wd0}\int}
\usepackage{scalerel,stackengine}
\stackMath
\newcommand\widecheck[1]{%
	\savestack{\tmpbox}{\stretchto{%
			\scaleto{%
				\scalerel*[\widthof{\ensuremath{#1}}]{\kern-.6pt\bigwedge\kern-.6pt}%
				{\rule[-\textheight/2]{1ex}{\textheight}}
			}{\textheight}%
		}{0.5ex}}%
	\stackon[1pt]{#1}{\scalebox{-1}{\tmpbox}}%
}

\newcommand{\coloneq}{\mathrel{\mathop:}=}

\newtheorem{thm}{Theorem}[section]

\newtheorem{cor}[thm]{Corollary}

\newtheorem{lem}[thm]{Lemma}
\newtheorem{prop}[thm]{Proposition}

\theoremstyle{definition}
\newtheorem{defn}[thm]{Definition}

\theoremstyle{remark}
\newtheorem{rem}[thm]{Remark}
\numberwithin{equation}{section}

\allowdisplaybreaks

\begin{document}
	\title{Dissipative solutions for 2D Onsager's conjecture \thanks {Mathematics Subject Classification: 35Q31,~35A02, ~35D30, ~76B03.}}
	\author{Lili $\mbox{Du}^{1}$\quad
		Xinliang $\mbox{Li}^{2}$ \quad
		Weikui $\mbox{Ye}^{3}$
	}
	\date{}
	\maketitle
	\begin{abstract}

In 1949, L. Onsager conjectured that weak solutions of the Euler equations in $C_{t,x}^{\gamma}$  are energy conservative when $\gamma>\frac{1}{3}$, whereas nonconservative solutions exist when $\gamma<\frac{1}{3}$.
Constantin, E and
	Titi [Commun. Math. Phys.,
165 (1), 1994] proved the conservation part, while Isett [Ann. of Math., 188 (3), 2018]
	constructed nontrivial $C^\gamma_{t,x}$ weak solutions with compact
	support in time in three dimensions for each $0<\gamma<\frac{1}{3}$.
	Building on this breakthrough
	work, Buckmaster, De Lellis,
	Sz\'{e}kelyhidi Jr. and Vicol [Commun. Pure Appl. Math., 72 (2), 2019] obtained strictly dissipative
	solutions below the Onsager threshold. In two dimensions, Giri and
	Radu [Invent. Math., 238 (2), 2024] first constructed a nontrivial weak solution in
	$C^{\gamma}(\mathbb T^2\times[0,T])$  for $\gamma<\frac{1}{3}$ and left open the question whether it would be possible to construct 
	weak solutions with prescribed  energy profiles below the Onsager threshold in two dimensions.
	
	In this paper, we resolve this question and thereby obtain infinitely
	many strictly dissipative weak solutions in
	$C^\gamma(\mathbb T^2\times[0,T])$ for every $0<\gamma<\frac{1}{3}$.
Moreover, we prove that the associated wild initial data are dense in
	the divergence-free subspace of $C^\gamma(\mathbb T^2)$
	endowed with the $C^{\gamma'}$ topology for every $0<\gamma'<\gamma<\frac13$. The key new ingredient is a family of staircase traveling waves, which we combine with a Picard iteration coupled to the Newton--Nash
	scheme to control both the Reynolds stress and the kinetic energy.

	\end{abstract}
	\textit{\emph{Keywords}}: Onsager’s conjecture; dissipative solutions; Euler equations;  convex integration.
	
	\vskip   0.2cm \footnotetext[1]{
		Department of Mathematics, Sichuan University, Chengdu 610064, P. R. China. Email address: dulili@scu.edu.cn.}
	
	\vskip   0.2cm \footnotetext[2]{
		Department of Mathematics, Shantou University, Shantou 515000, P. R. China.
		Email address: xinliangli@stu.edu.cn.}
	
	\vskip   0.2cm \footnotetext[3]{
		School of Mathematical Sciences, Shenzhen University, Shenzhen, 518060, P.R. China.
		Email address: yeweikuisz@szu.edu.cn.}
	
	\setlength{\baselineskip}{17pt}
	\begin{center}
		\tableofcontents
	\end{center} 
	\section{Introduction and main results}
	In this paper, we study strict energy dissipation in the context
of Onsager's conjecture for the two-dimensional incompressible
Euler equations:
	\begin{equation}
		\left\{
		\begin{aligned}
			\partial_t v+\operatorname{div}(v\otimes v)+\nabla p&=0,\\
			\operatorname{div}v&=0,\\
			v|_{t=0}&=\vin
		\end{aligned}
		\right.
		\label{e:E}
	\end{equation}
	on the torus $\mathbb T^{2}$, where $v:\mathbb T^2\times[0,T]\rightarrow\mathbb R^2$ is the velocity field of the fluid and $p:\mathbb T^2\times[0,T]\rightarrow\mathbb R$ the pressure field.
	
	 In the celebrated work \cite{Onsa}, Onsager  proposed  that turbulent energy dissipation could persist without viscosity and stated in the final paragraph on p.~286 that velocity fields with H\"older regularity above $\frac{1}{3}$ conserve energy:
	
	\begin{quotation}
		\textit{It is of some interest to note that in principle, turbulent dissipation as described could take place just as readily without the final assistance by viscosity. In the absence of viscosity, the standard proof of the conservation of energy does not apply, because the velocity field does not remain differentiable! In fact it is possible to show that the velocity field in such ``ideal" turbulence cannot obey any LIPSCHITZ condition of the form}
		$$
		\left|\vec{v}\left(\overrightarrow{r'}+\vec{r}\right)-\vec{v}\left(\overrightarrow{r'}\right)\right|<\text { (const.) } r^n,
		$$
		\textit{for any order $n$ greater than $1 / 3$; otherwise the energy is conserved. $\cdots$ }
	\end{quotation}
	
	The same exponent  $\frac13$ also arises from the scaling of velocity increments in Kolmogorov's  theory \cite{KOL41,KOL411}. These ideas led to the modern form of Onsager's conjecture in \cite{DS13} by De Lellis and Sz\'{e}kelyhidi:
	\begin{enumerate}
		\item[(a)] Any weak solution $v$ of the Euler equations belonging to the H\"{o}lder space $C_{t, x}^\gamma$  is energy conservative when $\gamma>\frac{1}{3}$ ;\label{con-a}
		\item[(b)] \label{con-b} There exist  dissipative weak solutions $v \in C_{t,x}^\gamma$ for any $\gamma<\frac{1}{3}$.
	\end{enumerate}

	The first progress on the rigid part~(a) was made by Eyink \cite{EGL94}, who proved energy conservation under a slightly stronger $*$-H\"{o}lder condition. Constantin, E and Titi \cite{CET} subsequently showed
	that weak solutions $v\in  L_{t}^3B^\gamma_{3,\infty}\cap C_{t}L_{x}^{2}$ conserve energy
	whenever $\gamma>\frac{1}{3}$, independently of the spatial dimension. At the endpoint, Cheskidov, Constantin, Friedlander and Shvydkoy \cite{CCFS} obtained the criterion in $L^3_t B^{\frac{1}{3}}_{3,c_0}$.
		The consideration for part~(b) began with Scheffer's construction  \cite{S93} of a nonzero weak solution in the plane with compact support in space--time and Shnirelman's construction  \cite{Sh97} of a weak solution on $\mathbb T^2$ with compact support in time. In their seminal paper \cite{DS09}, De Lellis and Sz\'{e}kelyhidi gave a groundbreaking differential-inclusion formulation of the multidimensional incompressible Euler equations, combining Baire-category arguments (see, e.g., \cite{DM97}) with convex integration (see, e.g., \cite{MS}). Following a series of fundamental technical advances \cite{DS13,DS14,DS10,B15,BDIS15,BDS16,Ise17,DS17,IO15}, Isett \cite{Ise18} constructed nontrivial $C^{\frac{1}{3}-}_{t,x}$ weak solutions with compact
		support in time, which
		therefore fail to conserve kinetic energy. Buckmaster, De Lellis, Sz\'{e}kelyhidi Jr. and Vicol \cite{BDSV} obtained strictly dissipative
		solutions in three dimensions below the Onsager threshold.  More recently, Isett \cite{Isett2024Endpoint} substantially refined the construction toward the still-open endpoint by optimizing the frequency growth and avoiding power losses in the iteration, thereby producing nonconservative three-dimensional Euler flows whose regularity approaches the critical exponent $\frac{1}{3}$ at small spatial scales. For the forced Euler equations on $\mathbb T^d$, $d\geq3$, Bulut, Huynh and Palasek \cite{BU23} constructed nonunique, almost-everywhere smooth solutions of H\"{o}lder regularity $\frac{1}{2}-$. Enciso, Pe\~{n}afiel-Tom\'{a}s, and Peralta-Salas \cite{EPP1} showed that a smooth Euler flow on $\Omega\times(0,T)$ can be extended to an admissible weak solution on $\mathbb R^3\times(0,T)$ of class $C^\beta$ for all $\beta<\frac{1}{3}$, with controlled spatial support. In two dimensions, Giri and Radu \cite{VGR23}  cleverly introduced an innovative Newton--Nash iteration to construct a nontrivial weak solution in $C^\gamma(\mathbb T^2\times[0,T])$ with compact support in time for any $\gamma<\frac{1}{3}$. For the forced 2D Euler equations, Castro, Faraco, Mengual and
	Solera \cite{CFMS} gave a simpler proof of Vishik's nonuniqueness
	theorem with zero initial vorticity. In the unforced periodic setting,
	Bru\`{e}, Colombo, and Kumar \cite{BCK24} extended the Lorentz-space
	result of Bru\`{e} and Colombo \cite{BC23} to nonconservative solutions
	with vorticity in $C_tL_x^p$ for some $p>1$. For each $2<p<\infty$, Mengual
	\cite{MF24} constructed an initial velocity on $\mathbb R^2$ with
	vorticity in $L_x^1\cap L_x^p$ that admits infinitely many bounded
	admissible weak solutions.

	In particular, Giri and Radu  proposed an open problem in Remark 2.3 of \cite{VGR23}(the reference numbers in quotation are those of  \cite{VGR23}):
		\begin{quotation}
		\noindent\textit{\textbf{Remark 2.3}}\
		\textit{In [6], the authors modify the iteration used in [29] and show that any
		given smooth energy profile $e:[0,1]\to\mathbb{R}_{>0}$ can be achieved
		by three-dimensional flexible Euler flows: i.e. the solution can be
		designed to satisfy
		$\int_{\mathbb{T}^3}|u(x,t)|^2\,dx=e(t)$.
		It is an interesting question whether the scheme in the present paper
		admits such a modification. In order to point out the difficulties,
		let us first briefly recall the the ideas used in [6].
		\ldots}
		
			\ldots\quad 	\ldots
		
		\textit{While we do believe that the energy profile can be prescribed, the considerations above show that constructing solutions satisfying
		this notion of flexibility raises important technical issues, and,
		as such, lies beyond the scope of this paper.}
		\end{quotation} 
	which can be concluded as:
		\begin{center}
{\bf Would it be possible to construct weak solutions with prescribed decreasing energy profiles in
$C^\gamma(\mathbb T^2\times[0,T])$
for every $0<\gamma<\frac13$?}
		\end{center} 
		
		Motivated by the problem raised by Giri and Radu in \cite{VGR23}, we give an affirmative answer in this paper by introducing a family of \textit{ staircase traveling waves}, the specific conclusions are as follows:
	
	\begin{thm}\label{main0}
		
	Let $T>0$, $0<\gamma<\frac13$, and let $e(t):[0,T]\to(0,\infty)$ be a smooth function. Then there exist a continuous vector field
			$
			v\in C^{\gamma}(\mathbb T^2\times[0,T])
			$, and a continuous scalar pressure $p$ which solve the incompressible Euler equations
           \eqref{e:E}
		 in the sense of distributions and such that $$e(t)=\int_{\mathbb T^2}|v(x,t)|^2\,\dd x$$
		for all $t\in[0,T]$.

	\end{thm}

		Since Theorem~\ref{main0} prescribes the energy exactly, choosing a strictly decreasing profile $e(t)$ yields
		\[
		\int_{\mathbb T^2}|v(x,t)|^2\,\dd x
		<
		\int_{\mathbb T^2}|v(x,s)|^2\,\dd x,
		\qquad 0\leq s<t\leq T.
		\]
This resulting solution of \eqref{e:E} satisfies the strong energy inequality and is strictly dissipative. Moreover, every solution on $\mathbb T^2$ can be viewed naturally as a solution on $\mathbb T^d$ for $d\geq2$. Theorem~\ref{main0} therefore yields strictly dissipative weak solutions in $C^\gamma(\mathbb T^d\times[0,T])$
for $0<\gamma<\frac13$ and $d\geq2$.
	
	On the other hand, Daneri, Runa and
	Sz\'{e}kelyhidi Jr. introduced the following definition in \cite{DRS21}:
	\begin{defn}\label{WID}
	For $d\geq2$, given a divergence-free vector field
	$v^{\mathrm{in}}\in C^{\gamma_0}(\mathbb T^d)$, $v^{\mathrm{in}}$ is a
	wild initial datum in $C^\gamma$ if the Euler equations on
	$\mathbb T^d\times[0,T]$ admit infinitely many weak solutions $v$
	with initial datum $v^{\mathrm{in}}$ such that
	\begin{equation}
	\int_{\mathbb T^d}|v(x,t)|^2\,\dd x
	\leq
	\int_{\mathbb T^d}|v^{\mathrm{in}}|^2\,\dd x
	\qquad\text{for every }t\in[0,T],
	\label{1.2}
	\end{equation}
	and
	\begin{equation}
	|v(x,t)-v(y,t)|
	\leq C|x-y|^\gamma
	\qquad
	\text{for every }t\in[0,T]\text{ and }x,y\in\mathbb T^d.
		\label{1.3}
	\end{equation}
	\end{defn}
	Daneri, Runa and Sz\'{e}kelyhidi Jr. \cite{DRS21} proved that, for every $0<\gamma<\frac{1}{3}$, the divergence-free wild initial data in $C^\gamma(\mathbb T^3)$  form an $L^2$-dense subset of the divergence-free subspace of $L^2(\mathbb T^3)$.	In this paper, Theorem~\ref{main0} provides the construction of strictly dissipative weak solutions, a natural further question is how large the corresponding class of wild initial data is.  We address this
	second question as follows.

\begin{thm}\label{main1}
	\label{cor0}
	For any $0<\gamma'<\gamma<\frac13$, the set of divergence-free
	vector fields $v_w^{\mathrm{in}}$ which are wild initial data
	in $C^\gamma$ in the sense of Definition~\ref{WID} is dense
	in the divergence-free subspace of
	$C^\gamma(\mathbb T^2)$ endowed with the $C^{\gamma'}$ topology.
\end{thm}

	\begin{rem}\label{rem1}
		\emergencystretch=1em\relax
			Daneri, Runa and Sz\'{e}kelyhidi Jr. \cite{DRS21} innovatively used a double convex-integration scheme to prove that for any $\epsilon>0$ and divergence-free vector field $v^{\mathrm{in}}\in L^2(\mathbb T^3)$, there exists a wild initial datum $v_w^{\mathrm{in}}$ such that $\|v_w^{\mathrm{in}}-v^{\mathrm{in}}\|_{L^2(\mathbb T^3)}
			\leq\epsilon.
			$
			In two dimensions, we establish density in $C^{\gamma'}$ topology, for any $\epsilon>0$, divergence-free vector field $v^{\mathrm{in}}\in C^\gamma(\mathbb T^2)$ and $0<\gamma'<\gamma<\frac13$, there exists a wild initial datum $v_w^{\mathrm{in}}$ such that
			$
			\|v_w^{\mathrm{in}}-v^{\mathrm{in}}\|_{C^{\gamma'}(\mathbb T^2)}
			\leq\epsilon.
			$
			This does not follow directly from interpolation, since the available $C^{\gamma}$ bound depends on $\epsilon$ and therefore does not yield the required $C^{\gamma'}$ estimate. Due to the nonseparability of $C^\gamma(\mathbb T^2)$,
our approximation argument yields density only in the weaker
$C^{\gamma'}$ topology.
An analogous density result may hold in the Besov spaces
$B^\gamma_{\infty,r}(\mathbb T^2)$, $1\leq r<\infty$,
without lowering the regularity exponent.

	\end{rem}

	The next section describes the principal ideas underlying the proofs of Theorems~\ref{main0} and~\ref{main1}.
	\section{Overview of the approach  }
	
To construct strictly dissipative weak solutions of the two-dimensional
Euler equations in
$C^\gamma(\mathbb T^2\times[0,T])$
for $0<\gamma<\frac13$,
we need to control the kinetic energy as well as the Reynolds stress
throughout the iteration. In three dimensions, the Mikado flows introduced
by Daneri and Sz\'{e}kelyhidi Jr. \cite{DS17}, whose supports can be chosen to
be disjoint, allow one to control both the Reynolds stress and the energy
increment. In two dimensions, however, such a spatial separation of
building blocks in different directions is not available, since
non-parallel tubes necessarily intersect. To overcome this difficulty,
Giri and Radu \cite{VGR23} skillfully introduced high-frequency time oscillations
into the Newton--Nash iteration, choosing disjoint time supports so that
different building blocks do not interact.

 As stated in Remark 2.3 of \cite{VGR23}, when one seeks to prescribe the energy profile,  it seems necessary to replace the purely time-dependent distribution $g_{k,n}(\mu_{q+1}t)$ with a traveling wave of the form $g_{k,n}(\mu_{q+1}t, x)$. A more serious difficulty then arises: \textit{how to control the energy effectively
throughout the iteration while preserving disjoint supports for
these waves at different Newton levels}. In this paper, we resolve the
energy-prescription problem posed by Giri and Radu in
\cite[Remark~2.3]{VGR23} by introducing the following staircase
profiles:

 Let $\Lambda$ be the finite set of
directions used in the geometric decomposition, $\xi\notin\Lambda$ be an
auxiliary direction, and let $h_\eta$, $\eta\in\Lambda\cup\{\xi\}$, be smooth
$1$-periodic functions satisfying
\[
\int_{\mathbb T}h_\eta^2(y)\,\dd y=1,
\qquad
\operatorname{supp}h_\eta
\cap
\operatorname{supp}h_{\eta'}
=\emptyset
\qquad
\text{whenever }\eta\neq\eta'.
\]
For $0\leq n\leq L-1$, $k\in\Lambda$ and $\rho\in \mathbb N$, where $L$ denotes the number of Newton levels,
we then define the $1$-periodic staircase profile
\[
g_{k,n}(y)
=
h_k(\rho^n y)
\prod_{\gamma=0}^{n-1}
h_\xi(\rho^\gamma y),
\qquad y\in\mathbb T,
\]
where the empty product is understood to be $1$.
The factor $h_\xi$ may be viewed as a continuation factor,
while $h_k$ plays the role of a termination factor.
The first few profiles therefore have the schematic form:

\begin{figure}[H]
	\centering
	\[
		\begin{array}{c|ccccc}
			&\textup{level $0$}&\textup{level $1$}
			&\textup{level $2$}&\textup{level $3$}&\cdots\\
			\hline
			g_{k,0}&h_k(y)&&&&\\
			g_{k,1}&h_\xi(y)&h_k(\rho y)&&&\\
			g_{k,2}&h_\xi(y)&h_\xi(\rho y)
			&h_k(\rho^2 y)&&\\
			g_{k,3}&h_\xi(y)&h_\xi(\rho y)
			&h_\xi(\rho^2 y)&h_k(\rho^3 y)&
		\end{array}
	\]
	\caption{Schematic structure of the staircase profiles.}
	\label{fig:staircase-profiles}
\end{figure}
As illustrated in Figure~\ref{fig:staircase-profiles}, the staircase
 	profiles have pairwise disjoint supports:
 	\[
 	\operatorname{supp}g_{k,n}
 	\cap
 	\operatorname{supp}g_{k',n'}
 	=\emptyset
 	\qquad\text{whenever }(k,n)\neq(k',n').
 	\]
   The traveling phase
	and staircase structure allow us to  control the energy at each fixed time and at every Newton level
	while preserving pointwise non-interaction, so that the energy
	corrector $h_q(t)$ can be chosen after the leading energy
	contributions have been identified. Thus, in our construction, the staircase mechanism plays a role
		similar to that of disjoint Mikado flows in higher dimensions.

It remains to realize the prescribed kinetic energy.
The principal perturbation $w_{q+1}^{(p_1)}(t,x)$ corresponding to the level $n=0$ has the leading energy contribution
\[
\int_{\mathbb T^2}
|w_{q+1}^{(p_1)}(t,x)|^2\,\dd x
=
2h_q(t)
\sum_i
\int_{\mathbb T^2}
\widetilde\chi_i^2(t,x)\,\dd x
+
\textup{lower-order terms}.
\]
Thus, the scalar energy corrector $h_q(t)$ is introduced to compensate for the gap
between the prescribed energy profile $e(t)$ and the kinetic energy of
glued velocity $\bar v_q$. However, $h_q$ cannot be chosen independently of the perturbation,
since the latter depends on $h_q$, which must in turn be determined
from the resulting energy balance. This mutual dependence leads
to the following nonlinear feedback loop
\[
h_q
\longmapsto
w_{q+1}
\longmapsto
\int_{\mathbb T^2}
|\bar v_q+w_{q+1}|^2\,\dd x
\longmapsto
h_q.
\]
To close this loop, we use the \textit{Picard-type iteration}
\[
h_q^{(m)}
\longmapsto
\bigl(
w_{q+1}^{(m)},
\mathring R_{q+1}^{(m)}
\bigr)
\longmapsto
h_q^{(m+1)}.
\]
For each fixed candidate function $h_q^{(m)}$, the corresponding perturbation is
constructed through the Newton--Nash scheme. Uniform estimates then allow
us to pass to the limiting triple 
\[
\bigl(
h_q,
w_{q+1},
\mathring R_{q+1}
\bigr),
\]
as $m\rightarrow\infty$, which simultaneously satisfies the next Reynolds-stress estimate and the
required energy estimate. 

Therefore, for a given energy corrector, the Newton-Nash iteration
	reduces the Reynolds stress, while the Picard iteration determines
	that corrector from the resulting energy balance. These two
	procedures form the multiple iteration scheme used in the paper. Please see the following Figure~\ref{MIS1}.

\begin{figure}[H]
	
	\begin{center}
		\tikzset{every picture/.style={line width=0.75pt}} 
		
		\begin{tikzpicture}[x=0.75pt,y=0.75pt,yscale=-1,xscale=1]
			
			\draw    (140.4,45.83) -- (296.72,46.11) ;
			\draw [shift={(298.72,46.11)}, rotate = 180.1] [color={rgb, 255:red, 0; green, 0; blue, 0 }  ][line width=0.75]    (6.56,-1.97) .. controls (4.17,-0.84) and (1.99,-0.18) .. (0,0) .. controls (1.99,0.18) and (4.17,0.84) .. (6.56,1.97)   ;
			\draw    (140.85,135.17) -- (297.17,135.44) ;
			\draw [shift={(299.17,135.45)}, rotate = 180.1] [color={rgb, 255:red, 0; green, 0; blue, 0 }  ][line width=0.75]    (6.56,-1.97) .. controls (4.17,-0.84) and (1.99,-0.18) .. (0,0) .. controls (1.99,0.18) and (4.17,0.84) .. (6.56,1.97)   ;
			\draw    (141.51,225.67) -- (297.83,225.94) ;
			\draw [shift={(299.83,225.95)}, rotate = 180.1] [color={rgb, 255:red, 0; green, 0; blue, 0 }  ][line width=0.75]    (6.56,-1.97) .. controls (4.17,-0.84) and (1.99,-0.18) .. (0,0) .. controls (1.99,0.18) and (4.17,0.84) .. (6.56,1.97)   ;
			\draw   (224.77,342.31) -- (225.77,375.83) -- (229.61,375.72) -- (222.58,398.29) -- (214.22,376.17) -- (218.07,376.06) -- (217.08,342.54) -- cycle ;
			\draw    (140.85,326.33) -- (297.17,326.61) ;
			\draw [shift={(299.17,326.61)}, rotate = 180.1] [color={rgb, 255:red, 0; green, 0; blue, 0 }  ][line width=0.75]    (6.56,-1.97) .. controls (4.17,-0.84) and (1.99,-0.18) .. (0,0) .. controls (1.99,0.18) and (4.17,0.84) .. (6.56,1.97)   ;
			\draw    (330.1,155.62) -- (133.73,210.43) ;
			\draw [shift={(131.8,210.97)}, rotate = 344.4] [color={rgb, 255:red, 0; green, 0; blue, 0 }  ][line width=0.75]    (10.93,-3.29) .. controls (6.95,-1.4) and (3.31,-0.3) .. (0,0) .. controls (3.31,0.3) and (6.95,1.4) .. (10.93,3.29)   ;
			\draw    (330.1,66.72) -- (133.73,121.53) ;
			\draw [shift={(131.8,122.07)}, rotate = 344.4] [color={rgb, 255:red, 0; green, 0; blue, 0 }  ][line width=0.75]    (10.93,-3.29) .. controls (6.95,-1.4) and (3.31,-0.3) .. (0,0) .. controls (3.31,0.3) and (6.95,1.4) .. (10.93,3.29)   ;
			
			\draw (94.33,30.83) node [anchor=north west][inner sep=0.75pt]    {$h_{q}^{( 0)}( t)$};
			\draw (153.2,30.43) node [anchor=north west][inner sep=0.75pt]   [align=left] {{\fontsize{0.67em}{0.8em}\selectfont {\fontfamily{ptm}\selectfont Newton-Nash iteration scheme}}};
			\draw (157.45,50.13) node [anchor=north west][inner sep=0.75pt]  [font=\fontsize{0.67em}{0.8em}\selectfont] [align=left] {{\fontfamily{ptm}\selectfont under the condition \eqref{hq}} };
			\draw (300,30) node [anchor=north west][inner sep=0.75pt]    {$\left( w_{q+1}^{( 0)} ,\RR_{q+1}^{( 0)}\right)$};
			\draw (95.03,121) node [anchor=north west][inner sep=0.75pt]    {$h_{q}^{(1)}( t)$};
			\draw (150.2,120.3) node [anchor=north west][inner sep=0.75pt]   [align=left] {{\fontsize{0.67em}{0.8em}\selectfont {\fontfamily{ptm}\selectfont Newton-Nash iteration scheme}}};
			\draw (160.28,140.97) node [anchor=north west][inner sep=0.75pt]  [font=\fontsize{0.67em}{0.8em}\selectfont] [align=left] {{\fontfamily{ptm}\selectfont under the condition \eqref{hq}} };
			\draw (299.67,119.13) node [anchor=north west][inner sep=0.75pt]    {$\left( w_{q+1}^{( 1)} ,\RR_{q+1}^{( 1)}\right)$};
			\draw (94.2,210.5) node [anchor=north west][inner sep=0.75pt]    {$h_{q}^{( 2)}( t)$};
			\draw (149.53,210.3) node [anchor=north west][inner sep=0.75pt]   [align=left] {{\fontsize{0.67em}{0.8em}\selectfont {\fontfamily{ptm}\selectfont Newton-Nash iteration scheme}}};
			\draw (159.12,230.3) node [anchor=north west][inner sep=0.75pt]  [font=\fontsize{0.67em}{0.8em}\selectfont] [align=left] {{\fontfamily{ptm}\selectfont under the condition \eqref{hq}} };
			\draw (299.83,209.97) node [anchor=north west][inner sep=0.75pt]    {$\left( w_{q+1}^{( 2)} ,\RR_{q+1}^{( 2)}\right)$};
			\draw (75.7,310) node [anchor=north west][inner sep=0.75pt]    {$h_{q}^{(m+1)}( t)$};
			\draw (300,310.63) node [anchor=north west][inner sep=0.75pt]    {$\left( w_{q+1}^{(m+1)} ,\RR_{q+1}^{( m+1)}\right)$};
			\draw (154.53,310.3) node [anchor=north west][inner sep=0.75pt]   [align=left] {{\fontsize{0.67em}{0.8em}\selectfont {\fontfamily{ptm}\selectfont Newton-Nash iteration scheme}}};
			\draw (160.12,330.47) node [anchor=north west][inner sep=0.75pt]  [font=\fontsize{0.67em}{0.8em}\selectfont] [align=left] {{\fontfamily{ptm}\selectfont under the condition \eqref{hq}} };
			\draw (235,365.27) node [anchor=north west][inner sep=0.75pt]    {$m\rightarrow\infty$};
			\draw (160.17,395.47) node [anchor=north west][inner sep=0.75pt]    {$\left(h_{q}(t), w_{q+1} ,\RR_{q+1}\right)$};
			\draw (105.7,240.67) node [anchor=north west][inner sep=0.75pt]   [align=left] {};
			\draw (107.7,239.17) node [anchor=north west][inner sep=0.75pt]    {$ \begin{array}{l}
					\cdot\\
					\cdot\\
					\cdot
				\end{array}$};
			\draw (222.2,239.67) node [anchor=north west][inner sep=0.75pt]    {$ \begin{array}{l}
					\cdot\\
					\cdot\\
					\cdot
				\end{array}$};
			\draw (332.2,240.67) node [anchor=north west][inner sep=0.75pt]    {$ \begin{array}{l}
					\cdot\\
					\cdot\\
					\cdot
				\end{array}$};
			\draw (198,76.71) node [anchor=north west][inner sep=0.75pt]   [align=left] {\eqref{energy-gap}};
			\draw (196,165.71) node [anchor=north west][inner sep=0.75pt]   [align=left] {\eqref{energy-gap}};

		\end{tikzpicture}
		\caption{Multiple iteration scheme}\label{MIS1}
	\end{center}
\end{figure}

\noindent{\bf Organization of this paper.}
	The remainder of this paper is organized as follows.
	Section~\ref{proof of negative} is devoted to the proof of
	Theorem~\ref{main0}. We first formulate the main iteration propositions
	and introduce the parameters used in the convex integration scheme.
	After mollifying the approximate solution, we glue exact Euler flows
	on short time intervals and obtain a new Reynolds stress with suitable
	temporal localization. We then introduce the staircase traveling waves
	and establish their estimates and exact pointwise non-interaction
	properties. These waves are incorporated into a multiple iteration
	scheme that combines the Newton--Nash iteration with a Picard-type
	iteration for the energy corrector. The resulting estimates for the
	velocity perturbation, the Reynolds stress, and the kinetic energy
	complete the iteration and yield a weak solution with the prescribed
	energy profile.
	
	Section~\ref{Proof of Theorem 2} proves Theorem~\ref{main1}. Starting
	from a smooth approximation of an arbitrary divergence-free initial
	datum, we apply the construction of Theorem~\ref{main0} to two strictly
	decreasing energy profiles that agree near the initial time and differ
	afterwards. This produces two distinct dissipative weak solutions with
	the same initial datum and gives the required density of wild initial
	data in the $C^{\gamma'}$ topology.
	
	The first appendix proves the
	\texorpdfstring{$L_t^\infty$}{L-infinity} contraction estimate for the
	Picard map and the corresponding stability estimates through the
	successive Newton levels. The second appendix records the H\"older and
	Besov notation and collects the mollification estimates, the
	anti-divergence operator, the geometric lemma, and the other analytic
	tools used throughout the construction.
	\vspace{5mm}
	
	\noindent{\bf Notation.} Throughout the paper, $x\lesssim y$ means that $x\le Cy$ for a universal constant whose value may change from line to line. The notation $\lesssim_N$ indicates that the implicit constant may depend on $N$. We write $x\ll y$ and $x\gg y$ to mean that $x$ is sufficiently smaller or sufficiently larger than $y$, respectively. The projection $\mathbb P_{\ne0}$ onto the nonzero Fourier modes
is defined by
$\mathbb P_{\ne0}f:=f-\fint_{\mathbb T^2}f(x)\,\dd x$.

\section{The construction of dissipative weak solutions}\label{proof of negative}

In this section, we prove Theorem~\ref{main0}
for an arbitrary fixed exponent $0<\gamma<\frac13$. The proof rests on the following two iteration results,
Propositions~\ref{p:main-prop} and~\ref{p:main-prop2}.

\subsection{Iteration proposition}\label{Iteration proposition}
To prove Theorem~\ref{main0} by convex integration method, we consider the following relaxation of \eqref{e:E}:
\begin{equation}
	\left\{
	\begin{alignedat}{2}
		&\del_t v_q + \div (v_q \otimes v_q) + \nabla p_q = \div \RR_q, \\
		&\nabla \cdot v_q = 0,
	\end{alignedat}
	\right. \label{e:subsol-B}
\end{equation}
where
\begin{equation}
	\int_{\TTT^2} p_q \dd x = 0, \label{p_q}
\end{equation}
and the stress $\RR_q$ is assumed to be a symmetric, trace-free $2 \times 2$ matrix.

We first specify the parameters used in the iteration. Fix $0 < \beta < \frac{1}{3}$ and set
$$
b_0 = \min\big\{1+\tfrac{1-3\beta}{12\beta}, 2\big\},
$$
and
$$
L \coloneq \max\Big\{ \big\lceil\frac{-10}{2\beta b^2-b(3\beta+\frac{1}{3})+ \beta+\frac{1}{3}}\big\rceil,\big\lceil\frac{100}{(\frac{1}{3}-\beta)(b-1) }\big\rceil,100\Big\} \in \ZZ_{+},
$$
where $b \in (1, b_0)$ and $\lceil x \rceil$ denotes the ceiling function. We further set
\begin{align}
	\epsilon_{0}:= \frac{-(2\beta b^2-b(3\beta+\frac{1}{3})+ \beta+\frac{1}{3})}{20L}, \label{e:params0}
\end{align}
and 
\begin{align}\label{N0}
	\alpha=\frac{\epsilon_{0}}{100},
	\qquad
	N_0\coloneq
		\left\lceil
		\frac{2+2\beta b^2+5\alpha b}{\epsilon_{0}}
		\right\rceil+1\in\ZZ_{+}.
\end{align}
In particular,
\begin{equation*}
	N_0\epsilon_{0}>2+2\beta b^2+5\alpha b.
\end{equation*}
Let $a$ be a sufficiently large positive parameter, depending on $M, \alpha, L, N_0, e(t)$ and $T$, where $M$ is a large constant.

Since the sets in Lemma~\ref{first S} are finite and consist of rational directions, fix \(\lambda_\ast\in\mathbb N\) such that \(\lambda_\ast k\in\mathbb Z^2\) for all \(k\in\Lambda_1\cup\Lambda_2\).
Define
\begin{align}
	\lambda_q \coloneq \lambda_\ast\left\lceil a^{b^q} \right\rceil, \quad \delta_q \coloneq \lambda_q^{-2\beta},\label{lambdaq}
\end{align}
and define $\ell_q$, $\tau_q$, and $\mu_{q+1}$ by
\begin{align}
	\ell_q^{-1} &\coloneq \lambda_q  \left\lceil\lambda_q^{\epsilon_{0}}\right\rceil, \label{e:ell} \\
	\tau_q &\coloneq \frac{\lambda_q^{-5\alpha}}{\delta_q^{1/2} \lambda_q}, \quad \mu_{q+1} \coloneq \left\lceil\lambda_{q+1}^{\frac{1}{3}} \lambda_q^{2/3+5\alpha} \delta_{q+1}^{1/2}\right\rceil.
\end{align}
These choices ensure that $M \ll \ell_q^{-\frac{\alpha}{2}}$ and $\ell_q^{-1} \ll \lambda_{q+1}$. Moreover,
\begin{equation*}
	\lambda_q^2(\lambda_q\ell_q)^{N_0}
	\leq \lambda_q^{2-N_0\epsilon_0}
	\ll \delta_{q+2}\lambda_{q+1}^{-5\alpha}
\end{equation*}
provided $a$ is sufficiently large. Finally, define
$$
\delta_{q+1,n} \coloneq (\tau_q \mu_{q+1} \ell_q^{10\alpha+L\epsilon_{0}})^{-n} \delta_{q+1},
$$
for $0 \leq n \leq L$. For all $q \geq 1$, the following estimate holds once $a$ is sufficiently large:
\begin{equation}
	\lambda_{q}^{2}(\lambda_q\ell_q)^{N_0} + \delta_{q+1,L} + \ell_q^{-2}\left\lceil\ell_q^{-\epsilon_0}\right\rceil^{L+1} \lambda_q^{-2} \Big(\frac{\lambda_q}{\lambda_{q+1}}\Big)^{\frac{1}{3}} \Big(\frac{\delta_{q+1}}{\delta_q}\Big)^{1/2} \delta_{q+1} \leq \delta_{q+2} \lambda_{q+1}^{-5\alpha}. \label{lambdaN} 
\end{equation}

    We now state the main iteration proposition for the system \eqref{e:subsol-B}.

\begin{prop}[Main iteration proposition]
	\label{p:main-prop}
	Let $b \in (1, b_0)$ and $M$ be a large constant. Suppose that $(v_q, p_q, \RR_q)$ satisfies \eqref{e:subsol-B}, \eqref{p_q}, and
	\begin{align}
		&\|v_q\|_0 \leq M\Big(1 + \sum_{i=1}^q \delta_i^{1/2}\Big), \label{e:vq-C0} \\
		&\|v_q\|_N \leq M \delta_q^{1/2} \lambda_q^N, \quad 1 \leq N \leq N_0+2, \label{e:vq-C1} \\
		&\|\RR_q\|_0 \leq \delta_{q+1} \lambda_q^{-4\alpha}, \label{e:RR_q-C0}\\
        &\|\partial_t\RR_q\|_0 \leq \lambda_q^{5}\delta_{q+1} , \label{t:RR_q-C0}
	\end{align}
	and
 	\begin{align}
		10\delta_{q+1} \leq e(t) - \|v_q(t)\|_{L^2(\TTT^2)}^2 \leq 20\delta_{q+1}.\label{e:initial}
	\end{align}
	Then there exists a solution $(v_{q+1}, p_{q+1}, \RR_{q+1})$ satisfying \eqref{e:subsol-B}, \eqref{p_q}, and \eqref{e:vq-C0}--\eqref{e:initial} with $q$ replaced by $q+1$. In addition,
	\begin{align}
		\|v_{q+1} - v_q\|_0 + \frac{1}{\lambda_{q+1}} \|v_{q+1} - v_q\|_1 \leq M \delta_{q+1}^{1/2}. \label{e:velocity-diff}
	\end{align}
\end{prop}
The following proposition, which is mainly motivated by \cite{MNY}(Proposition 3.3), will be used solely to establish the non-uniqueness result.
\begin{prop}
	\label{p:main-prop2}
	Let $T > 0$. Suppose that $(v_q, p_q, \RR_q)$ and $(\widetilde{v}_q, \widetilde{p}_q, \widetilde{\RR}_q)$ are two solutions of the relaxed system \eqref{e:subsol-B}, and that $e(t), \widetilde{e}(t)$ satisfy the hypotheses of Proposition~\ref{p:main-prop}, with $e(t) = \widetilde{e}(t)$ on $[0, 10\tau_1]$. If $(v_q,p_q,\RR_q)=(\widetilde v_q,\widetilde p_q,\widetilde{\RR}_q)$ on $[0,8\tau_q]$, then $(v_{q+1},p_{q+1},\RR_{q+1})=(\widetilde v_{q+1},\widetilde p_{q+1},\widetilde{\RR}_{q+1})$ on $[0,8\tau_{q+1}]$.
\end{prop}

\noindent \textbf{Propositions~\ref{p:main-prop} and~\ref{p:main-prop2} imply Theorem~\ref{main0}.}
Without loss of generality, we set $T=1$. Choose $M$ sufficiently large
depending on $\|e\|_{C^0([0,1])}$, and then choose $a$ sufficiently large
that $15\delta_2<\frac12\min_{[0,1]}e$. Set
\[
K:=\left\lceil\lambda_1^{2/3}\right\rceil,
\qquad A(t):=\sqrt{2\bigl(e(t)-15\delta_2\bigr)},
\]
and initialize the iteration with
\begin{align*}
	v_1(t,x)&=\bigl(A(t)\sin(2\pi Kx_2),0\bigr),\qquad p_1=0,\\
	\RR_1(t,x)&=-\frac{A'(t)}{2\pi K}\cos(2\pi Kx_2)
	\begin{pmatrix}0&1\\1&0\end{pmatrix}.
\end{align*}
This smooth triple satisfies \eqref{e:subsol-B}, \eqref{p_q}, and
\[
e(t)-\|v_1(t)\|_{L^2(\TTT^2)}^2=15\delta_2.
\]
Since $2\beta b+4\alpha<\frac23<1-\beta$,  one can easily deduce that
\begin{align*}
	\|v_1\|_N&\lesssim_{e,N}K^N
	\ll\delta_1^{1/2}\lambda_1^N,
	\qquad 1\leq N\leq N_0+2,\\
	\|\RR_1\|_0+\|\partial_t\RR_1\|_0
	&\lesssim_e K^{-1}\ll\delta_2\lambda_1^{-4\alpha}.
\end{align*}
Together with $\|v_1\|_0\leq\sqrt{2\|e\|_{C^0([0,1])}}\leq M$,
these bounds verify \eqref{e:vq-C0}--\eqref{e:initial} for $q=1$.
 Then, we obtain a sequence of solutions $\{(v_q, p_q, \RR_q)\}_{q \geq 1}$ to system \eqref{e:subsol-B} satisfying \eqref{p_q} and \eqref{e:vq-C0}--\eqref{e:initial}. According to \eqref{p_q}, \eqref{e:vq-C0}, \eqref{e:RR_q-C0}, and \eqref{e:velocity-diff}, the sequence $v_q$ converges uniformly to some continuous function $v$, and $p_q$ converges to $p$ in $L^m$ for some $m \in (1, \infty)$. The pair $(v, p)$ satisfies the Euler system \eqref{e:E} in a weak sense since $\|\RR_q\|_0 \to 0$ as $q \to \infty$.  By interpolation in \eqref{e:velocity-diff}, the time-regularity
argument in \cite[Section~2.2]{BDSV}, and Schauder estimates, we obtain
$v\in C^\gamma(\TTT^2\times[0,1])$ and the continuous normalized pressure
$p=-\Delta^{-1}\operatorname{div}\operatorname{div}(v\otimes v)
\in C([0,1];C^\eta(\TTT^2))$
for every $0<\eta<\gamma<\beta$. Furthermore, let two smooth positive functions $e(t), \widetilde{e}(t)$ satisfy \eqref{e:initial} with $e(t) = \widetilde{e}(t)$ in $[0, 10\tau_1]$ and $e(t) \neq \widetilde{e}(t)$ in $[20\tau_1, 1]$. Then, by Propositions~\ref{p:main-prop} and~\ref{p:main-prop2}, there exist two weak solutions such that
$$
v, \widetilde{v} \in C^{\gamma}(\TTT^2 \times [0, 1]), \quad \|v\|_{L^2}^2 = e(t),  \quad \|\widetilde{v}\|_{L^2}^2 = \widetilde{e}(t).
$$
Hence, non-uniqueness is evident, and we complete the proof of Theorem~\ref{main0}.

The remainder of this section is devoted to the proofs of Propositions~\ref{p:main-prop} and~\ref{p:main-prop2}. The argument consists of a mollification step, gluing step, and the construction of the perturbations.
\subsection{Mollification step and gluing procedure}\label{3.2}
We first mollify the subsolution $(v_q,p_q,\RR_q)$ to obtain $(v_{\ell_q},p_{\ell_q},\RR_{\ell_q})$, and then use the mollified triple to construct the glued Euler flow $(\vv_q,\pp_q,\RRR_q)$ at level $q$.

\subsubsection{Mollification step}
The mollification step prevents a loss of derivatives in the convex integration scheme. With $\ell_q$ as in \eqref{e:ell}, we define $(v_{\ell_q},p_{\ell_q},\RR_{\ell_q})$ using the spatial mollifier $\psi_{\ell_q}$ introduced in \eqref{e:defn-mollifier-x}:
\begin{align*}
	& {v_{\ell_q} \coloneq v_q * \psi_{\ell_q}, \quad p_{\ell_q} \coloneq p_q * \psi_{\ell_q} + \tfrac12(|v_q|^2 * \psi_{\ell_q}) - \tfrac12 |v_{\ell_q}|^2,} \\
	& \RR_{\ell_q} \coloneq \RR_q * \psi_{\ell_q} - (v_q \ootimes v_q) * \psi_{\ell_q} + v_{\ell_q} \ootimes v_{\ell_q},
\end{align*}
then
\begin{equation}
	\left\{
	\begin{alignedat}{2}
		&\del_t v_{\ell_q} + \div (v_{\ell_q} \otimes v_{\ell_q}) + \nabla p_{\ell_q} = \div \RR_{\ell_q}, \\
		&\nabla \cdot v_{\ell_q} = 0.
	\end{alignedat}
	\right. \label{e:mollified-euler}
\end{equation}

The standard arguments of \cite{BDSV,MNY}, together with the assumption $\lambda_q^{-3\alpha} \leq \ell_q^{\frac{9}{4}\alpha}$, yield the following estimates.

\begin{prop}[Estimates for mollified functions]\label{p:estimates-for-mollified}
	Let $\alpha$ satisfy \eqref{e:params0}. Then the following estimates hold uniformly for $t\in[0,T]$:
	\begin{align}
		&\|v_{\ell_q}-v_q\|_2
			\lesssim \ell_q^{N_0}\|v_q\|_{N_0+2}
			\lesssim M\delta_q^{1/2}\lambda_q^2(\lambda_q\ell_q)^{N_0}
			\lesssim \lambda_q^2(\lambda_q\ell_q)^{N_0}, \label{e:v_ell-vq-high}\\
		&\|v_{\ell_q}\|_{N+1} \lesssim \delta_q^{1/2} \lambda_q \ell_q^{-N}, && 0\leq N\leq N_0+3, \label{e:v_ell-CN+1} \\
		&\|\RR_{\ell_q}\|_{N+\alpha} \lesssim \delta_{q+1} \ell_q^{-N+\frac{5}{4}\alpha}, && 0\leq N\leq N_0+3, \label{e:R_ell} \\
        &\|\partial_t\RR_{\ell_q}\|_N \lesssim \lambda_q^{5}\ell_q^{-N }\delta_{q+1} , && 0\leq N\leq N_0+3,\\
		&9\delta_{q+1} \leq e(t) - \|v_{\ell_q}\|_{L^2(\TTT^2)}^2 \leq 21\delta_{q+1}. \label{e:initial2}
	\end{align}
\end{prop}

\subsubsection{Classical exact  flows}
Set
\[
t_0 = 0, \quad t_i \coloneq t_0 + i\tau_q, \quad i \in \mathbb{N},
\]
and
\[
i_{\max} = \sup\{i \geq 1 \mid \tau_q \leq T - t_i < 2\tau_q\} + 1.
\]

For each $0\leq i\leq i_{\max}$, let $(\vex_i,\pex_i)$ be the unique smooth
solution of the two-dimensional Euler equations
\begin{equation}
	\left\{
	\begin{alignedat}{2}
		&\del_t\vex_i+(\vex_i\cdot\nabla)\vex_i+\nabla\pex_i=0,\\
		&\nabla\cdot\vex_i=0,\\
		&\vex_i|_{t=t_i}=v_{\ell_q}(\cdot,t_i)
	\end{alignedat}
	\right.
	\label{e:exact-B}
\end{equation}
on its maximal interval of
existence. 

Subtracting equation~\eqref{e:mollified-euler} from equation~\eqref{e:exact-B}, we obtain
\[
\begin{aligned}
\del_t(\vex_i-v_{\ell_q})
+(\vex_i-v_{\ell_q})\cdot\nabla\vex_i
+\nabla(\pex_i-p_{\ell_q})
=-\operatorname{div}\RR_{\ell_q}, \ (\vex_i-v_{\ell_q})\big|_{t=t_i}=0.
\end{aligned}
\]
Moreover, Proposition~\ref{p:estimates-for-mollified} gives
$
\tau_q\|v_{\ell_q}\|_{1+\alpha}
\lesssim\lambda_q^{-5\alpha}\ell_q^{-\alpha}\ll1
$,
and the stability argument in \cite{MNY1} therefore ensures that
$(\vex_i,\pex_i)$ is well-defined on
$[t_i-\tau_q,t_i+\tau_q]$ and yields the estimates
in the following proposition.

\begin{prop}[Local well-posedness]\label{p:stability}
	For each $i\geq0$, the system \eqref{e:exact-B} has a unique smooth solution $(\vex_i,\pex_i)$ on $|t-t_i|\leq\tau_q$. Moreover, for $0\leq N\leq N_0+3$,
	\begin{align}
		\|\vex_i\|_{N+1+\alpha} &\lesssim \ell^{-({N+\alpha})}_q\lambda_{q}\delta^{1/2}_q, \label{e:vex-bound} \\
		\|\vex_i - v_{\ell_q}\|_{N+\alpha} &\lesssim \tau_q \delta_{q+1} \ell_q^{-N-1+\frac{5}{4}\alpha}, \label{e:stability-v} \\
		\|\nabla \pex_i - \nabla p_{\ell_q}\|_{N+\alpha} &\lesssim \delta_{q+1} \ell_q^{-N-1+\frac{5}{4}\alpha}, \label{e:stability-p} \\
		\|\matd{v_{\ell_q}} (\vex_i - v_{\ell_q})\|_{N+\alpha} &\lesssim \delta_{q+1} \ell_q^{-N-1+\frac{5}{4}\alpha}, \label{e:stability-matd} \\
		\|\mathcal{R}(\vex_i - v_{\ell_q})\|_{\alpha} &\lesssim \tau_q \delta_{q+1} \ell_q^{\frac{5}{4}\alpha}, \label{e:stability-v-1}
	\end{align}
	where the inverse-divergence operator $\mathcal{R}$ is defined in Definition~\ref{def.R}.
\end{prop}

\subsubsection{Gluing procedure}
We construct the glued velocity and pressure $(\vv_q,\ppp_q)$ from the exact flows $(\vex_i,\pex_i)$ by setting
\begin{align*}
	\vv_q(x,t) \coloneq \sum_{i=0}^{i_\text{max}} \eta_i(t) \vex_i(x,t),\quad
	\ppp_q(x,t) \coloneq \sum_{i=0}^{i_\text{max}} \eta_i(t) \pex_i(x,t),
\end{align*}
where $\{\eta_i\}_{i=0}^{i_\text{max}}$ is a partition of unity on $[0,T]$ satisfying, for all $N\geq0$,
\begin{align*}
	&\text{supp}_t \eta_0 \subset[-1, \tfrac{2\tau_q}{3}], \quad\quad\ \ \quad\quad\quad\eta_0|_{J_0} = 1, \quad\quad \| \del_t^N \eta_0 \|_0 \lesssim \tau_q^{-N}, \\
	&\text{supp}_t \eta_i = I_{i-1} \cup J_i \cup I_i,  \quad \quad\quad\quad  \eta_i|_{J_i} = 1, \quad\quad \| \del_t^N \eta_i \|_0 \lesssim \tau_q^{-N}, \quad 1 \leq i \leq i_\text{max}-1, \\
	&\text{supp}_t \eta_{i_\text{max}} \subset [t_{i_\text{max}-1} + \tfrac{\tau_q}{3}, T], \quad \eta_{i_\text{max}}|_{J_{i_\text{max}}} = 1, \quad \| \del_t^N \eta_{i_\text{max}} \|_0 \lesssim \tau_q^{-N}.
\end{align*}
The intervals $I_i$ and $J_i$ are defined by
\begin{align*}
	&I_i \coloneq [t_i + \tfrac{\tau_q}{3}, t_i + \tfrac{2\tau_q}{3}], \quad 0 \leq i \leq i_\text{max}-1, \\
	&J_i \coloneq (t_i - \tfrac{\tau_q}{3}, t_i + \tfrac{\tau_q}{3}), \quad 1 \leq i \leq i_\text{max}-1, \\
	&J_0 \coloneq [0, \tfrac{\tau_q}{3}), \quad J_{i_\text{max}} \coloneq (t_{i_\text{max}-1} + \tfrac{2\tau_q}{3}, T].
\end{align*}

A direct computation shows that $(\vv_q,\pp_q)$ satisfies the following system on $[0,T]$:
\begin{equation}
	\left\{
	\begin{alignedat}{2}
		&\del_t \vv_q + \div (\vv_q \otimes \vv_q) + \nabla \pp_q = \div \RRR_q, \\
		&\nabla \cdot \vv_q = 0,
	\end{alignedat}
	\right. \label{e:glu-B}
\end{equation}
where
\begin{align}
	\RRR_q &\coloneq \sum_{i=0}^{i_\text{max}-1} \del_t \eta_i \mathcal{R}(\vex_i - \vex_{i+1}) - \sum_{i=0}^{i_\text{max}-1} \eta_i(1-\eta_i)(\vex_i - \vex_{i+1}) \ootimes (\vex_i - \vex_{i+1}), \label{e:def-RRR_q} \\
	\pp_q &\coloneq \ppp_q - \int_{\TTT^2} \ppp_q \dd x + {\frac12\sum_{i=0}^{i_\text{max}-1} \eta_i(1-\eta_i)\bigg(|\vex_i - \vex_{i+1}|^2 - \int_{\mathbb{T}^2} |\vex_i - \vex_{i+1}|^2 \dd x\bigg)}. \label{ppp}
\end{align}
In particular, $\vv_q$ is an exact solution of the two-dimensional Euler system \eqref{e:E} on each $J_i$, and $\textup{supp}_t \RRR_q \subset \bigcup\limits_{i=0}^{i_\text{max}-1} I_i$.

We have the following estimates for $(\vv_q,\RRR_q)$.
\begin{prop}[Estimates for $(\vv_q, \RRR_q)$]\label{est-gluing}
	For $0 \leq N \leq N_0+3$,
	\begin{align}
		&\| \vv_q - v_{\ell_q} \|_\alpha \leq \delta_{q+1}^{1/2} \ell_q^\alpha, \label{e:stability-vv_q} \\
		&\| \vv_q \|_0 \leq M\Big(1+\sum_{i=1}^q \delta_i^{1/2}\Big) + \delta_{q+1}^{1/2} \ell_q^\alpha, \label{e:stability-vv_q-N} \\
		&\| \vv_q \|_{N+1} \lesssim \delta_q^{1/2} \lambda_q \ell_q^{-N}, \label{e:vv_q-bound} \\
		&\| \RRR_q \|_{N+\alpha} \lesssim \delta_{q+1} \ell_q^{-N+\frac{5}{4}\alpha}, \label{e:RRR_q-N+alpha-bd} \\
		&\| \matd{\vv_q} \RRR_q \|_{N+\alpha} \lesssim \tau^{-1}_q \delta_{q+1}\ell_q^{-N+\frac{5}{4}\alpha}, \label{e:matd-RRR_q} \\
		&8\delta_{q+1} \leq e(t) - \|\vv_q\|_{L^2(\TTT^2)}^2 \leq 23\delta_{q+1}, \label{e:initial1}\\
	& {|\partial_t(e(t) - \|\vv_q\|_{L^2(\TTT^2)}^2) |\leq \tau^{-1}_q\delta_{q+1}\ell^{\alpha}_q,} \label{e:glued-energy-derivative}\\
     &\|\partial_t\RRR_{q}\|_N\lesssim \ell_q^{-N -1}\delta_{q+1} , \label{t:matd-RRR_q} \\
        & \|\partial_t\matd{\vv_q}\RRR_{ q}\|_N \lesssim \lambda_q^{3}\ell_q^{-N }\delta_{q+1} , \label{tt:matd-RRR_q} \\
        &|\partial^2_t(e(t) - \|\vv_q\|_{L^2(\TTT^2)}^2) |\leq \tau^{-1}_q\ell^{-1-\alpha}_{q}\delta_{q+1} . \label{2te:matd-RRR_q}  
	\end{align}
\end{prop}
The proof of Proposition \ref{est-gluing} could be found in   \cite{MNY1}, while \eqref{tt:matd-RRR_q} and \eqref{2te:matd-RRR_q} are valid by using \eqref{e:exact-B} and \eqref{e:glu-B}.

Let ${\Phi}_i$ denote the inverse flow map, defined as the solution of
\begin{equation}
	\left\{
	\begin{alignedat}{2}
		&\del_t {\Phi}_i + \vv_q\cdot \nabla {\Phi}_i = 0, \\
		&{\Phi}_i(t_i,x) = x.
	\end{alignedat}
	\right. \label{phi}
\end{equation}

The following estimates will be used repeatedly.

\begin{prop}[Estimates for $\Phi_i$, \cite{BDSV}]\label{p:estimates-for-inverse-flow-map}
	For $a\gg1$, $0 \leq N \leq N_0+3$, and $t\in[t_i-\tau_q,t_i+2\tau_q]\cap[0,T]$, one has
	\begin{align}
		\|\nabla \Phi_i - {\rm Id}_{2 \times 2}\|_0 &\leq \frac{1}{10}, \label{e:nabla-phi-i-minus-I3x3} \\
		\|(\nabla \Phi_i)^{-1}\|_N + \|\nabla \Phi_i\|_N &\lesssim \ell_q^{-N}, \label{e:nabla-phi-i-CN} \\
		\|\matd{\vv_q} \nabla \Phi_i\|_N &\lesssim \delta_q^{1/2} \lambda_q \ell_q^{-N}. \label{e:nabla-phi-i-matd}
	\end{align}
\end{prop}

This completes the passage from $v_q$ to $\vv_q$. We next construct a perturbation $w_{q+1}$ of $\vv_q$ and the corresponding Reynolds stress $\RR_{q+1}$ at level $q+1$. More precisely,
\[\vv_q \mapsto v_{q+1} := \vv_q + w_{q+1}.\]
\subsection{Staircase traveling waves }\label{Traveling waves}
We now introduce staircase traveling waves that modulate the kinetic energy on arbitrary time intervals.

Choose temporal cutoffs $\{\bar\chi_i(t)\}_i$ with the following properties:
\begin{enumerate}
	\item $ \bar{\chi}_i \in C^\infty_c(  [t_i+\tfrac{1}{12}\tau_q, t_{i+1} - \tfrac{1}{12}\tau_q] ; [0,1])$ and $ \bar{\chi}_{i }(t) \equiv 1$ for $t\in I_i$. Moreover,
	\begin{align}
		\|\del^n_t  \bar{\chi}_{i}\|_{L^\infty_tC^\gamma_x    }  \lesssim_{n,\gamma} \tau_q^{-n},~~~n,\gamma\ge 0. \notag
	\end{align}
	\item $\supp_{t}~\bar{\chi}_{i} \cap \supp_{t} ~\bar{\chi}_{j}=\emptyset$ whenever $i\neq j$.
\end{enumerate}
To correct the energy gap, we also choose space--time cutoffs
	$\widetilde\chi_i\in C^\infty(\mathbb T^2\times[0,T];[0,1])$,
$0\leq i\leq i_{\max}-1$, as in \cite{BDSV}, with the following properties:
\begin{enumerate}
	\item 
For $0\leq i\leq i_{\max}-2$, we set
\[
\operatorname{supp}_t\widetilde\chi_i
\subset
[t_i-\tfrac{1}{48}\tau_q,\,
 t_{i+1}+\tfrac{1}{48}\tau_q]\cap[0,T],
\]
and $\widetilde\chi_i(\cdot,t)\equiv1$ for $t\in I_i$.
For the last index $i=i_{\max}-1$, we choose
\[
\operatorname{supp}_t\widetilde\chi_{i_{\max}-1}
\subset
[t_{i_{\max}-1}-\tfrac{1}{48}\tau_q,\,T],
\]
and $\widetilde\chi_{i_{\max}-1}(\cdot,t)\equiv1$ for
$t\in[t_{i_{\max}-1}+\tfrac13\tau_q,\,T]$.
    Moreover,
	\begin{align}
		\|\del^n_t  \widetilde{\chi}_{i}\|_{L^\infty_tC^\gamma_x    }  \lesssim_{n,\gamma} \tau_q^{-n},~~~n,\gamma\ge 0. \notag
	\end{align}
	\item $\supp~\widetilde{\chi}_{i} \cap \supp~\widetilde{\chi}_{j}=\emptyset$~and~ $\supp~\bar{\chi}_{i} \cap \supp~\widetilde{\chi}_{j}=\emptyset$ whenever $i\neq j$.
	\item  For all $t\in[0,T]$,
	$ \frac{1}{4} \le \sum_{i=0}^{i_{max}-1} \int_{\mathbb T^2}  \widetilde{\chi}_{i}^2(x,t) \dd x \le 1.$
\end{enumerate}

\begin{prop}\label{time-g}
	Let $\Lambda \subset \mathbb{S} \cap \mathbb{Q}^2$ be the set supplied by Lemma~\ref{first S}, and fix $\xi\in \mathbb{S} \cap \mathbb{Q}^2$ with $\xi\notin\Lambda$. Choose a nonzero $h\in C^{\infty}(\mathbb{T})$ with sufficiently small support in {$[\frac{1}{100},\frac{2}{100}]$}. After normalization and suitable translations of $h$, there exist smooth $1$-periodic functions $h_{\eta } : \mathbb{T} \to \mathbb{R}$, indexed by $\eta\in\Lambda\cup\{\xi\}$, such that
	$$
	\int_0^1 h^2_{\eta }(y) dy = 1 \quad \textup{and} \quad {\textup{supp}\,h_{\eta } \cap \textup{supp}\,h_{\eta'  } = \emptyset},
	$$
	whenever $\eta\neq\eta'$.
\end{prop}

For $0\leq n\leq L-1$ and $k\in\Lambda$, define the $1$-periodic function
\[
g_{k,n}(y):=
h_k\!\left(\left\lceil\ell_q^{-\epsilon_0}\right\rceil^n y\right)
\prod_{0\leq \gamma\leq n-1}
h_\xi\!\left(\left\lceil\ell_q^{-\epsilon_0}\right\rceil^\gamma y\right),
\qquad y\in\mathbb T,
\]
where
$g_{k,0}(y)=h_k(y)$. Given $i$, let
\[
y=\mu_{q+1}t-\ell_q^{-1}\left\lceil\ell_q^{-\epsilon_0}\right\rceil\Phi_i^1(t,x),
\]
and define the staircase traveling waves by
\[
g_{k,i,n}(t,x):=g_{k,n}(y).
\]
	For $0\leq n,n'\leq L-1$ and $k,k'\in\Lambda$, whenever
		$(k,i,n)\neq(k',i',n')$, by Proposition \ref{time-g}   we deduce
	\begin{align}
		\bar{\chi}_i g_{k,i,n}\,\bar{\chi}_{i'}g_{k',i',n'}=0,\qquad
		\widetilde{\chi}_i g_{k,i,n}\,\widetilde{\chi}_{i'}g_{k',i',n'}=0,\qquad
		\bar{\chi}_i g_{k,i,n}\,\widetilde{\chi}_{i'}g_{k',i',n'}=0. \label{bujiao}
	\end{align}
Indeed, for $i\neq i'$, the cutoff supports are disjoint. For a fixed $i$, if $n\neq n'$, the staircase structure illustrated in
	Figure~\ref{fig:staircase-profiles} ensures that $g_{k,n}$ and
	$g_{k',n'}$ have disjoint supports.
	If $n=n'$ but $k\neq k'$, the same conclusion follows from
	the terminal factors $h_k$ and $h_{k'}$, which have disjoint
	supports and are evaluated at the same argument.
	
	Next, we give the decomposition of $g_{k,i,n}^2(t,x)$. Since $h_\eta^2$ is smooth and $1$-periodic, its Fourier series
	\[
	h_\eta^2(y)=\sum_{j\in\mathbb Z}b_j^{(\eta)}e^{2\pi i j y},
	\qquad b_0^{(\eta)}=\int_{\mathbb T}h_\eta^2(y)\,\dd y=1,
	\qquad \eta\in\Lambda\cup\{\xi\},
	\]
	converges absolutely. Hence, by the definition of $g_{k,n}$, for
		$1\leq n\leq L-1$,
	\begin{align*}
		g_{k,n}^2(y)
		={}&h_k^2\!\left(\left\lceil\ell_q^{-\epsilon_0}\right\rceil^n y\right)
		\prod_{\gamma=0}^{n-1}
		h_\xi^2\!\left(\left\lceil\ell_q^{-\epsilon_0}\right\rceil^\gamma y\right)\\
		={}&\left(\sum_{j_n\in\mathbb Z}b_{j_n}^{(k)}
		e^{2\pi i j_n\left\lceil\ell_q^{-\epsilon_0}\right\rceil^n y}\right)
		\prod_{\gamma=0}^{n-1}
		\left(\sum_{j_\gamma\in\mathbb Z}b_{j_\gamma}^{(\xi)}
		e^{2\pi i j_\gamma\left\lceil\ell_q^{-\epsilon_0}\right\rceil^\gamma y}\right)\\
		= &\sum_{j_0,\ldots,j_n\in\mathbb Z}\Big(
		b_{j_n}^{(k)}\prod_{\gamma=0}^{n-1}b_{j_\gamma}^{(\xi)}\Big)
		e^{2\pi i\left(j_0+j_1\left\lceil\ell_q^{-\epsilon_0}\right\rceil
			+\cdots+j_n\left\lceil\ell_q^{-\epsilon_0}\right\rceil^n\right)y}.
	\end{align*}
	Setting
	\[
	{\bf J}_n:=
	j_0+j_1\left\lceil\ell_q^{-\epsilon_0}\right\rceil+\cdots+
	j_n\left\lceil\ell_q^{-\epsilon_0}\right\rceil^n ,
	\]
	we obtain
	\begin{align}\label{fourier-g}
		g_{k,i,n}^2(t,x)
		={}&\sum_{j_0,\ldots,j_n\in\mathbb Z}
		b_{j_n}^{(k)}\prod_{\gamma=0}^{n-1}b_{j_\gamma}^{(\xi)}
		e^{2\pi i{\bf J}_n
			\left(\mu_{q+1}t-\ell_q^{-1}
			\left\lceil\ell_q^{-\epsilon_0}\right\rceil\Phi_i^1\right)}
		\notag\\
		={}&{\sum_{\substack{j_0,\ldots,j_n\in\mathbb Z\\
					{\bf J}_n=0}}}
		b_{j_n}^{(k)}\prod_{\gamma=0}^{n-1}b_{j_\gamma}^{(\xi)}
		+\sum_{\substack{j_0,\ldots,j_n\in\mathbb Z\\
				{\bf J}_n\neq0}}
		b_{j_n}^{(k)}\prod_{\gamma=0}^{n-1}b_{j_\gamma}^{(\xi)}
		e^{2\pi i{\bf J}_n
			\left(\mu_{q+1}t-\ell_q^{-1}
			\left\lceil\ell_q^{-\epsilon_0}\right\rceil\Phi_i^1\right)}
		\notag\\
		:= {}&\int_{\mathbb T}g_{k,n}^2(y)\,\dd y-f_{k,i,n}(t,x),
	\end{align}
	where
	\begin{align}
		\label{fkin-traveling}
		f_{k,i,n}(t,x)
		={}&
		-\sum_{\substack{j_0,\ldots,j_n\in\mathbb Z\\
				{\bf J}_n\neq0}}
		b_{j_n}^{(k)}
		\prod_{\gamma=0}^{n-1}b_{j_\gamma}^{(\xi)}
		e^{2\pi i{\bf J}_n
			\left(\mu_{q+1}t-\ell_q^{-1}
			\left\lceil\ell_q^{-\epsilon_0}\right\rceil\Phi_i^1\right)}
		\notag\\
		={}&
		-\sum_{\substack{j_0,\ldots,j_n\in\mathbb Z\\
				{\bf J}_n\neq0}}
		\frac{
			b_{j_n}^{(k)}
			\prod_{\gamma=0}^{n-1}b_{j_\gamma}^{(\xi)}
		}{
			2\pi i{\bf J}_n\mu_{q+1}}
		e^{-2\pi i{\bf J}_n\ell_q^{-1}
			\left\lceil\ell_q^{-\epsilon_0}\right\rceil\Phi_i^1}
		\partial_t e^{2\pi i{\bf J}_n\mu_{q+1}t}.
	\end{align}

	 We claim that $\int_{\mathbb T}g_{k,n}^2(y)\,\dd y\approx 1$. Indeed, for $1\leq n\leq L-1$, using
	$h_k^2=1+\mathbb P_{\ne0}h_k^2$, we have
	\begin{align*}
		\int_{\mathbb T}g_{k,n}^2(y)\,\dd y
		={}&\int_{\mathbb T}
		\prod_{\gamma=0}^{n-1}
		h_\xi^2\!\left(\left\lceil\ell_q^{-\epsilon_0}\right\rceil^\gamma y\right)\,\dd y+\int_{\mathbb T}
		\big(\mathbb P_{\ne0}h_k^2\big)\!\left(\left\lceil\ell_q^{-\epsilon_0}\right\rceil^n y\right)
		\prod_{\gamma=0}^{n-1}
		h_\xi^2\!\left(\left\lceil\ell_q^{-\epsilon_0}\right\rceil^\gamma y\right)\,\dd y\\
		=:{}&c_{n,\xi}
		+\sum_{j\neq0}b_j^{(k)}\int_{\mathbb T}
		e^{2\pi i j\left\lceil\ell_q^{-\epsilon_0}\right\rceil^n y}
		\prod_{\gamma=0}^{n-1}
		h_\xi^2\!\left(\left\lceil\ell_q^{-\epsilon_0}\right\rceil^\gamma y\right)\,\dd y,
	\end{align*}
	where $c_{n,\xi}=\int_{\mathbb T}\prod_{\gamma=0}^{n-1}
		h_\xi^2\!\left(\left\lceil\ell_q^{-\epsilon_0}\right\rceil^\gamma y\right)\,\dd y$. Integrating by parts $N_0$ times and using the rapid decay
		of $b_j^{(k)}$  yields
	\begin{align}
		\left|\int_{\mathbb T}g_{k,n}^2(y)\,\dd y-c_{n,\xi}\right|
		&\lesssim
		\ell_q^{nN_0\epsilon_0}
		\left\|\partial_y^{N_0}\prod_{\gamma=0}^{n-1}
		h_\xi^2\!\left(\left\lceil\ell_q^{-\epsilon_0}\right\rceil^\gamma y\right)
		\right\|_{L^1(\mathbb T)}\notag\\
		&\lesssim
		\ell_q^{nN_0\epsilon_0}
		\ell_q^{-(n-1)N_0\epsilon_0}
		\lesssim
		\ell_q^{N_0\epsilon_0}.
		\label{zero-mode-k}
	\end{align}
	For $2\leq n\leq L-1$, by the definition of $c_{n,\xi}$ and the decomposition
	$h_\xi^2=1+\mathbb P_{\ne0}(h_\xi^2)$, we  have
	\begin{align*}
		c_{n,\xi}
		={}&
		\int_{\mathbb T}
		\prod_{\gamma=0}^{n-2}
		h_\xi^2\left(
		\left\lceil\ell_q^{-\epsilon_0}\right\rceil^\gamma y
		\right)
		\left(
		1+\mathbb P_{\ne0}(h_\xi^2)
		\left(
		\left\lceil\ell_q^{-\epsilon_0}\right\rceil^{n-1}y
		\right)
		\right)\,\dd y\\
		={}&c_{n-1,\xi}
		+\int_{\mathbb T}
		\prod_{\gamma=0}^{n-2}
		h_\xi^2\left(
		\left\lceil\ell_q^{-\epsilon_0}\right\rceil^\gamma y
		\right)
		\mathbb P_{\ne0}(h_\xi^2)
		\left(
		\left\lceil\ell_q^{-\epsilon_0}\right\rceil^{n-1}y
		\right)\,\dd y .
	\end{align*}
	 Similar to \eqref{zero-mode-k}, we also obtain
	\begin{align}
		|c_{n,\xi}-c_{n-1,\xi}|
		{\lesssim
			\ell_q^{(n-1)N_0\epsilon_0}
			\ell_q^{-(n-2)N_0\epsilon_0}}{\lesssim
			\ell_q^{N_0\epsilon_0}.}
	\end{align}
	
	Since $c_{1,\xi}=1$, summing this recursive estimate yields
	\begin{align}
		|c_{n,\xi}-1|
		&\leq\sum_{r=2}^{n}|c_{r,\xi}-c_{r-1,\xi}|\lesssim_{L}
		\ell_q^{N_0\epsilon_0}.
		\label{zero-mode-xi}
	\end{align}
	
	Combining \eqref{zero-mode-k} and \eqref{zero-mode-xi}, we obtain, for the choice of $N_0$ in \eqref{N0} and $1\leq n\leq L-1$,
	\begin{align}
		\left|\int_{\mathbb T}g_{k,n}^2(y)\,\dd y-1\right|
		&\leq
		\left|\int_{\mathbb T}g_{k,n}^2(y)\,\dd y-c_{n,\xi}\right|
		+|c_{n,\xi}-1|\notag\\
		&\lesssim_{L}
		\ell_q^{N_0\epsilon_0}
		=o(\ell_q).
		\label{zero-mode-final}
	\end{align}
	Since $\int_{\mathbb T}g_{k,0}^2(y)\,\dd y=1$, estimate \eqref{zero-mode-final} also holds when $n=0$. Consequently, we decompose   $g_{k,i,n}^2(t,x)$ such that
	\begin{align}\label{intg}
		g_{k,i,n}^2(t,x)
		=1+\left(\int_{\mathbb T}g_{k,n}^2(y)\,\dd y-1\right)-f_{k,i,n}(t,x).
	\end{align}

	\subsection{Multiple iteration scheme}\label{MIS}
	We now introduce the multiple iteration scheme that generates the energy corrector incorporated into the perturbation. The construction combines a Newton--Nash iteration with a Picard iteration to eliminate the energy gap.
	
	\subsubsection{Newton-Nash iteration}
	We begin with the Newton--Nash iteration summarized in Proposition~\ref{assertion} and the proof is divided into three steps: the construction of the perturbations, the perturbation estimates, and the analysis of the Reynolds stress.
	
	\begin{prop}\label{assertion}
		Let $0 \leq N \leq N_0+2$, and suppose that $h_q(t)$ satisfies
		\begin{align}
			h_q(t) &\in \Big[ \delta_{q+1}, 100 \delta_{q+1}\Big], \quad    \|\partial_t h_q\|_{0} \lesssim \tau^{-1}_{q} \delta_{q+1}, \quad    \|\partial_{tt} h_q\|_{0} \lesssim \lambda^5_{q} \delta_{q+1}.\label{hq}
		\end{align}
		Then there exists a corresponding pair $(w_{q+1},\RR_{q+1})$, with $w_{q+1} := w^{(p_1)}_{q+1} + w^{(p_2)}_{q+1} + w^{(c)}_{q+1} + w^{(t)}_{q+1}$, such that
		\begin{align}
			&\lambda^{-N}_{q+1}\|\RR_{q+1}\|_N + \lambda^{-N-5}_{q+1} \|\partial_t \RR_{q+1}\|_N \leq \delta_{q+2} \lambda^{-4\alpha}_{q+1}, \label{M-energy-gap-m1} \\
			&\|w^{(p_1)}_{q+1}\|_0 + \lambda^{-N}_{q+1} \|w^{(p_1)}_{q+1}\|_N + \lambda^{-N-1}_{q+1} \|\partial_t w^{(p_1)}_{q+1}\|_N \leq \frac{M}{2} \delta_{q+1}^{1/2}, \label{M-energy-gap-m2} \\
			&\|w^{(p_2)}_{q+1}\|_0 + \|w^{(c)}_{q+1}\|_0 + \|w^{(t)}_{q+1}\|_0 \leq \delta^{1/2}_{q+1}, \label{M-energy-gap-m3} \\
			&\|w^{(p_2)}_{q+1}\|_N + \|w^{(c)}_{q+1}\|_N + \|w^{(t)}_{q+1}\|_N \leq \lambda^N_{q+1} \delta^{1/2}_{q+1}, \label{M-energy-gap-m4} \\
			&\|\partial_t w^{(p_2)}_{q+1}\|_N + \|\partial_t w^{(c)}_{q+1}\|_N + \|\partial_t w^{(t)}_{q+1}\|_N \leq \lambda^{N+1}_{q+1} \delta^{1/2}_{q+1}. \label{M-energy-gap-m5}
		\end{align}
	\end{prop}

	\noindent  \textbf{Step 1:\,\,Construction of the perturbations}

	(1) Set $\RRR_{q,0}=\RRR_q$. Given $\RRR_{q,n}$ with $0\leq n\leq L-1$, we define
	\begin{equation}
		a_{k,i,n} := \begin{cases}
			\widetilde{\chi}_i h^{1/2}_q a_k \Big(\nabla \Phi_i (\text{Id} - \frac{\RRR_q}{h_q}) \nabla \Phi_i^\text{T}\Big), & n = 0, \\
			\bar{\chi}_i \delta^{1/2}_{q+1,n}  a_k \Big(\nabla \Phi_i (\text{Id} - \frac{\RRR_{q,n}}{\delta_{q+1,n}}) \nabla \Phi_i^\text{T}\Big), & n = 1, 2, \dots, L-1,
		\end{cases} \label{av12}
	\end{equation}
	and
	\begin{align}
		A_{k,i,n} := (a_{k,i,n})^2 \nabla \Phi_i^{-1} (k \otimes k) \nabla \Phi_i^\text{-T}. \label{Av12}
	\end{align}
   Choose a further family of temporal cutoffs
    $\{\bar\eta_i\}_{0\leq i\leq i_{\max}-1}$ with
    $\bar\eta_i\in C^\infty([0,T];[0,1])$.
    For $0\leq i\leq i_{\max}-2$, we set
    \[
    \operatorname{supp}_t\bar\eta_i
    \subset
    [t_i-\tfrac23\tau_q,\,t_{i+1}+\tfrac23\tau_q]\cap[0,T],
    \]
    and
   $
    \bar\eta_i\equiv1$ on
    $[t_i-\tfrac13\tau_q,\,t_{i+1}+\tfrac13\tau_q]\cap[0,T].
    $
    For the last index $i=i_{max}-1$, we choose
    \[
    \operatorname{supp}_t\bar\eta_{i_{\max}-1}
    \subset
    [t_{i_{\max}-1}-\tfrac23\tau_q,\,T],
    \]
    and
$
    \bar\eta_{i_{\max}-1}\equiv1$
    on
    $[t_{i_{\max}-1}-\tfrac13\tau_q,\,T].
    $
  All these cutoffs satisfy
    $
    \|\partial_t^N\bar\eta_i\|_{L^\infty(0,T)}
    \lesssim_r\tau_q^{-N},
 $
	so that
	\[
	\bar{\eta}_i\widetilde{\chi}_i=\widetilde{\chi}_i,
	\qquad
	\bar{\eta}_i\bar{\chi}_i=\bar{\chi}_i,
	\qquad
	\partial_t\bar{\eta}_i\widetilde{\chi}_i
	=
	\partial_t\bar{\eta}_i\bar{\chi}_i=0.
	\]
	The cutoffs are also chosen so that
	\[
	\sum_i\widetilde{\chi}_i^2\RRR_q=\RRR_q,
	\qquad
	\sum_i\bar{\chi}_i^2\RRR_{q,n}=\RRR_{q,n},
	\quad 1\leq n\leq L-1,
	\]
	on the supports of the corresponding stresses.
	
	Define $\RRR_{q,n+1}$ by
	\begin{align}
		\RRR_{q,n+1} := \sum_i \partial_t \bar{\eta}_i \mathcal{R} w_{i,n+1}, \label{next Rqq}
	\end{align}
	where $w_{i,n+1}$ solves the following system on $[t_i-\tau_q,t_i+2\tau_q]\cap[0,T]$, 
	\begin{equation}
		\left\{
		\begin{alignedat}{2}
			&\del_t w_{i,n+1} + \vv_q \cdot\nabla w_{i,n+1} + w_{i,n+1}\cdot \nabla \vv_q + \nabla p_{i,n+1}  = \sum_{k \in \Lambda} \mathbb{P}_H\div(f_{k,i,n}  A_{k,i,n}), \\
			&w_{i,n+1}|_{t=t_i} = 0,\\
			&\operatorname{div}w_{i,n+1}=0.
		\end{alignedat}
		\right. \label{time}
	\end{equation}
	Here $f_{k,i,n}$ is defined by \eqref{fkin-traveling}, and $\fint_{\mathbb{T}^2}p_{i,n+1}\,\mathrm{d}x=0$.
	
	Proposition~\ref{est-RM} shows that $A_{k,i,n}$ and $w_{i,n+1}$ are well defined for $n = 0, 1, 2, \dots, L-1$. We define the Newton perturbation $\wtq$ by
	$$
	\wtq := \sum_{n=0}^{L-1} \sum_{i=0}^{i_\text{max}-1} \bar{\eta}_i(t) w_{i,n+1}, \quad \ptq := \sum_{n=0}^{L-1} \sum_{i=0}^{i_\text{max}-1} \bar{\eta}_i(t) p_{i,n+1}.
	$$
	
	(2) Let $\psi:\mathbb T\to\mathbb R$ be a smooth mean-zero cutoff supported in $[0,\frac1{100}]$ and normalized so that $\int_{\mathbb T}|\psi'(s)|^2\,\dd s=1$. Set $\phi=-\psi'$ and define the shear profile
	\begin{align}
		\phi_k(x) := \phi(\lambda_{q+1} k^{\perp} \cdot x), \quad \int_{\mathbb{T}^2} \phi_k^2 \dd x = 1. \label{phi-xi}
	\end{align}
	
	Let $\bar{\Phi}_i$ solve  the transport system
	\begin{equation}
		\left\{
		\begin{alignedat}{2}
			\del_t \bar{\Phi}_i + (\vv_q + \wtq) \cdot\nabla \bar{\Phi}_i &= 0, \\
			\bar{\Phi}_i(t_i, x) &= x,
		\end{alignedat}
		\right. \label{barphi}
	\end{equation}
	for $0 \leq i \leq i_\text{max}-1$.
	
	Set
	\[
	\bar g_{k,i,n}(t,x)
	:=g_{k,n}\!\left(\mu_{q+1}t-
	\ell_q^{-1}\left\lceil\ell_q^{-\epsilon_0}\right\rceil\bar\Phi_i^1(t,x)\right),
	\qquad 0\leq n\leq L-1.
	\]
	Similar to \eqref{bujiao}, we also have
	\begin{align}
		\bar{\chi}_i \bar{g}_{k,i,n}\,\bar{\chi}_{i'}\bar{g}_{k',i',n'}=0,\qquad
		\widetilde{\chi}_i \bar{g}_{k,i,n}\,\widetilde{\chi}_{i'}\bar{g}_{k',i',n'}=0,\qquad
		\bar{\chi}_i \bar{g}_{k,i,n}\,\widetilde{\chi}_{i'}\bar{g}_{k',i',n'}=0. \label{bujiaobar}
	\end{align}

	To cancel the Reynolds stresses \(\RRR_{q,n}\) for \(n = 0, 1, \dots, L-1\), define the principal perturbation \(\wpq\) by
	\begin{align}
		\wpq &:= \sum_{i; k \in \Lambda} \bar g_{k,i,0} \bar{a}_{k,i,0} \nabla \bar{\Phi}_i^{-1} \phi_k(\bar{\Phi}_i) k + \sum_{n=1}^{L-1} \sum_{i; k \in \Lambda} \bar g_{k,i,n} \bar{a}_{k,i,n} \nabla \bar{\Phi}_i^{-1} \phi_k(\bar{\Phi}_i) k \notag \\
		&:= w^{(p_1)}_{q+1} + w^{(p_2)}_{q+1}, \label{defi-wpq}
	\end{align}
where the term $w^{(p_1)}_{q+1}$, corresponding to $n=0$, provides the
	leading energy contribution, while $w^{(p_2)}_{q+1}$ collects
	the lower-order corrections from levels $n\geq1$, and
	\begin{align*}
		&
		\bar{a}_{k,i,0} = \widetilde{\chi}_i h_q^{1/2} a_k \Big(\nabla \bar{\Phi}_i \big(\operatorname{Id} - \frac{\RRR_{q,0}}{h_q}\big) \nabla \bar{\Phi}_i^\text{T}\Big),\\
		&
		\bar{a}_{k,i,n} = \bar{\chi}_i \delta_{q+1,n}^{1/2} a_k \Big(\nabla \bar{\Phi}_i \big(\operatorname{Id} - \frac{\RRR_{q,n}}{\delta_{q+1,n}}\big) \nabla \bar{\Phi}_i^\text{T}\Big),
		\qquad n\geq1.
	\end{align*}

	In analogy with \eqref{av12}, set
	\begin{align*}
		\bar{A}_{k,i,n} := \bar{a}_{k,i,n}^2 \nabla \bar{\Phi}_i^{-1} (k \otimes k) \nabla \bar{\Phi}_i^\text{-T}, \quad n \geq 0.
	\end{align*}
	
	(3) We choose a corrector $\wcq$ so that $\div(\wpq + \wcq) = 0$. Observe that
	\[
	\phi_k(\bar\Phi_i)\,\nabla\bar\Phi_i^{-1}k
	=
	\lambda_{q+1}^{-1}
	\nabla^\perp
	\left[
	\psi(\lambda_{q+1}k^{\perp}\cdot\bar\Phi_i)
	\right],
	\]
	where $\nabla^\perp = (-\partial_{x_2}, \partial_{x_1})$ and $\phi=-\psi'$. Hence
	\begin{align*}
		\wpq
		=  \lambda_{q+1}^{-1} \sum_{n=0}^{L-1} \sum_{i; k \in \Lambda} \bar{g}_{k,i,n} \bar{a}_{k,i,n} \nabla^\perp [\psi(\lambda_{q+1} k^{\perp} \cdot \bar{\Phi}_i)].
	\end{align*}
	
	Accordingly, define the divergence-free corrector $\wcq$ by
	\begin{align*}
		\wcq &:=  \lambda_{q+1}^{-1} \sum_{n=0}^{L-1} \sum_{i; k \in \Lambda} \nabla^\perp (\bar{g}_{k,i,n} \bar{a}_{k,i,n}) \psi(\lambda_{q+1} k^{\perp} \cdot \bar{\Phi}_i).
	\end{align*}
	Then
	\begin{align}
		\wpq + \wcq &= \lambda_{q+1}^{-1} \nabla^\perp \sum_{n=0}^{L-1} \sum_{i; k \in \Lambda} \bar{g}_{k,i,n} \bar{a}_{k,i,n} \psi(\lambda_{q+1} k^{\perp} \cdot \bar{\Phi}_i), \label{wpq-wcq}
	\end{align}
	and therefore $\div(w_{q+1}^{(p)} + w_{q+1}^{(c)}) = 0$.
	
	Finally, the full perturbation is defined by
	\begin{align}
		w_{q+1} := \wpq + \wcq + \wtq. \label{per-oper}
	\end{align}
	
	\noindent{\textbf{Step 2:\,\,Preliminaries for the perturbation estimates}}
	
	Since $\bar\Phi_i$ solves \eqref{barphi} and preserves measure, the function $\phi_k(\bar\Phi_i)$ is $\TTT^2$-periodic and has zero mean.
	
	A minor modification of Proposition~C.2 in \cite{BDSV} gives the following statement.
	
	\begin{cor}\label{-1jie}
		Let $H$ be a smooth vector field satisfying
		\[
		\|H\|_{N+\alpha}
		\lesssim C_H\ell_q^{-N}\ell_q^{-LN\epsilon_0},
		\qquad N\geq0,
		\]
		and suppose that \(\ell_q^{-1}\left\lceil\ell_q^{-\epsilon_0}\right\rceil^{L}\ll\lambda_{q+1}\). Let \(\bar\Phi_i\), \(\Phi_i\), and \(\phi_k\) be defined as above. Then, for \(t\in[t_i-\tau_q,t_i+2\tau_q]\cap[0,T]\),
		\begin{align*}
			\|\mathcal{R}(H \cdot \phi_k(\bar{\Phi}_i))\|_\alpha+ \|\mathcal{R}(H \cdot \phi_k(\Phi_i))\|_\alpha\lesssim C_H \lambda_{q+1} ^{-1+\alpha}.
		\end{align*}
	\end{cor}
	
	We estimate $\RRR_{q,n}$ by induction on $n$.
	
	\begin{prop}[Estimates for $\RRR_{q,n}$]\label{est-RM}
		For \(0\leq N\leq N_0+3\), estimate \eqref{estimate-wt} holds whenever \(0\leq n\leq L-1\), whereas \eqref{estimate-Rv} and \eqref{estimate-DRv} hold whenever \(0\leq n\leq L\),
		\begin{align}
			\|w_{i,n+1}\|_{N+\alpha}
			&\lesssim \delta_{q+1,n}\mu_{q+1}^{-1}\ell_q^{-N-1}\ell_q^{-(n+1)(N+1)\epsilon_0}\ell_q^{-5\alpha}, \label{estimate-wt} \\
			\|\RRR_{q,n}\|_{N+\alpha}
			&\lesssim \delta_{q+1,n}\ell_q^{-N}\ell_q^{-nN\epsilon_0}\ell_q^{\alpha}, \label{estimate-Rv} \\
			\|\matd{\vv_q} \RRR_{q,n}\|_{N+\alpha}
			&\lesssim \tau_q^{-1}\delta_{q+1,n}\ell_q^{-N}\ell_q^{-nN\epsilon_0}\ell_q^{\alpha}. \label{estimate-DRv}
		\end{align}
	\end{prop}
	
	\begin{proof}
		Since $\RRR_{q,0}=\RRR_q$, the bounds \eqref{estimate-Rv}--\eqref{estimate-DRv} at $n=0$ follow directly from Proposition~\ref{est-gluing}.
		
		Suppose that \eqref{estimate-Rv} and \eqref{estimate-DRv} hold at level \(n\), where \(0\leq n\leq L-1\). We prove the corresponding estimates at level \(n+1\). From \eqref{av12}--\eqref{Av12}, for $n=0$ we obtain
		\begin{align}
			\|a_{k,i,0}\|_{N+\alpha}&\lesssim \delta_{q+1}^{1/2} \ell_q^{-N-2\alpha}, \quad \|\matd{\vv_q} a_{k,i,0}\|_{N+\alpha} \lesssim \tau_q^{-1} \delta_{q+1}^{1/2} \ell_q^{-N-2\alpha}, \label{the estimate of coefficient 1} \\
			\|A_{k,i,0}\|_{N+\alpha}&\lesssim \delta_{q+1} \ell_q^{-N-2\alpha}, \quad \|\matd{\vv_q} A_{k,i,0}\|_{N+\alpha} \lesssim \tau_q^{-1} \delta_{q+1} \ell_q^{-N-2\alpha}, \label{the estimate of coefficient 2}
		\end{align}
		whereas, for \(1\leq n\leq L-1\) and \(0\leq N\leq N_0+3\),
		\begin{align}
			\|a_{k,i,n}\|_{N+\alpha}
			&\lesssim \delta_{q+1,n}^{1/2}\ell_q^{-N}\ell_q^{-nN\epsilon_0}\ell_q^{-2\alpha},
			\quad
			\|\matd{\vv_q}a_{k,i,n}\|_{N+\alpha}
				\lesssim \tau_q^{-1}\delta_{q+1,n}^{1/2}\ell_q^{-N}\ell_q^{-nN\epsilon_0}\ell_q^{-2\alpha}, \label{the estimate of coefficient 1-n} \\
			\|A_{k,i,n}\|_{N+\alpha}
			&\lesssim \delta_{q+1,n}\ell_q^{-N}\ell_q^{-nN\epsilon_0}\ell_q^{-2\alpha},
			\quad
			\|\matd{\vv_q}A_{k,i,n}\|_{N+\alpha}
				\lesssim \tau_q^{-1}\delta_{q+1,n}\ell_q^{-N}\ell_q^{-nN\epsilon_0}\ell_q^{-2\alpha}. \label{the estimate of coefficient 2-n}
		\end{align}
		
		Since $\|\vv_q\|_{1+\alpha}\lesssim\delta_q^{1/2}\lambda_q\ell_q^{-\alpha}$, the system \eqref{time} has a unique solution $(w_{i,n+1},p_{i,n+1})$ on $[t_i-\tau_q,t_i+2\tau_q]\cap[0,T]$. Let $y(t,x)$ denote the Lagrangian flow generated by $\vv_q$, defined by $\partial_t y(t,x)=\vv_q(t,y(t,x))$ and $y(t_i,x)=x$. Then
		\begin{equation} \label{Lar}
			\begin{alignedat}{2}
				& w_{i,n+1}(t, y(t,x)) = -\int_{t_i}^t \big(\mathbb{P}_H(w_{i,n+1} \cdot\nabla \vv_q)+( - \Delta)^{-1}(\nabla \div)(\vv_q \cdot\nabla w_{i,n+1})\big) (s, y(s,x)) ds \\
				& \quad\quad\quad\quad\quad\quad\quad\quad + \sum_{k \in \Lambda} \int_{t_i}^t \big[\mathbb{P}_H \div(f_{k,i,n} A_{k,i,n})\big](s, y(s,x)) ds,
			\end{alignedat}
		\end{equation}
		where $\mathbb{P}_H := {\rm Id} - \Delta^{-1} \nabla \div$ is the Helmholtz projection.
		
		Since $y(t,\cdot)$ and $y^{-1}(t,\cdot)$ are measure-preserving diffeomorphisms on $[t_i-\tau_q,t_i+2\tau_q]\cap[0,T]$, the following composition estimates hold for $-1<\gamma<1$,
		\begin{align}
			&\|f(t,\cdot)\|_{B_{\infty,\infty}^{\gamma}}\lesssim e^{\int_{t_i}^{t}\|\nabla\vv_q(s)\|_{\alpha}\,\dd s}\|f(t,y(t,\cdot))\|_{B_{\infty,\infty}^{\gamma}}\lesssim\|f(t,y(t,\cdot))\|_{B_{\infty,\infty}^{\gamma}},\label{e:f}\\
			&\|f(t,y(t,\cdot))\|_{B_{\infty,\infty}^{\gamma}}\lesssim e^{\int_{t_i}^{t}\|\nabla\vv_q(s)\|_{\alpha}\,\dd s}\|f(t,\cdot)\|_{B_{\infty,\infty}^{\gamma}}\lesssim\|f(t,\cdot)\|_{B_{\infty,\infty}^{\gamma}}.\label{e:f1}
		\end{align}
		For the proof of  \eqref{e:f}--\eqref{e:f1}, please see Lemma~2.7 of \cite{CNY24} or Lemma~A.1 of \cite{XZ15}.
		Using
			\begin{equation*}\label{Fkin-def}
				f_{k,i,n}
				=-\sum_{\substack{j_0,\ldots,j_n\in\mathbb Z\\ {\bf J}_n\neq0}}
				\frac{b_{j_n}^{(k)}\prod_{\gamma=0}^{n-1}b_{j_\gamma}^{(\xi)}}{2\pi i{\bf J}_n\mu_{q+1}}
				e^{-2\pi i{\bf J}_n\ell_q^{-1}\left\lceil\ell_q^{-\epsilon_0}\right\rceil\Phi_i^1}
				\partial_t e^{2\pi i{\bf J}_n\mu_{q+1}t},
			\end{equation*}
			and integrating by parts in time, we obtain
			\begin{align}
				&\Big\|\int_{t_i}^t
				\big[\mathbb P_H\div(f_{k,i,n}A_{k,i,n})\big](s,y(s,x))\,\dd s\Big\|_\alpha
				\notag\\
				={}&\left\|\begin{aligned}&-
				\sum_{\substack{j_0,\ldots,j_n\in\mathbb Z\\ {\bf J}_n\neq0}}\frac{b_{j_n}^{(k)}\prod_{\gamma=0}^{n-1}b_{j_\gamma}^{(\xi)}}{2\pi i{\bf J}_n\mu_{q+1}}\times\int_{t_i}^t
				\partial_s e^{2\pi i{\bf J}_n\mu_{q+1}s}
				\Big[\mathbb P_H\div\!\left(
				e^{-2\pi i{\bf J}_n\ell_q^{-1}\left\lceil\ell_q^{-\epsilon_0}\right\rceil\Phi_i^1}
				A_{k,i,n}
				\right)\Big](s,y(s,x))\,\dd s\end{aligned}\right\|_\alpha
				\notag\\
				\lesssim{}&\frac{1}{2\pi\mu_{q+1}}\Big\|\Big[
				e^{2\pi i{\bf J}_n\mu_{q+1}s}
				\cdot\mathbb P_H\div\!\left(
				e^{-2\pi i{\bf J}_n\ell_q^{-1}\left\lceil\ell_q^{-\epsilon_0}\right\rceil\Phi_i^1}
				A_{k,i,n}
				\right)(s,y(s,x))
				\Big]_{s=t_i}^{s=t}\Big\|_\alpha
				\notag\\
				&+\frac{1}{2\pi\mu_{q+1}}
				\Big\|\int_{t_i}^t e^{2\pi i{\bf J}_n\mu_{q+1}s}
				\Big[D_{t,\vv_q}\mathbb P_H\div\!\left(
				e^{-2\pi i{\bf J}_n\ell_q^{-1}\left\lceil\ell_q^{-\epsilon_0}\right\rceil\Phi_i^1}
				A_{k,i,n}
				\right)\Big](s,y(s,x))\,\dd s\Big\|_\alpha
				\notag\\
				\lesssim{}&\frac{1}{2\pi\mu_{q+1}}
				\Big\|\Big[
				e^{2\pi i{\bf J}_n\mu_{q+1}s}
				\mathbb P_H\div\!\left(
				e^{-2\pi i{\bf J}_n\ell_q^{-1}\left\lceil\ell_q^{-\epsilon_0}\right\rceil\Phi_i^1}
				A_{k,i,n}
				\right)(s,y(s,x))
				\Big]_{s=t_i}^{s=t}\Big\|_\alpha
				\notag\\
				&+\frac{1}{2\pi\mu_{q+1}}
				\Big\|\int_{t_i}^t e^{2\pi i{\bf J}_n\mu_{q+1}s}
				\Big[\mathbb P_H\div\!\left(
				e^{-2\pi i{\bf J}_n\ell_q^{-1}\left\lceil\ell_q^{-\epsilon_0}\right\rceil\Phi_i^1}
				D_{t,\vv_q}A_{k,i,n}
				\right)\Big](s,y(s,x))\,\dd s\Big\|_\alpha
				\notag\\
				&+\frac{1}{2\pi\mu_{q+1}}
				\Big\|\int_{t_i}^t e^{2\pi i{\bf J}_n\mu_{q+1}s}
				\Big([\vv_q\cdot\nabla,\mathbb P_H\div]\!\left(
				e^{-2\pi i{\bf J}_n\ell_q^{-1}\left\lceil\ell_q^{-\epsilon_0}\right\rceil\Phi_i^1}
				A_{k,i,n}
				\right)\Big)(s,y(s,x))\,\dd s\Big\|_\alpha
				\notag\\
				=:{}& I_1+I_2+I_3.\label{A:e}
			\end{align}
			For the boundary term $I_1$,
			\begin{align}
				I_1
				\lesssim&
				\mu_{q+1}^{-1}
				\left\|
				\mathbb P_H\div
				\left(
				e^{-2\pi i{\bf J}_n\ell_q^{-1}
					\left\lceil\ell_q^{-\epsilon_0}\right\rceil
					\Phi_i^1}
				A_{k,i,n}
				\right)
				\right\|_\alpha
				\notag	\\
				\lesssim&
				\mu_{q+1}^{-1}
				\left\|
				e^{-2\pi i{\bf J}_n\ell_q^{-1}
					\left\lceil\ell_q^{-\epsilon_0}\right\rceil
					\Phi_i^1}
				A_{k,i,n}
				\right\|_{1+\alpha}
				\notag	\\
				\lesssim&
				\mu_{q+1}^{-1}
				\Big(
				\left\|
				e^{-2\pi i{\bf J}_n\ell_q^{-1}
					\left\lceil\ell_q^{-\epsilon_0}\right\rceil
					\Phi_i^1}
				\right\|_{1+\alpha}
				\|A_{k,i,n}\|_0
				+
				\|A_{k,i,n}\|_{1+\alpha}
				\Big)
				\notag	\\
				\lesssim&
				\mu_{q+1}^{-1}
				\delta_{q+1,n}
				\ell_q^{-1-2\alpha}
				\ell_q^{-(n+1)\epsilon_0}.\label{I_{1}}
			\end{align}
			Likewise, for $I_2$,
			\begin{align}
				I_2
				\lesssim&
				\mu_{q+1}^{-1}
				\int_{t_i}^{t}
				\left\|
				\mathbb P_H\div
				\left(
				e^{-2\pi i{\bf J}_n\ell_q^{-1}
					\left\lceil\ell_q^{-\epsilon_0}\right\rceil
					\Phi_i^1}
				D_{t,\vv_q}A_{k,i,n}
				\right)
				\right\|_\alpha
				\,\mathrm ds
				\notag\\
				\lesssim&
				\mu_{q+1}^{-1}
				\int_{t_i}^{t}
				\left\|
				e^{-2\pi i{\bf J}_n\ell_q^{-1}
					\left\lceil\ell_q^{-\epsilon_0}\right\rceil
					\Phi_i^1}
				D_{t,\vv_q}A_{k,i,n}
				\right\|_{1+\alpha}
				\,\mathrm ds
				\notag\\
				\lesssim&
				\mu_{q+1}^{-1}
				\int_{t_i}^{t}
				\Big(
				\left\|
				e^{-2\pi i{\bf J}_n\ell_q^{-1}
					\left\lceil\ell_q^{-\epsilon_0}\right\rceil
					\Phi_i^1}
				\right\|_{1+\alpha}
				\|D_{t,\vv_q}A_{k,i,n}\|_0
				+
				\|D_{t,\vv_q}A_{k,i,n}\|_{1+\alpha}
				\Big)
				\,\mathrm ds
				\notag\\
				\lesssim&
				\mu_{q+1}^{-1}
				\delta_{q+1,n}
				\ell_q^{-1-2\alpha}
				\ell_q^{-(n+1)\epsilon_0}.\label{I_{2}}
			\end{align}
			We shall also use the commutator estimate (see, \cite{Ba,MNY}),
			\begin{align}\label{rm-T-commutator-est}
				\|[\vv_q\cdot\nabla,\mathbb P_H\div]F\|_{N+\alpha}
				\lesssim
				\|\nabla\vv_q\|_\alpha\|F\|_{N+1+\alpha}
				+\|\nabla\vv_q\|_{N+\alpha}\|F\|_{1+\alpha}.
			\end{align}
			It follows that
			\begin{align}
				I_3
				\lesssim&
				\mu_{q+1}^{-1}
				\int_{t_i}^{t}
				\left\|
				[\vv_q\cdot\nabla,\mathbb P_H\div]
				\left(
				e^{-2\pi i{\bf J}_n\ell_q^{-1}
					\left\lceil\ell_q^{-\epsilon_0}\right\rceil
					\Phi_i^1}
				A_{k,i,n}
				\right)
				\right\|_\alpha
				\,\mathrm ds
				\notag\\
				\lesssim&
				\mu_{q+1}^{-1}
				\int_{t_i}^{t}
				\|\nabla\vv_q\|_\alpha
				\left\|
				e^{-2\pi i{\bf J}_n\ell_q^{-1}
					\left\lceil\ell_q^{-\epsilon_0}\right\rceil
					\Phi_i^1}
				A_{k,i,n}
				\right\|_{1+\alpha}
				\,\mathrm ds
				\notag\\
				\lesssim&
				\mu_{q+1}^{-1}\delta_{q+1,n}\ell_q^{-1}\ell_q^{-(n+1)\epsilon_0}\ell_q^{-5\alpha}.\label{I_{3}}
			\end{align}Combining \eqref{I_{1}}--\eqref{I_{3}}, we obtain
			\begin{align}\label{rm-source-CN}
				\left\|
				\int_{t_i}^t
				\big[\mathbb P_H\operatorname{div}(f_{k,i,n}A_{k,i,n})\big](s,y(s,\cdot))\,\dd s
				\right\|_\alpha
				\lesssim
				\mu_{q+1}^{-1}\delta_{q+1,n}\ell_q^{-1}\ell_q^{-(n+1)\epsilon_0}\ell_q^{-5\alpha}.
		\end{align}
		
		The standard commutator estimates (see, \cite{Ba,MNY}) and the identity $\div w_{i,n+1}=0$ give
		\begin{align}
			&
			\|\mathbb{P}_H(w_{i,n+1}\cdot\nabla \vv_q)(s, y(s,x))\|_\alpha
			+\|(-\Delta)^{-1}(\nabla\div)(\vv_q\cdot\nabla w_{i,n+1})(s, y(s,x))\|_\alpha
			\notag\\
			\lesssim&
			\|w_{i,n+1}\cdot\nabla\vv_q\|_\alpha
			+
			\|(-\Delta)^{-1}(\nabla\div)(\vv_q\cdot\nabla w_{i,n+1})\|_\alpha
			\notag\\
			\lesssim&
			\|w_{i,n+1}\|_\alpha\|\nabla\vv_q\|_0
			+
			\|w_{i,n+1}\|_0\|\nabla\vv_q\|_\alpha+
			\left\|
			\left[
			(-\Delta)^{-1}(\nabla\div),
			\vv_q\cdot\nabla
			\right]
			w_{i,n+1}
			\right\|_\alpha
			\notag\\
			\lesssim&
			\|w_{i,n+1}\|_\alpha\|\nabla\vv_q\|_0
			+
			\|w_{i,n+1}\|_0\|\nabla\vv_q\|_\alpha+
			\|\nabla\vv_q\|_0
			\|w_{i,n+1}\|_\alpha
			+
			\|\nabla\vv_q\|_\alpha
			\|w_{i,n+1}\|_0
			\notag\\
			\lesssim&
			\|\vv_q\|_{1+\alpha}
			\|w_{i,n+1}\|_\alpha .
			\label{e:a}
		\end{align}
		
		Consequently, \eqref{Lar} yields
		\begin{align}
			\|w_{i,n+1}(x,t)\|_\alpha
			\lesssim&
			\int_{t_i}^{t}
			\Big(
			\|\mathbb{P}_H(w_{i,n+1}\cdot\nabla\vv_q)(s,y(s,x))\|_\alpha
			\notag\\
			&
			+
			\|(-\Delta)^{-1}(\nabla\div)
			(\vv_q\cdot\nabla w_{i,n+1})(s,y(s,x))\|_\alpha
			\Big)\,\dd s
			\notag\\
			&+
			\left\|
			\int_{t_i}^{t}
			\sum_{k\in\Lambda}
			\mathbb P_H\div
			(f_{k,i,n}A_{k,i,n})
			(s,y(s,x))\,\dd s
			\right\|_\alpha
			\notag\\
			\lesssim&
			\int_{t_i}^{t}
			\|\vv_q\|_{1+\alpha}
			\|w_{i,n+1}\|_\alpha
			\,\dd s+
			\mu_{q+1}^{-1}
			\delta_{q+1,n}
			\ell_q^{-1}
			\ell_q^{-(n+1)\epsilon_0}
			\ell_q^{-5\alpha}.
			\label{e:Lar1}
		\end{align}
		Gronwall's inequality, together with \(\int_{[t_i-\tau_q,t_i+2\tau_q]\cap[0,T]}\|\vv_q\|_{1+\alpha}\,\dd s\lesssim1\), gives
		\begin{align}\label{est-wn-0}
				\|w_{i,n+1}\|_\alpha
				\lesssim \mu_{q+1}^{-1}\delta_{q+1,n}\ell_q^{-1}\ell_q^{-(n+1)\epsilon_0}\ell_q^{-5\alpha}.
		\end{align}

		We next estimate \(\|\mathcal{R}w_{i,n+1}\|_\alpha\).  Repeating the temporal integration by parts from \eqref{A:e}, now in \(B^{-1+\alpha}_{\infty,\infty}\) (please refer to \cite{MNY1} for the definition of  Besov space), yields
			\begin{align}\label{A:e-besov}
				&
				\Big\|
				\int_{t_i}^t
				\big[\mathbb P_H\div(f_{k,i,n}A_{k,i,n})\big]
				(s,y(s,\cdot))\,\dd s
				\Big\|_{B^{-1+\alpha}_{\infty,\infty}}
				\notag\\
				\lesssim&
				\mu_{q+1}^{-1}
				\Big\|
				\Big[
				e^{2\pi i{\bf J}_n\mu_{q+1}s}
				\mathbb P_H\div
				\left(
				e^{-2\pi i{\bf J}_n\ell_q^{-1}
					\left\lceil\ell_q^{-\epsilon_0}\right\rceil
					\Phi_i^1}
				A_{k,i,n}
				\right)
				(s,y(s,\cdot))
				\Big]_{s=t_i}^{s=t}
				\Big\|_{B^{-1+\alpha}_{\infty,\infty}}
				\notag\\
				&+
				\mu_{q+1}^{-1}
				\Big\|
				\int_{t_i}^{t}
				e^{2\pi i{\bf J}_n\mu_{q+1}s}
				D_{t,\vv_q}
				\mathbb P_H\div
				\left(
				e^{-2\pi i{\bf J}_n\ell_q^{-1}
					\left\lceil\ell_q^{-\epsilon_0}\right\rceil
					\Phi_i^1}
				A_{k,i,n}
				\right)
				(s,y(s,\cdot))
				\,\dd s
				\Big\|_{B^{-1+\alpha}_{\infty,\infty}}
				\notag\\
				\lesssim&
				\mu_{q+1}^{-1}
				\Big\|
				\Big[
				e^{2\pi i{\bf J}_n\mu_{q+1}s}
				\mathbb P_H\div
				\left(
				e^{-2\pi i{\bf J}_n\ell_q^{-1}
					\left\lceil\ell_q^{-\epsilon_0}\right\rceil
					\Phi_i^1}
				A_{k,i,n}
				\right)
				(s,y(s,\cdot))
				\Big]_{s=t_i}^{s=t}
				\Big\|_{B^{-1+\alpha}_{\infty,\infty}}
				\notag\\
				&+
				\mu_{q+1}^{-1}
				\Big\|
				\int_{t_i}^{t}
				e^{2\pi i{\bf J}_n\mu_{q+1}s}
				\mathbb P_H\div
				\left(
				e^{-2\pi i{\bf J}_n\ell_q^{-1}
					\left\lceil\ell_q^{-\epsilon_0}\right\rceil
					\Phi_i^1}
				D_{t,\vv_q}A_{k,i,n}
				\right)
				(s,y(s,\cdot))
				\,\dd s
				\Big\|_{B^{-1+\alpha}_{\infty,\infty}}
				\notag\\
				&+
				\mu_{q+1}^{-1}
				\Big\|
				\int_{t_i}^{t}
				e^{2\pi i{\bf J}_n\mu_{q+1}s}
				[\vv_q\cdot\nabla,\mathbb P_H\div]
				\left(
				e^{-2\pi i{\bf J}_n\ell_q^{-1}
					\left\lceil\ell_q^{-\epsilon_0}\right\rceil
					\Phi_i^1}
				A_{k,i,n}
				\right)
				(s,y(s,\cdot))
				\,\dd s
				\Big\|_{B^{-1+\alpha}_{\infty,\infty}} \notag\\
				\mathrel{:=}{}&I_4+I_5+I_6.
			\end{align}
			For $I_4$,
			\begin{align}
				I_4
				\lesssim&
				\mu_{q+1}^{-1}
				\left\|
				\mathbb P_H\div
				\left(
				e^{-2\pi i{\bf J}_n\ell_q^{-1}
					\left\lceil\ell_q^{-\epsilon_0}\right\rceil
					\Phi_i^1}
				A_{k,i,n}
				\right)
				\right\|_{B^{-1+\alpha}_{\infty,\infty}}
				\notag\\
				\lesssim&
				\mu_{q+1}^{-1}
				\left\|
				e^{-2\pi i{\bf J}_n\ell_q^{-1}
					\left\lceil\ell_q^{-\epsilon_0}\right\rceil
					\Phi_i^1}
				A_{k,i,n}
				\right\|_{B^{\alpha}_{\infty,\infty}}
				\notag\\
				\lesssim&
				\mu_{q+1}^{-1}
				\Big(
				\left\|
				e^{-2\pi i{\bf J}_n\ell_q^{-1}
					\left\lceil\ell_q^{-\epsilon_0}\right\rceil
					\Phi_i^1}
				\right\|_{B^\alpha_{\infty,\infty}}
				\|A_{k,i,n}\|_0
				+
				\|A_{k,i,n}\|_{B^\alpha_{\infty,\infty}}
				\Big)
				\notag\\
				\lesssim&
				\mu_{q+1}^{-1}
				\delta_{q+1,n}
				\ell_q^{-5\alpha}.\label{I_{4}}
			\end{align}
			Similarly, $I_5$ satisfies
			\begin{align}
				I_5\lesssim&
				\mu_{q+1}^{-1}
				\int_{t_i}^{t}
				\left\|
				\mathbb P_H\div
				\left(
				e^{-2\pi i{\bf J}_n\ell_q^{-1}
					\left\lceil\ell_q^{-\epsilon_0}\right\rceil
					\Phi_i^1}
				D_{t,\vv_q}A_{k,i,n}
				\right)
				\right\|_{B^{-1+\alpha}_{\infty,\infty}}
				\,\dd s
				\notag\\
				\lesssim&
				\mu_{q+1}^{-1}
				\int_{t_i}^{t}
				\left\|
				e^{-2\pi i{\bf J}_n\ell_q^{-1}
					\left\lceil\ell_q^{-\epsilon_0}\right\rceil
					\Phi_i^1}
				D_{t,\vv_q}A_{k,i,n}
				\right\|_{B^\alpha_{\infty,\infty}}
				\,\dd s
				\notag\\
				\lesssim&
				\mu_{q+1}^{-1}
				\int_{t_i}^{t}
				\Big(
				\left\|
				e^{-2\pi i{\bf J}_n\ell_q^{-1}
					\left\lceil\ell_q^{-\epsilon_0}\right\rceil
					\Phi_i^1}
				\right\|_{B^\alpha_{\infty,\infty}}
				\|D_{t,\vv_q}A_{k,i,n}\|_0
				+
				\|D_{t,\vv_q}A_{k,i,n}\|_{B^\alpha_{\infty,\infty}}
				\Big)
				\,\dd s
				\notag\\
				\lesssim&
				\mu_{q+1}^{-1}
				\delta_{q+1,n}
				\Big(
				\ell_q^{-\alpha}
				\ell_q^{-\alpha(n+1)\epsilon_0}
				+
				\ell_q^{-2\alpha}
				\Big)
				\notag\\
				\lesssim&
				\mu_{q+1}^{-1}
				\delta_{q+1,n}
				\ell_q^{-5\alpha}.\label{I_5}
			\end{align}
			Finally, applying the commutator estimate (see, \cite{Ba,MNY}) once more, we obtain
			\begin{align}
				I_6
				\lesssim&
				\mu_{q+1}^{-1}
				\int_{t_i}^{t}
				\Big\|
				[\vv_q\cdot\nabla,\mathbb P_H\div]
				\left(
				e^{-2\pi i{\bf J}_n\ell_q^{-1}
					\left\lceil\ell_q^{-\epsilon_0}\right\rceil
					\Phi_i^1}
				A_{k,i,n}
				\right)
				\Big\|_{B^{-1+\alpha}_{\infty,\infty}}
				\,\dd s
				\notag\\
				\lesssim&
				\mu_{q+1}^{-1}
				\int_{t_i}^{t}
				\|\nabla\vv_q\|_\alpha
				\left\|
				e^{-2\pi i{\bf J}_n\ell_q^{-1}
					\left\lceil\ell_q^{-\epsilon_0}\right\rceil
					\Phi_i^1}
				A_{k,i,n}
				\right\|_{B^\alpha_{\infty,\infty}}
				\,\dd s
				\notag\\
				\lesssim&
				\mu_{q+1}^{-1}
				\delta_{q+1,n}
				\ell_q^{-5\alpha}.\label{I_6}
			\end{align}
			Substituting the estimates \eqref{I_{4}}--\eqref{I_6} into \eqref{A:e-besov}, we obtain
		\begin{align}\label{rm-source-Besov}
			\left\|
			\int_{t_i}^t
			\big[\mathbb P_H\operatorname{div}(f_{k,i,n}A_{k,i,n})\big](s,y(s,\cdot))\,\dd s
			\right\|_{B^{-1+\alpha}_{\infty,\infty}}
			\lesssim
			\mu_{q+1}^{-1}\delta_{q+1,n}\ell_q^{-5\alpha}.
		\end{align}
		By the Bony decomposition (see \cite{Ba}) and the divergence-free condition on  $w_{i,n+1}$, we have
		\begin{align}
			&
			\|\mathbb{P}_H(w_{i,n+1}\cdot\nabla \vv_q)(s, y(s,x))\|_{B^{-1+\alpha}_{\infty,\infty}}
			+
			\|(-\Delta)^{-1}(\nabla\div)
			(\vv_q\cdot\nabla w_{i,n+1})(s, y(s,x))\|_{B^{-1+\alpha}_{\infty,\infty}}
			\notag\\
			\lesssim&
			\|w_{i,n+1}\cdot\nabla\vv_q\|_{B^{-1+\alpha}_{\infty,\infty}}
			+
			\|(-\Delta)^{-1}(\nabla\div)
			(\vv_q\cdot\nabla w_{i,n+1})\|_{B^{-1+\alpha}_{\infty,\infty}}
			\notag\\
			\lesssim&
			\|\vv_q\|_{1+\alpha}
			\|w_{i,n+1}\|_{B^{-1+\alpha}_{\infty,\infty}} .
			\label{e:a-besov}
		\end{align}
		
		It follows from \eqref{Lar} that
		\begin{align}
			&\quad\|w_{i,n+1} (x,t)\|_{B^{-1+\alpha}_{\infty,\infty}} \notag\\
			&\lesssim\int_{t_{i}}^{t}\|
			\Big(\mathbb{P}_H(w_{i,n+1}\cdot \nabla \vv_q)+(-\Delta)^{-1}(\nabla\div)(\vv_q\cdot \nabla w_{i,n+1})\Big)(s,y(s,x))\|_{B^{-1+\alpha}_{\infty,\infty}}\dd s\notag\\
			&\quad+\left\|\sum_{k\in\Lambda}\int_{t_i}^{t}[\mathbb P_H\div(f_{k,i,n}A_{k,i,n})](s,y(s,x))\,\dd s\right\|_{B^{-1+\alpha}_{\infty,\infty}}\notag\\
			&\lesssim
			\int_{t_{i}}^{t}\|w_{i,n+1}\|_{B^{-1+\alpha}_{\infty,\infty}}\|\vv_q\|_{1+\alpha}\dd s+\mu_{q+1}^{-1}\delta_{q+1,n}\ell_q^{-5\alpha}. \label{e:Lar1-besov}
		\end{align}
		Applying the integral form of Gr\"{o}nwall's inequality and using $\int_{t_i}^{t}\|\vv_q\|_{1+\alpha}\,\dd s\ll1$, we obtain
		\begin{align}\label{est-wn-minus}
			\|w_{i,n+1}\|_{B^{-1+\alpha}_{\infty,\infty}}\lesssim \mu_{q+1}^{-1}\delta_{q+1,n}\ell_q^{-5\alpha}e^{\int_{t_{i}}^{t}(\|\vv_q\|_{1+\alpha}+1)\,\dd s}\lesssim \mu_{q+1}^{-1}\delta_{q+1,n}\ell_q^{-5\alpha}.
		\end{align}
		
		Consequently,
		\begin{align}\label{est-wn--1}
			\|\mathcal{R}w_{i,n+1}\|_{\alpha} \lesssim \|w_{i,n+1}\|_{B^{-1+\alpha}_{\infty,\infty}}\lesssim \mu^{-1}_{q+1} \delta_{q+1,n}\ell_{q}^{-5\alpha}.
		\end{align}
		
		For \(1\leq N\leq N_0+3\), applying \(\partial^\gamma\) with \(|\gamma|=N\) to \eqref{time} yields
		\begin{align}
			&\del_t \partial^\gamma w_{i,n+1} + \vv_q \cdot\nabla  \partial^\gamma w_{i,n+1}
			+ \sum_{\gamma_1 + \gamma_2 = \gamma,\, \gamma_1 \neq \gamma}
			\partial^{\gamma_2} \vv_q \cdot\nabla \partial^{\gamma_1} w_{i,n+1} \notag \\
			&\quad+ \sum_{\gamma_1 + \gamma_2 = \gamma}
			\partial^{\gamma_2} w_{i,n+1}\cdot \nabla \partial^{\gamma_1} \vv_q
			+ \nabla  \partial^\gamma p_{i,n+1}
			= {\sum_{k \in \Lambda} \partial^\gamma \mathbb P_H\div (f_{k,i,n}A_{k,i,n}).}
		\end{align}
		Applying \(\mathbb P_H\) and repeating the preceding Lagrangian estimate yields
		\begin{align}\label{est-wn-N}
				\|w_{i,n+1}\|_{N+\alpha}
				\lesssim \mu_{q+1}^{-1}\delta_{q+1,n}
				\ell_q^{-N-1}\ell_q^{-(n+1)(N+1)\epsilon_0}\ell_q^{-5\alpha},
				\quad 1\leq N\leq N_0+3.
		\end{align}

		Recalling \eqref{next Rqq} and using \eqref{est-wn--1} and \eqref{est-wn-N}, we obtain
		\begin{align}\label{RMn+1}
			\|\RRR_{q,n+1}\|_{N+\alpha}
			&\lesssim \tau_q^{-1}\mu_{q+1}^{-1}\delta_{q+1,n}
			\ell_q^{-N}\ell_q^{-(n+1)N\epsilon_0}\ell_q^{-5\alpha}\notag\\
			&=\delta_{q+1,n+1}\ell_q^{-N}\ell_q^{-(n+1)N\epsilon_0}\ell_q^{5\alpha+L\epsilon_0}\notag\\
			&\lesssim \delta_{q+1,n+1}\ell_q^{-N}\ell_q^{-(n+1)N\epsilon_0}\ell_q^\alpha.
		\end{align}
		 This proves  \eqref{estimate-Rv}  at level  \(n+1\). Together with  \eqref{est-wn-0}  and  \eqref{est-wn-N}, it also gives  \eqref{estimate-wt}  at the current Newton step.
		
		 Finally, to estimate  \(\|\matd{\vv_q}\RRR_{q,n+1}\|_{N+\alpha}\), we write
		\begin{align*}
			\matd{\vv_q} \RRR_{q,n+1}
			={}& \sum_i \partial_t^2 \bar{\eta}_i \mathcal{R} w_{i,n+1}
			+ \sum_i \partial_t \bar{\eta}_i \mathcal{R} \matd{\vv_q} w_{i,n+1} \\
			&+ \sum_i \partial_t \bar{\eta}_i [\vv_q\cdot \nabla,\mathcal{R}] w_{i,n+1}.
		\end{align*}
		 Moreover, since 
		\(
		\partial_t\bar\eta_i\widetilde\chi_i
		=
		\partial_t\bar\eta_i\bar\chi_i=0
		\), multiplication by  \(\partial_t\bar\eta_i\)  gives
		\begin{align*}
			\partial_t\bar\eta_i\matd{\vv_q}w_{i,n+1}
			={}&-\partial_t\bar\eta_i\mathbb P_H\big((w_{i,n+1}\cdot\nabla)\vv_q\big){-\partial_t\bar\eta_i(-\Delta)^{-1}\nabla\div\big((\vv_q\cdot\nabla)w_{i,n+1}\big).}
		\end{align*}
		 The source term vanishes because 
		\(\partial_t\bar\eta_iA_{k,i,n}=0\). Hence the preceding estimates and the standard commutator bound yield
		\begin{align}\label{DRn+1}
			&\quad\|\matd{\vv_q}\RRR_{q,n+1}\|_{N+\alpha}\notag\\
			&{\lesssim
				\sum_i\Big(
				\|\partial_t^2\bar\eta_i\|_0\|\mathcal Rw_{i,n+1}\|_{N+\alpha}
				+\|\partial_t\bar\eta_i\|_0\|\mathcal R\matd{\vv_q}w_{i,n+1}\|_{N+\alpha}
				+\|\partial_t\bar\eta_i\|_0
				\|[\vv_q\cdot\nabla,\mathcal R]w_{i,n+1}\|_{N+\alpha}
				\Big)}\notag\\
			&\lesssim \tau_q^{-1}\delta_{q+1,n+1}
			\ell_q^{-N}\ell_q^{-(n+1)N\epsilon_0}\ell_q^\alpha.
		\end{align}
		 This proves  \eqref{estimate-DRv}  at level  \(n+1\).
		
		This completes the proof of Proposition \ref{est-RM}.
	\end{proof}
	
	 The preceding estimates for  $w_{i,n+1}$  immediately yield the following bounds for  $w^{(t)}_{q+1}$.
	\begin{cor}\label{est-wtq}
		 For  \(0 \leq N \leq N_0+3\), the definition  \(\wtq := \sum\limits_{n=0}^{L-1} \sum\limits_{i=0}^{i_{\max}-1} \bar{\eta}_i(t) w_{i,n+1}\)  gives
		\begin{align}
			\|\mathcal{R} \wtq\|_0
			&\lesssim \mu_{q+1}^{-1} \delta_{q+1}\ell_q^{-5\alpha}, \\
			\|\wtq\|_N
			&\lesssim \mu_{q+1}^{-1}\ell_q^{-N-1}\ell_q^{-5\alpha}
				\sum_{n=0}^{L-1}\delta_{q+1,n}\ell_q^{-(n+1)(N+1)\epsilon_0}
				\ll \lambda_{q+1}^N\delta_{q+1}^{1/2}, \\
			\|\matd{\vv_q + \wtq} \RRR_{q,n}\|_N
			&\lesssim \tau_q^{-1}\delta_{q+1,n}\ell_q^{-N}\ell_q^{-nN\epsilon_0}.
		\end{align}
		In particular, for \(N=0,1\), the geometric decay of \(\delta_{q+1,n}\) yields
			\[
			\|\wtq\|_N
			\lesssim \mu_{q+1}^{-1}\delta_{q+1}\ell_q^{-N-1}\ell_q^{-(N+1)\epsilon_0}\ell_q^{-5\alpha}.
			\]
	\end{cor}

	\begin{prop}[ Energy interaction of the transport perturbation]
		\label{prop:transport-energy}
		The transport perturbation satisfies
		\begin{align}
			\left|
			\partial_t
			\int_{\mathbb T^2}
			2\wtq\cdot\bar v_q\,\dd x
			\right|
			&\lesssim
			\ell_q^\alpha\tau_q^{-1}\delta_{q+1},
			\label{transport-energy-time}\\
            \left|
			\partial^2_t
			\int_{\mathbb T^2}
			2\wtq\cdot\bar v_q\,\dd x
			\right|
			&\lesssim
			\lambda^3_{q}\tau_q^{-1}\delta_{q+1},
			\label{2transport-energy-time}   \\
			\left|
			\int_{\mathbb T^2}
			2\wtq\cdot\bar v_q\,\dd x
			\right|
			&\lesssim
			\ell_q^\alpha\delta_{q+1}.
			\label{transport-energy}
		\end{align}
	\end{prop}
	\begin{proof}
		 We first prove  \eqref{transport-energy-time}. Recall that
		\[
		\wtq
		=
		\sum_{n=0}^{L-1}\sum_{i=0}^{i_{\max}-1}
		\bar{\eta}_i(t)w_{i,n+1}.
		\]
		 For each  $w_{i,n+1}$, the construction yields
		\[
		\partial_t\int_{\mathbb T^2}w_{i,n+1}\bar v_q\,\dd x
		=
		\int_{\mathbb T^2}\div\RRR_q\cdot w_{i,n+1}\,\dd x
		+
		\int_{\mathbb T^2}\sum_{k\in\Lambda}
		\bar v_q\div(f_{k,i,n}A_{k,i,n})\,\dd x .
		\]
		
		 For the first term on the right-hand side of the above equality, Proposition  \ref{est-RM}  and integration by parts on  $\mathbb T^2$  give
		\[
		\left|
		\int_{\mathbb T^2}\div\RRR_q\cdot w_{i,n+1}\,\dd x
		\right|
		\lesssim
		\|\RRR_q\|_0\|\nabla w_{i,n+1}\|_0 \lesssim
		\mu_{q+1}^{-1}\delta_{q+1}\delta_{q+1,n}
		\ell_q^{-2}\ell_q^{-2(n+1)\epsilon_0}\ell_q^{-5\alpha}.
		\]
		 For the second term, the identity  $\div\bar v_q=0$  and integration by parts give
		\[
		\left|
		\int_{\mathbb T^2}\sum_{k\in\Lambda}
		\bar v_q\div(f_{k,i,n}A_{k,i,n})\,\dd x
		\right|
		\lesssim
		\|\nabla\bar v_q\|_0
		\sum_{k\in\Lambda}
		\|f_{k,i,n}A_{k,i,n}\|_0\lesssim\tau_q^{-1}\ell_q^\alpha\delta_{q+1,n}.
		\]
		Summing over  $n$  and  $i$  proves  \eqref{transport-energy-time}. By using Proposition \ref{est-gluing}, a similar calculation yields \eqref{2transport-energy-time}. 
		
		By Corollary~\ref{est-wtq} and Hölder's inequality,
		\[
		\begin{aligned}
			\left|
			\int_{\mathbb T^2}
			2\wtq\cdot\bar v_q\,\dd x
			\right|=\left|
			\int_{\mathbb T^2}2\nabla \bar v_q:\mathcal{R} 
			\wtq\,\dd x
			\right|
			\lesssim \delta_{q}^{\frac{1}{2}}\lambda_{q}
			\mu_{q+1}^{-1}
			\delta_{q+1}
			\ell_q^{-5\alpha}\lesssim
			\ell_q^\alpha\delta_{q+1},
		\end{aligned}
		\]
		which proves \eqref{transport-energy}.
	\end{proof}
	
	 The preceding estimates for  $w_{q+1}^{(t)}$  immediately yield the following bounds for  $\bar{\Phi}_i$.
	
	\begin{prop}[Estimates for $\bar{\Phi}_i$]\label{p:estimates-for-inverse-bar-flow-map}
		 For  \(a \gg 1\),  \(0 \leq N \leq N_0+3\), and  \(t\in[t_i-\tau_q,t_i+2\tau_q]\cap[0,T]\), the following estimates hold:
		\begin{align}
			\|\nabla \bar{\Phi}_i^{\pm1} - {\rm Id}_{2 \times 2}\|_0 &\leq \frac{1}{10}, \label{e:nabla-barphi-i-minus-I3x3} \\
			\|(\nabla \bar{\Phi}_i)^{-1}\|_N + \|\nabla \bar{\Phi}_i\|_N
			&\lesssim \ell_q^{-N}\ell_q^{-LN\epsilon_0}, \label{e:nabla-barphi-i-CN} \\
			\|\matd{\vv_q + \wtq} \nabla \bar{\Phi}_i^{\pm1}\|_N
			&\lesssim \delta_q^{1/2} \lambda_q \ell_q^{-N}\ell_q^{-LN\epsilon_0}, \label{e:nabla-barphi-i-matd} \\
			\| \bar{\Phi}_i^{\pm1} - \Phi_i^{\pm1}\|_0
			&\lesssim \ell_q^{-1}\ell_q^{-\epsilon_0}\lambda_q^{-2} \Big(\frac{\lambda_q}{\lambda_{q+1}}\Big)^{\frac{1}{3}} \Big(\frac{\delta_{q+1}}{\delta_q}\Big)^{1/2}, \label{e:barphi-phi-C0}  \\
			\|\nabla \bar{\Phi}_i^{\pm1} - \nabla \Phi_i^{\pm1}\|_0
			&\lesssim \ell_q^{-2}\ell_q^{-2\epsilon_0}\lambda_q^{-2} \Big(\frac{\lambda_q}{\lambda_{q+1}}\Big)^{\frac{1}{3}} \Big(\frac{\delta_{q+1}}{\delta_q}\Big)^{1/2}. \label{e:barphi-phi-C1}
		\end{align}
	\end{prop}
	
	\begin{proof}
		It suffices to treat  \(\bar\Phi_i\), since  \(\bar\Phi_i^{-1}\) is estimated in the same way. The standard flow estimates,  \eqref{barphi}, and Corollary~\ref{est-wtq}  yield  \eqref{e:nabla-barphi-i-minus-I3x3}--\eqref{e:nabla-barphi-i-matd}. In particular,
		\begin{align*}
			\|\matd{\vv_q + \wtq} \nabla \bar{\Phi}_i\|_N
			&\lesssim \|\nabla(\vv_q+\wtq)\nabla\bar\Phi_i\|_N
				\lesssim \delta_q^{1/2}\lambda_q\ell_q^{-N}\ell_q^{-LN\epsilon_0}.
		\end{align*}
		 Subtracting the two flow equations gives
		\begin{equation*}
			\left\{
			\begin{alignedat}{2}
				\del_t (\bar{\Phi}_i - \Phi_i) + \vv_q \cdot \nabla (\bar{\Phi}_i - \Phi_i) + \wtq \cdot \nabla \bar{\Phi}_i &= 0, \\
				(\bar{\Phi}_i - \Phi_i)(t_i, x) &= 0.
			\end{alignedat}
			\right.
		\end{equation*}
		Gronwall's inequality along the  \(\vv_q\)-flow therefore gives
		\begin{align}\label{barphi-simple-diff-C0-raw}
			\|\bar\Phi_i-\Phi_i\|_0
				\lesssim \tau_q\mu_{q+1}^{-1}\delta_{q+1}\ell_q^{-1}\ell_q^{-\epsilon_0}\ell_q^{-5\alpha}.
		\end{align}
		 By the definitions of the parameters,
		\begin{equation}\label{barphi-simple-tau-mu-compute}
			\tau_q\mu_{q+1}^{-1}\delta_{q+1}
				\lesssim
				\lambda_q^{-2-10\alpha}
				\Big(\frac{\lambda_q}{\lambda_{q+1}}\Big)^{\frac{1}{3}}
				\Big(\frac{\delta_{q+1}}{\delta_q}\Big)^{1/2}.
		\end{equation}
		 This proves  \eqref{e:barphi-phi-C0}. Differentiating once in space and repeating the argument gives
		\begin{align}\label{barphi-simple-diff-grad-raw}
			\|\nabla\bar\Phi_i-\nabla\Phi_i\|_0
				\lesssim \tau_q\mu_{q+1}^{-1}\delta_{q+1}\ell_q^{-2}\ell_q^{-2\epsilon_0}\ell_q^{-5\alpha}.
		\end{align}
		Together with  \eqref{barphi-simple-tau-mu-compute}, this proves  \eqref{e:barphi-phi-C1}.
	\end{proof}
	
	Recall that  \(\bar{a}_{k,i,n}\)  and  $\bar{A}_{k,i,n}$  are obtained from  ${a}_{k,i,n}$ and  ${A}_{k,i,n}$, respectively, by replacing  ${\Phi}_i$ with  $\bar{\Phi}_i$. Propositions  \ref{est-RM}--\ref{p:estimates-for-inverse-bar-flow-map}  therefore yield the following estimates.
	
	\begin{cor}\label{barA-A}
		 For  \(0\leq n\leq L-1\),  \(0\leq N\leq N_0+3\), and  \(t\in[t_i-\tau_q,t_i+2\tau_q]\cap[0,T]\), we have
		\begin{align*}
			\|\bar{a}_{k,i,0}\|_N
			&\leq \frac{M}{8}\delta_{q+1}^{1/2}\ell_q^{-N} \ell_q^{-LN\epsilon_0} ,
			\quad
			\|\bar{a}_{k,i,n}\|_N
				\lesssim \delta_{q+1,n}^{1/2}\ell_q^{-N}\ell_q^{-LN\epsilon_0}; \\
			\|\matd{\vv_q + \wtq}\bar{a}_{k,i,0}\|_N
			&\leq \frac{M}{8}\tau_q^{-1}\delta_{q+1}^{1/2}\ell_q^{-N} \ell_q^{-LN\epsilon_0} \ell_q^{-\frac32\alpha};
			\\
			\|\matd{\vv_q + \wtq}\bar{a}_{k,i,n}\|_N
			&\lesssim \tau_q^{-1}\delta_{q+1,n}^{1/2}\ell_q^{-N}\ell_q^{-LN\epsilon_0}\ell_q^{-\frac32\alpha};\\
			\|A_{k,i,n}-\bar{A}_{k,i,n}\|_0
			&\lesssim \ell_q^{-2}\ell_q^{-2\epsilon_0}\lambda_q^{-2}
				\Big(\frac{\lambda_q}{\lambda_{q+1}}\Big)^{\frac{1}{3}}
			\Big(\frac{\delta_{q+1}}{\delta_q}\Big)^{1/2}\delta_{q+1,n}.
		\end{align*}
	\end{cor}


	 We now estimate the perturbation  $w_{q+1}$  using the bounds for its coefficients.
	
	\begin{prop}[Estimates for $w_{q+1}$]\label{estimate-wq+1}
		For \(0 \leq N \leq N_0+2\), there exists a constant $M$ such that
		\begin{align}
			&\|w^{(p_1)}_{q+1}\|_0 + \tfrac{1}{\lambda^N_{q+1}} \|w^{(p_1)}_{q+1}\|_N \leq \frac{M}{2} \delta^{1/2}_{q+1}, \label{estimate-wp}\\
			&\|w^{(p_2)}_{q+1}\|_0 + \tfrac{1}{\lambda^N_{q+1}} \|w^{(p_2)}_{q+1}\|_N \lesssim \delta^{1/2}_{q+1,1} \ll \delta^{1/2}_{q+1}, \label{estimate-wp2}\\
			&\|w^{(c)}_{q+1}\|_0 + \tfrac{1}{\lambda^N_{q+1}} \|w^{(c)}_{q+1}\|_N
				\lesssim \lambda_{q+1}^{-1}\delta_{q+1}^{1/2}\ell_q^{-1}\ell_q^{-L\epsilon_0}, \label{estimate-wc}\\
			&\|\wtq\|_N
				\lesssim \mu_{q+1}^{-1}\ell_q^{-N-1}\ell_q^{-5\alpha}
				\sum_{n=0}^{L-1}\delta_{q+1,n}\ell_q^{-(n+1)(N+1)\epsilon_0}
				\ll\lambda_{q+1}^N\delta_{q+1}^{1/2},\\
			&\|w_{q+1}\|_0 + \tfrac{1}{\lambda^N_{q+1}} \|w_{q+1}\|_N \leq \frac{2M}{3} \delta^{1/2}_{q+1}. \label{estimate-w}
		\end{align}
	\end{prop}
	
	\begin{rem}
		 Interpolating  \eqref{estimate-w}  and using  $\gamma < \beta < \tfrac{1}{3}$, we obtain
		\begin{align}\label{Cbeta}
			\|w_{q+1}\|_{C^{\gamma}} \lesssim (\lambda_{q+1} \delta^{1/2}_{q+1})^{\gamma} (\delta^{1/2}_{q+1})^{1-\gamma}
			\lesssim \lambda^{\gamma-\beta}_{q+1},
		\end{align}
		Hence $\{v_q\}_{q \geq 1}$ is a Cauchy sequence in  $C_tC^{\gamma}_x$.
	\end{rem}

	\noindent{\textbf{Step 3:\,\,Estimates for Reynolds stress}}
	
	 Set  $v_{q+1} = \vv_q + w_{q+1}$. Then  $(v_{q+1}, P_{q+1}, \RR_{q+1})$  satisfies  \eqref{e:subsol-B}  at level  $q+1$, namely,
	\begin{equation}
		\left\{
		\begin{alignedat}{2}
			&\del_t v_{q+1} + \div (v_{q+1} \otimes v_{q+1}) + \nabla P_{q+1} = \div \RR_{q+1}, \\
			&\div v_{q+1} = 0,
		\end{alignedat}
		\right. \label{e:q+1-B}
	\end{equation}
	where
	\begin{align*}
		&\div \RR_{q+1} \\
		=& \div \underbrace{\mathcal{R}(\del_t (\wpq + \wcq) + (\vv_q + \wtq) \cdot \nabla (\wpq + \wcq))}_{\Rtransport} \\
		&+ \div \underbrace{\mathcal{R}((\wpq + \wcq) \cdot \nabla (\vv_q + \wtq))}_{\Rnash} \\
		&+ \div \big(\wcq \mathring\otimes (\wpq+\wcq) + (\wpq+\wcq) \mathring\otimes \wcq - \wcq\mathring\otimes \wcq + \wtq \mathring\otimes\wtq\big)\\
		&+ \div \mathcal{R} \sum_{i;n} \big(\partial_t \bar{\eta}_i w_{i,n+1} + \bar{\eta}_i \partial_t w_{i,n+1}\big) \\
		&+ \div \big(\vv_q \otimes \wtq + \wtq \otimes \vv_q + \mathcal{R} \nabla P^{(v)}_{q+1} + \wpq \otimes \wpq + \RRR_q\big).
	\end{align*}
	 The pressure is defined by
	\begin{align*}
			P_{q+1}={}&\bar p_q+P_{q+1}^{(v)}-\frac12\operatorname{tr}\Big(
			w_{q+1}^{(c)}\otimes(w_{q+1}^{(p)}+w_{q+1}^{(c)})
			+(w_{q+1}^{(p)}+w_{q+1}^{(c)})\otimes w_{q+1}^{(c)}\\
			&\hspace{42mm}-w_{q+1}^{(c)}\otimes w_{q+1}^{(c)}
			+w_{q+1}^{(t)}\otimes w_{q+1}^{(t)}\Big).
	\end{align*}
	 The auxiliary pressure is given by
	\begin{align}\label{Pu1}
		P^{(v)}_{q+1} :=  P^{(t)}_{q+1} + \sum_{n;i;k \in \Lambda} \frac{\div}{\Delta}\div(f_{k,i,n}  A_{k,i,n})- h_q\sum_i\widetilde\chi_i^2 .
	\end{align}
	Finally, after subtracting its spatial mean if necessary, we normalize  \(P_{q+1}\)  by requiring 
	\(\int_{\mathbb T^2}P_{q+1}\,\dd x=0\). This does not affect  \eqref{e:q+1-B}.

	\begin{prop}[ Estimate for  $\Rtransport$]\label{proptrans}
		\begin{equation}\label{Trans}
			\|\Rtransport\|_0
				\lesssim \ell_q^{-2}\ell_q^{-(L+1)\epsilon_0}\lambda_q^{-2}
				\Big(\frac{\lambda_q}{\lambda_{q+1}}\Big)^{\frac{1}{3}}
				\Big(\frac{\delta_{q+1}}{\delta_q}\Big)^{1/2}\delta_{q+1}.
		\end{equation}
	\end{prop}
	
	\begin{proof}
		 Since  $\matd{\vv_q + \wtq} \phi_k(\bar{\Phi}_i) = 0$, we have
		\begin{align*}
			&\matd{\vv_q + \wtq} (\wpq + \wcq) \\
			=& \lambda^{-1}_{q+1} \sum_{n;i;k \in \Lambda} \matd{\vv_q + \wtq} \nabla^\perp \big(\bar{g}_{k,i,n} \bar{a}_{k,i,n} \psi(\lambda_{q+1} k^{\perp} \cdot \bar{\Phi}_i)\big) \\
			=& \lambda^{-1}_{q+1}\sum_{n;i;k\in\Lambda}\Big[\nabla^\perp\Big(\matd{\vv_q+\wtq}(\bar g_{k,i,n}\bar a_{k,i,n})\,\psi(\lambda_{q+1}k^\perp\cdot\bar\Phi_i)\Big)\notag\\
			&\hspace{32mm}+\big[\matd{\vv_q+\wtq},\nabla^\perp\big]\Big(\bar g_{k,i,n}\bar a_{k,i,n}\psi(\lambda_{q+1}k^\perp\cdot\bar\Phi_i)\Big)\Big].
		\end{align*}
		 To justify the commutator term, set  $D_{t,v}=\partial_t+v\cdot\nabla$  and define
			\[
				J:=
				\begin{pmatrix}
					0&-1\\
					1&0
				\end{pmatrix},
				\qquad
				J\nabla
				=
				\begin{pmatrix}
					-\partial_{x_2}\\
					\partial_{x_1}
				\end{pmatrix}
				=\nabla^\perp.
			\]
			One deduces that 
			\[
			[D_{t,v},\nabla^\perp]F
			=J[D_{t,v},\nabla]F
			=-J(\nabla v)^{\!T}\nabla F,
			\]
			So we have
			\[
			\|[D_{t,v},\nabla^\perp]F\|_N
			\lesssim \|\nabla v\|_0\|\nabla F\|_N
			+\|\nabla v\|_N\|\nabla F\|_0.
			\] 
			 Corollaries  \ref{-1jie}  and  \ref{barA-A}  now give
			\begin{align}\label{tranw}
				\|\Rtransport\|_0
				&\lesssim \lambda_{q+1}^{-1+\alpha}\mu_{q+1}\ell_q^{-(L-1)\epsilon_0}\delta_{q+1}^{1/2}\notag\\
				&\lesssim \ell_q^{-2}\ell_q^{-(L+1)\epsilon_0}\lambda_q^{-2}
				\Big(\frac{\lambda_q}{\lambda_{q+1}}\Big)^{\frac{1}{3}}
				\Big(\frac{\delta_{q+1}}{\delta_q}\Big)^{1/2}\delta_{q+1}.
			\end{align}
			 This proves Proposition  \ref{proptrans}.
		\end{proof}
		
		 We next estimate the Nash error.
		
		\begin{prop}[ Estimate for  $\Rnash$]\label{Nash}
			\begin{equation}\label{nash}
				\|\Rnash\|_0 \lesssim \frac{\delta^{1/2}_{q+1} \delta^{1/2}_q \lambda_q}{\lambda^{1-\alpha}_{q+1}}.
			\end{equation}
		\end{prop}
		
		\begin{proof}
			 The definition of  $\wpq + \wcq$  in  \eqref{wpq-wcq}, together with Corollaries  \ref{-1jie},  \ref{est-wtq}, and  \ref{barA-A}, yields
			\begin{align*}
				\|\Rnash\|_0 = \|\mathcal{R} [(\wpq + \wcq) \nabla (\vv_q + \wtq)]\|_0 \lesssim \frac{\delta^{1/2}_{q+1} \delta^{1/2}_q \lambda_q}{\lambda^{1-\alpha}_{q+1}}.
			\end{align*}
			 This proves Proposition  \ref{Nash}.
		\end{proof}
		 We conclude by estimating the oscillation error.
		\begin{prop}[ Estimate for  $\Rosc$]\label{propMosc}
			\begin{align*}
				\|\Rosc\|_0
					\lesssim \ell_q^{-2}\ell_q^{-(L+1)\epsilon_0}\lambda_q^{-2}
					\Big(\frac{\lambda_q}{\lambda_{q+1}}\Big)^{\frac{1}{3}}
					\Big(\frac{\delta_{q+1}}{\delta_q}\Big)^{1/2}\delta_{q+1}
					+\lambda_{q+1}^{-5\alpha}\delta_{q+2}.
			\end{align*}
		\end{prop}
		\begin{proof}
			 We first recall the decomposition of the Reynolds stress:
			\begin{align}\label{Rosc}
				&\div\Rosc\notag\\
				={}&\div\Big(\wcq\mathring\otimes(\wpq+\wcq)
					+(\wpq+\wcq)\mathring\otimes\wcq
					-\wcq\mathring\otimes\wcq
					+\wtq\mathring\otimes\wtq\Big)\notag\\
				&+\partial_t\wtq+\vv_q\cdot\nabla\wtq+\wtq\cdot\nabla\vv_q
				+\nabla P^{(v)}_{q+1}+\div(\wpq\otimes\wpq+\RRR_q)\notag\\
				={}&\div(R_{osc,1}+R_{osc,2}).
			\end{align}
			Here
			\begin{align*}
				\div R_{osc,1}
				={}&\div\Big(\wcq\mathring\otimes(\wpq+\wcq)
					+(\wpq+\wcq)\mathring\otimes\wcq
					-\wcq\mathring\otimes\wcq
					+\wtq\mathring\otimes\wtq\Big),
			\end{align*}
			and
			\begin{align*}
				\div R_{osc,2}
				={}&\div\mathcal R\Big(\partial_t\wtq+\vv_q\cdot\nabla\wtq
				+\wtq\cdot\nabla\vv_q+\nabla P^{(v)}_{q+1}
				+\div(\wpq\otimes\wpq+\RRR_q)\Big).
			\end{align*}
			
			By Corollary~\ref{est-wtq} and Proposition~\ref{estimate-wq+1}, we obtain
			\begin{align}\label{Rosc1}
				\|R_{osc,1}\|_0
				&\lesssim \|\wcq\|_0\big(\|\wpq\|_0+\|\wcq\|_0\big)+\|\wtq\|_0^2\notag\\
				&\lesssim \lambda_{q+1}^{-1}\ell_q^{-1}
				\ell_q^{-L\epsilon_0}\delta_{q+1}
				+\mu_{q+1}^{-2}\ell_q^{-2}
				\ell_q^{-2\epsilon_0}
				\ell_q^{-10\alpha}\delta_{q+1}^2\notag\\
				&\lesssim \ell_q^{-2}
				\ell_q^{-(L+1)\epsilon_0}\lambda_q^{-2}
				\Big(\frac{\lambda_q}{\lambda_{q+1}}\Big)^{\frac{1}{3}}
				\Big(\frac{\delta_{q+1}}{\delta_q}\Big)^{1/2}\delta_{q+1}.
			\end{align}
			
			 Next, we turn to estimate $R_{osc,2}$. Applying \eqref{bujiao}  and   \eqref{intg}   gives 
			\begin{align}\label{Pu}
				&\div (\RRR_q + \wpq \otimes \wpq) \notag\\
				=& \div \Big(\RRR_q + \sum_{n;i;k \in \Lambda}   {g}^2_{k,i,n} A_{k,i,n} + \sum_{n;i;k \in \Lambda}  \bar{g}^2_{k,i,n} \mathbb{P}_{\ne0}(\phi^2_k) (\bar{\Phi}_{i})\bar{A}_{k,i,n}+ \sum_{n;i;k \in \Lambda}   ( \bar{g}^2_{k,i,n}\bar{A}_{k,i,n} -   {g}^2_{k,i,n}A_{k,i,n})\Big) \notag\\
				=& \div \Big((\RRR_q +\sum_{n;i;k \in \Lambda}  A_{k,i,n} ) + \sum_{n;i;k \in \Lambda}\left[\left(\int_{\mathbb{T}}g^2_{k,n}(y)\dd y-1\right)-f_{k,i,n}\right] A_{k,i,n} \notag\\
				&+ \sum_{n;i;k \in \Lambda}  \bar{g}^2_{k,i,n} \mathbb{P}_{\ne0}(\phi^2_k)(\bar{\Phi}_{i}) \bar{A}_{k,i,n}+ \sum_{n;i;k \in \Lambda}   ( \bar{g}^2_{k,i,n}\bar{A}_{k,i,n} -   {g}^2_{k,i,n}A_{k,i,n})\Big) \notag\\
				=& \div \Big(\Big(h_q\sum_i\widetilde\chi_i^2\Big)\operatorname{Id}
				-  \sum_{n=0}^{L-2} \sum_{i} \partial_t \bar{\eta}_i \mathcal{R} w_{i,n+1} +\sum_{n;i;k \in \Lambda}  \bar{g}^2_{k,i,n} \mathbb{P}_{\ne0}(\phi^2_k) (\bar{\Phi}_{i})\bar{A}_{k,i,n} \notag\\
				&+ \sum_{n;i;k \in \Lambda}   ( \bar{g}^2_{k,i,n}\bar{A}_{k,i,n} -   {g}^2_{k,i,n}A_{k,i,n})
				+\sum_{n;i;k \in \Lambda}\left(\int_{\mathbb{T}}g^2_{k,n}(y)\dd y-1\right)A_{k,i,n}\Big) \notag\\
				&-\sum_{n;i;k \in \Lambda}\mathbb{P}_H\div(f_{k,i,n}A_{k,i,n})
				-\sum_{n;i;k \in \Lambda}\frac{\nabla\div}{\Delta}\div(f_{k,i,n}A_{k,i,n}).
			\end{align}
			 Here we used the identity
			\begin{align}
				&\div \sum_{n;i;k \in \Lambda} A_{k,i,n} \notag\\
				=& \div \sum_{n;i;k \in \Lambda} (a_{k,i,n})^2 \nabla \Phi^{-1}_i (k \otimes k) \nabla \Phi^\text{-T}_i \notag\\
				=& \nabla\Big(h_q\sum_i\widetilde\chi_i^2(t,x)+\sum_{n=1}^{L-1}\delta_{q+1,n}\sum_i\bar\chi_i^2(t)\Big)
				-\sum_{n=0}^{L-1} \div \RRR_{q,n} \notag\\
				=& \nabla\Big(h_q\sum_i\widetilde\chi_i^2\Big)
				-\div \RRR_q-\div\sum_{n=0}^{L-2}\sum_i\partial_t\bar\eta_i\mathcal Rw_{i,n+1}.
			\end{align}
			
			 Since
			\begin{align}\label{Pu1-osc}
				\nabla P^{(v)}_{q+1} := \nabla \Big(P^{(t)}_{q+1} + \sum_{n;i;k \in \Lambda} \frac{\div}{\Delta}\div(f_{k,i,n}  A_{k,i,n})
				-h_q\sum_i\widetilde\chi_i^2\Big),
			\end{align}
			 it follows that
			\begin{align}\label{Pu-osc2}
				&\div R_{osc,2}\notag\\
				={}&\partial_t \wtq+\vv_q\cdot\nabla\wtq+\wtq\cdot\nabla\vv_q+\nabla P^{(v)}_{q+1}
				+\div(\wpq\otimes\wpq+\RRR_q)\notag\\
				={}&\div\Big(-\sum_{n=0}^{L-2}\sum_i\partial_t\bar\eta_i\mathcal Rw_{i,n+1}
				+\sum_{n=0}^{L-1}\sum_i\partial_t\bar\eta_i\mathcal Rw_{i,n+1}\notag\\
				&\quad+\sum_{n;i;k\in\Lambda}\left(\int_{\mathbb T}g^2_{k,n}(y)\dd y-1\right)A_{k,i,n}\notag\\
				&\quad+\sum_{n;i;k\in\Lambda}\bar g^2_{k,i,n}\mathbb P_{\ne0}(\phi_k^2)(\bar{\Phi}_{i})\bar A_{k,i,n}
				+\sum_{n;i;k\in\Lambda}(\bar g^2_{k,i,n}\bar A_{k,i,n}-g^2_{k,i,n}A_{k,i,n})\Big)\notag\\
				={}&\div\Big(\sum_i\partial_t\bar\eta_i\mathcal Rw_{i,L}
				+\sum_{n;i;k\in\Lambda}\left(\int_{\mathbb T}g^2_{k,n}(y)\dd y-1\right)A_{k,i,n}\notag\\
				&\quad+\sum_{n;i;k\in\Lambda}\bar g^2_{k,i,n}\mathbb P_{\ne0}(\phi_k^2)(\bar{\Phi}_{i})\bar A_{k,i,n}
				+\sum_{n;i;k\in\Lambda}(\bar g^2_{k,i,n}\bar A_{k,i,n}-g^2_{k,i,n}A_{k,i,n})\Big)\notag\\
				={}&\div\Big[\mathcal{R}\div\Big(\sum_i\partial_t\bar\eta_i\mathcal Rw_{i,L}
				+\sum_{n;i;k\in\Lambda}\left(\int_{\mathbb T}g^2_{k,n}(y)\dd y-1\right)A_{k,i,n}\notag\\
				&\quad
				+\sum_{n;i;k\in\Lambda}(\bar g^2_{k,i,n}\bar A_{k,i,n}-g^2_{k,i,n}A_{k,i,n})\Big)+\mathcal{R}\big(\sum_{n;i;k\in\Lambda}\mathbb P_{\ne0}(\phi_k^2)(\bar{\Phi}_{i})\div(\bar g^2_{k,i,n}\bar A_{k,i,n})\big)\Big].
			\end{align}
			Moreover, by \eqref{zero-mode-final},
			\[
				\left\|\left(\int_{\mathbb T}g^2_{k,n}(y)\dd y-1\right)A_{k,i,n}\right\|_0
				\lesssim \ell_q^{N_0\epsilon_0}\delta_{q+1,n}.
			\]
			
			 For each fixed
				\(0\leq n\leq L-1\),
				\begin{align*}
					\bar g_{k,i,n}^2\bar A_{k,i,n}
					-g_{k,i,n}^2A_{k,i,n}
					={}&\bar g_{k,i,n}^2
					\big(\bar A_{k,i,n}-A_{k,i,n}\big)
					+\big(\bar g_{k,i,n}^2-g_{k,i,n}^2\big)A_{k,i,n}.
				\end{align*}
				The mean value theorem and \eqref{e:barphi-phi-C0} imply
				\begin{align*}
					\|\bar g_{k,i,n}-g_{k,i,n}\|_0
					&\lesssim
					\ell_q^{-1}\ell_q^{-(n+1)\epsilon_0}
					\|\bar\Phi_i^1-\Phi_i^1\|_0\notag\\
					&\lesssim
					\ell_q^{-2}\ell_q^{-(n+2)\epsilon_0}\lambda_q^{-2}
					\Big(\frac{\lambda_q}{\lambda_{q+1}}\Big)^{\frac{1}{3}}
					\Big(\frac{\delta_{q+1}}{\delta_q}\Big)^{1/2}\notag\\
					&\lesssim
					\ell_q^{-2}\ell_q^{-(L+1)\epsilon_0}\lambda_q^{-2}
					\Big(\frac{\lambda_q}{\lambda_{q+1}}\Big)^{\frac{1}{3}}
					\Big(\frac{\delta_{q+1}}{\delta_q}\Big)^{1/2}.
				\end{align*}
				Hence
				\begin{align*}
					\|\bar g_{k,i,n}^2-g_{k,i,n}^2\|_0
					&\leq
					\big(\|\bar g_{k,i,n}\|_0+\|g_{k,i,n}\|_0\big)
					\|\bar g_{k,i,n}-g_{k,i,n}\|_0\notag\\
					&\lesssim
					\ell_q^{-2}\ell_q^{-(L+1)\epsilon_0}\lambda_q^{-2}
					\Big(\frac{\lambda_q}{\lambda_{q+1}}\Big)^{\frac{1}{3}}
					\Big(\frac{\delta_{q+1}}{\delta_q}\Big)^{1/2}.
				\end{align*}
				Applying  Corollary~\ref{barA-A} and the zeroth-order estimate
				\(\|A_{k,i,n}\|_0\lesssim\delta_{q+1,n}\) established in the proof of
				Proposition~\ref{est-RM}, we conclude that
				\begin{align}
					&\|\bar g_{k,i,n}^2\bar A_{k,i,n}
					-g_{k,i,n}^2A_{k,i,n}\|_0\notag\\
					&\quad\leq
					\|\bar g_{k,i,n}\|_0^2
					\|\bar A_{k,i,n}-A_{k,i,n}\|_0
					+\|\bar g_{k,i,n}^2-g_{k,i,n}^2\|_0
					\|A_{k,i,n}\|_0\notag\\
					&\quad\lesssim
					\ell_q^{-2}\ell_q^{-2\epsilon_0}\lambda_q^{-2}
					\Big(\frac{\lambda_q}{\lambda_{q+1}}\Big)^{\frac{1}{3}}
					\Big(\frac{\delta_{q+1}}{\delta_q}\Big)^{1/2}\delta_{q+1,n}\notag\\
					&\qquad\quad+
					\ell_q^{-2}\ell_q^{-(L+1)\epsilon_0}\lambda_q^{-2}
					\Big(\frac{\lambda_q}{\lambda_{q+1}}\Big)^{\frac{1}{3}}
					\Big(\frac{\delta_{q+1}}{\delta_q}\Big)^{1/2}\delta_{q+1,n}\notag\\
					&\quad\lesssim
					\ell_q^{-2}\ell_q^{-(L+1)\epsilon_0}\lambda_q^{-2}
					\Big(\frac{\lambda_q}{\lambda_{q+1}}\Big)^{\frac{1}{3}}
					\Big(\frac{\delta_{q+1}}{\delta_q}\Big)^{1/2}\delta_{q+1,n}.
				\label{gA}
				\end{align}
			 The preceding estimate and Corollary  \ref{barA-A}  yield
			\begin{align}\label{Rocs2}
				&\|R_{osc,2}\|_0
				\notag\\
				\lesssim& \lambda_{q+1}^{-1+\alpha} \|\bar{g}^2_{k,i,n}\bar{A}_{k,i,n}\|_1
				+ \|\bar{g}^2_{k,i,n}\bar{A}_{k,i,n} -   {g}^2_{k,i,n}A_{k,i,n}\|_0\notag\\
				&\quad+\left(\left\|\sum_{i} \partial_t \bar{\eta}_i \mathcal{R} w_{i,L}\right\|_0
				+\left\|\left(\int_{\mathbb{T}}g^2_{k,n}(y)\dd y-1\right)A_{k,i,n}\right\|_0\right)\notag\\
				\lesssim& \lambda_{q+1}^{-1+\alpha}\ell_q^{-1}\ell_q^{-L\epsilon_0}\delta_{q+1}
				+\ell_q^{-2}\ell_q^{-(L+1)\epsilon_0}\lambda_q^{-2}
					\Big(\frac{\lambda_q}{\lambda_{q+1}}\Big)^{\frac{1}{3}}
					\Big(\frac{\delta_{q+1}}{\delta_q}\Big)^{1/2}\delta_{q+1}
				+\delta_{q+1,L}+\ell_q^{N_0\epsilon_0}\delta_{q+1}\notag\\
				\lesssim& \ell_q^{-2}\ell_q^{-(L+1)\epsilon_0}\lambda_q^{-2}
					\Big(\frac{\lambda_q}{\lambda_{q+1}}\Big)^{\frac{1}{3}}
					\Big(\frac{\delta_{q+1}}{\delta_q}\Big)^{1/2}\delta_{q+1}
				+\lambda_{q+1}^{-5\alpha}\delta_{q+2}.
			\end{align}
			
			 Combining  \eqref{Rosc1}  and  \eqref{Rocs2}, we conclude that
			\begin{align*}
				\|\Rosc\|_0
					\lesssim \ell_q^{-2}\ell_q^{-(L+1)\epsilon_0}\lambda_q^{-2}
					\Big(\frac{\lambda_q}{\lambda_{q+1}}\Big)^{\frac{1}{3}}
					\Big(\frac{\delta_{q+1}}{\delta_q}\Big)^{1/2}\delta_{q+1}
					+\lambda_{q+1}^{-5\alpha}\delta_{q+2}.
			\end{align*}
			 This proves Proposition~\ref{propMosc}.
		\end{proof}

		Combining Propositions~\ref{proptrans}--\ref{propMosc} and recalling the decomposition
			\(
			\RR_{q+1}=\Rtransport+\Rnash+\Rosc,
			\)
			we obtain
		\begin{align}
			\|\RR_{q+1}\|_0
			&\lesssim
				\ell_q^{-2}\ell_q^{-(L+1)\epsilon_0}\lambda_q^{-2}
				\Big(\frac{\lambda_q}{\lambda_{q+1}}\Big)^{\frac{1}{3}}
				\Big(\frac{\delta_{q+1}}{\delta_q}\Big)^{1/2}\delta_{q+1}
			\notag\\
			&\quad+
				\lambda_{q+1}^{\alpha}
				\frac{\delta_{q+1}^{1/2}\delta_q^{1/2}\lambda_q}{\lambda_{q+1}}
				+\lambda_{q+1}^{-5\alpha}\delta_{q+2}
			\notag\\
			&\lesssim \delta_{q+2}\lambda_{q+1}^{-4\alpha}.
			\label{Reynolds-C0-closure}
		\end{align}
		Here we used \eqref{lambdaN}  and the inequality
		$$\lambda_{q+1}^{\alpha}
		\frac{\delta_{q+1}^{1/2}\delta_q^{1/2}\lambda_q}{\lambda_{q+1}}
		\leq \ell_q^{-2}\ell_q^{-(L+1)\epsilon_0}\lambda_q^{-2}
		\Big(\frac{\lambda_q}{\lambda_{q+1}}\Big)^{\frac{1}{3}}
		\Big(\frac{\delta_{q+1}}{\delta_q}\Big)^{1/2}\delta_{q+1}.$$

		The same decomposition, combined with the corresponding product, commutator, and oscillatory anti-divergence estimates, yields the higher-order spatial and temporal bounds below.  Indeed, the following proposition is included solely to establish the differentiability of 
$(w_{q+1}, \RR_{q+1})$ in time and space. Since it is not essential to the proof that $\|\RR_{q+1}\|_0\leq \delta_{q+2}\lambda^{-4\alpha}_{q+1}$, we omit the detailed term-by-term calculations here.
		\begin{prop}[ Higher-order and time-derivative estimates]\label{bound-partialt}
			For $0 \leq N \leq N_0+2$, we have 
			\begin{align}
				&\|\RR_{q+1}\|_N+\lambda_{q+1}^{-5}\|\partial_t \RR_{q+1}\|_N \leq \lambda_{q+1}^{N-4\alpha} \delta_{q+2}, \quad \|\partial_t w^{(p_1)}_{q+1}\|_N \leq \frac{M}{2} \lambda^{N+1}_{q+1} \delta^{1/2}_{q+1}, \\
				&\|\partial_t w^{(p_2)}_{q+1}\|_N + \|\partial_t w^{(c)}_{q+1}\|_N + \|\partial_t w^{(t)}_{q+1}\|_N \leq \lambda^{N+1}_{q+1} \delta^{1/2}_{q+1}.
			\end{align}
		\end{prop}
		
		\vspace{2mm}
		Proposition \ref{estimate-wq+1}--Proposition \ref{bound-partialt} yield Proposition   \ref{assertion} directly. 
		\subsubsection{A Picard-type iteration}\label{3.5}
		 We initialize the scheme by setting  $h^{(0)}_q =10 \delta_{q+1}$, which satisfies  \eqref{hq}. Proposition  \ref{assertion}  then yields a perturbation  $(w^{(0)}_{q+1}, P^{(0)}_{q+1}, \RR^{(0)}_{q+1})$  with
		$$
		w^{(0)}_{q+1} :=w^{(p_1)}_{q+1,0} + w^{(p_2)}_{q+1,0} + w^{(c)}_{q+1,0} + w^{(t)}_{q+1,0}.
		$$
		
		 For  $m \in \mathbb{N}$, define
		\begin{align}\label{energy-gap}
				h_q^{(m+1)}(t)
				:=\frac{1}{X(t)}
				\Bigg(
				e(t)-\int_{\mathbb T^2}\big(|\bar v_q|^2+2\bar v_q\cdot w_{q+1,m}^{(t)}\big)\,\dd x
				-2\sum_{n=1}^{L-1}\sum_i\bar\chi_i^2c_{n,\xi}\delta_{q+1,n}
				-15\delta_{q+2}
				\Bigg),
		\end{align}
		where
		\(
		X(t):=2\sum_i\int_{\mathbb T^2}\widetilde\chi_i^2(t,x)\,\dd x,\ c_{n,\xi}=\int_{\mathbb T}\prod_{\gamma=0}^{n-1}
		h_\xi^2\!\left(\left\lceil\ell_q^{-\epsilon_0}\right\rceil^\gamma y\right)\,\dd y
		\).
      
        We now iterate Proposition   \ref{assertion}  through the following Picard scheme.
		
		\begin{prop}\label{repeatedly iteration proposition}
			 Let  $m \geq 0$, and $h^{(m)}_q$  be defined by  \eqref{energy-gap}, with  $h^{(0)}_q = 10 \delta_{q+1}$. Suppose that  $(w^{(m)}_{q+1}, \RR^{(m)}_{q+1})$, where
			$$
			w^{(m)}_{q+1} := w^{(p_1)}_{q+1,m} + w^{(p_2)}_{q+1,m} + w^{(c)}_{q+1,m} + w^{(t)}_{q+1,m},
			$$
			 satisfies the following estimates for  $0 \leq N \leq N_0+2$:
			\begin{align}
				& \| \RR^{(m)}_{q+1} \|_N + \lambda^{-5}_{q+1}\| \partial_t \RR^{(m)}_{q+1} \|_N \leq \lambda^{N-4\alpha}_{q+1} \delta_{q+2}, \label{energy-gap-m1} \\
				& \| w^{(p_1)}_{q+1,m} \|_0 + \lambda^{-N}_{q+1} \| w^{(p_1)}_{q+1,m} \|_N + \lambda^{-N-1}_{q+1} \| \partial_t w^{(p_1)}_{q+1,m} \|_N \leq \frac{M}{2} \delta_{q+1}^{1/2}, \label{energy-gap-m2} \\
				& \| w^{(p_2)}_{q+1,m} \|_0 + \| w^{(c)}_{q+1,m} \|_0 + \| w^{(t)}_{q+1,m} \|_0 \leq \delta^{1/2}_{q+1}, \label{energy-gap-m3} \\
				& \| w^{(p_2)}_{q+1,m} \|_N + \| w^{(c)}_{q+1,m} \|_N + \| w^{(t)}_{q+1,m} \|_N \leq \lambda^N_{q+1} \delta^{1/2}_{q+1}, \label{energy-gap-m4} \\
				& \| \partial_t w^{(p_2)}_{q+1,m} \|_N + \| \partial_t w^{(c)}_{q+1,m} \|_N + \| \partial_t w^{(t)}_{q+1,m} \|_N \leq \lambda^{N+1}_{q+1} \delta^{1/2}_{q+1}, \label{energy-partialt} \\
				& h^{(m)}_q(t) \in [ \delta_{q+1}, 100 \delta_{q+1}], \quad \|\partial_t h^{(m)}_q\|_{0}\lesssim\tau_q^{-1}\delta_{q+1}, \quad\|\partial_{tt} h^{(m)}_q\|_{0} \lesssim \lambda^5_{q} \delta_{q+1}.\label{energy-gap-m5}
			\end{align} 
 Then the proposition \ref{assertion}, applied with  $h^{(m+1)}_q$, yields a perturbation  $(w^{(m+1)}_{q+1}, \RR^{(m+1)}_{q+1})$  with
			$$
			w^{(m+1)}_{q+1} := w^{(p_1)}_{q+1,m+1} + w^{(p_2)}_{q+1,m+1} + w^{(c)}_{q+1,m+1} + w^{(t)}_{q+1,m+1},
			$$
			 and  \eqref{energy-gap-m1}--\eqref{energy-gap-m5}  hold at level  $m+1$.
		\end{prop}
		\begin{proof}
			 For  $m=0$, the bounds  \eqref{energy-gap-m1}--\eqref{energy-gap-m5}  follow directly from Proposition  \ref{assertion}  and the choice  $h^{(0)}_q = 10\delta_{q+1}$.
			
			 Assume inductively that  \eqref{energy-gap-m1}--\eqref{energy-gap-m5}  hold at level  $m$. It remains only to verify  \eqref{energy-gap-m5}  at level  $m+1$, since the argument proving Proposition  \ref{assertion}  in Section  \ref{MIS}  then yields  \eqref{energy-gap-m1}--\eqref{energy-partialt}.
			
			 Combining  \eqref{transport-energy}
			 with  \eqref{energy-gap-m4}  gives
			\begin{align}\label{fenzi1}
					\Bigg(
					e(t)-\int_{\mathbb T^2}\big(|\bar v_q|^2+2\bar v_q\cdot w_{q+1,m}^{(t)}\big)\,\dd x
					-2\sum_{n=1}^{L-1}\sum_i\bar\chi_i^2c_{n,\xi}\delta_{q+1,n}
					-15\delta_{q+2}
					\Bigg)
					\in[\tfrac{7}{2}\delta_{q+1},48\delta_{q+1}].
			\end{align}
			 Since
			\(
			\frac12\leq2\sum_i\int_{\mathbb T^2}\widetilde\chi_i^2\,\dd x\leq2,
			\)
			 it follows that
			\(
			h_q^{(m+1)}(t)\in[ \delta_{q+1},100\delta_{q+1}].
			\)
			
			 Differentiating  \eqref{energy-gap}  in time gives
			\begin{align}\label{partialt-hqm-start}
				\partial_t h_q^{(m+1)}
				={}&-\frac{\partial_tX(t)}{X^2(t)}
				\Bigg(
				e(t)-\int_{\mathbb T^2}\big(|\vv_q|^2+2\vv_q\cdot w_{q+1,m}^{(t)}\big)\,\dd x\notag\\
				&\hspace{28mm}-2\sum_{n=1}^{L-1}\sum_i\bar\chi_i^2c_{n,\xi}\delta_{q+1,n}
				-15\delta_{q+2}
				\Bigg)\notag\\
				&+\frac1{X(t)}\Bigg(
				\partial_t\Big(e(t)-\int_{\mathbb T^2}|\vv_q|^2\,\dd x\Big)
				-2\partial_t\int_{\mathbb T^2}\vv_q\cdot w_{q+1,m}^{(t)}\,\dd x\notag\\
				&\hspace{42mm}-2\partial_t\sum_{n=1}^{L-1}\sum_i\bar\chi_i^2c_{n,\xi}\delta_{q+1,n}
				\Bigg).
			\end{align}
			
			Making use of   \eqref{e:glued-energy-derivative}, \eqref{2te:matd-RRR_q}  and  \eqref{transport-energy-time}, \eqref{2transport-energy-time}, \eqref{fenzi1}
            , we deduce that
			\begin{equation}\label{hq-m-bounds}
					\delta_{q+1}\leq h_q^{(m+1)}(t)\leq 100\delta_{q+1},
					\qquad
					\|\partial_t h_q^{(m+1)}\|_0\lesssim\tau_q^{-1}\delta_{q+1}, \quad\|\partial_{tt} h^{(m+1)}_q\|_0 \lesssim \lambda^5_{q } \delta_{q+1}.
			\end{equation}
			 This proves  \eqref{energy-gap-m5}  at level  \(m+1\).
			
			 Applying Proposition  \ref{assertion}  with  $h^{(m+1)}_q$, we obtain the perturbation  $(w^{(m+1)}_{q+1},\allowbreak\ \RR^{(m+1)}_{q+1},\allowbreak\ P^{(m+1)}_{q+1})$  defined by
			$$
			w^{(m+1)}_{q+1} := w^{(p_1)}_{q+1,m+1} + w^{(p_2)}_{q+1,m+1} + w^{(c)}_{q+1,m+1} + w^{(t)}_{q+1,m+1},
			$$
			 which satisfies  \eqref{energy-gap-m1}--\eqref{energy-gap-m5}  at level  $m+1$.
			
			 This completes the proof of Proposition  \ref{repeatedly iteration proposition}.
		\end{proof}
		
		 Combining Propositions  \ref{assertion}  and  \ref{repeatedly iteration proposition}, for each fixed  $q$  we obtain a sequence  $\{(h^{(m+1)}_q,\allowbreak w^{(m+1)}_{q+1},\allowbreak P^{(m+1)}_{q+1},\allowbreak \RR^{(m+1)}_{q+1})\}_{m \geq 0}$  satisfying  \eqref{energy-gap-m1}--\eqref{energy-gap-m5}. The Ascoli--Arzel\`a theorem, together with the strict contraction property established in Appendix~\ref{app:picard-C0}, yields
		$$
		(h^{(m)}_q, w^{(m)}_{q+1}, P^{(m)}_{q+1}, \RR^{(m)}_{q+1}) \to (h_q, w_{q+1}, P_{q+1}, \RR_{q+1}) \quad \textup{in } L^\infty_{t,x}, \quad m \to \infty, 
		$$
		 where the limit satisfies  \eqref{e:subsol-B}  and  \eqref{e:vq-C0}--\eqref{e:RR_q-C0}  at level  $q+1$, while  \eqref{e:velocity-diff}  follows upon setting  $v_{q+1}=\bar{v}_{q}+w_{q+1}$. Passing to the limit in  \eqref{energy-gap}  gives
		\begin{align}
			h_q(t)\coloneq{}&\frac{1}{2\sum_i\int_{\mathbb{T}^2}\widetilde{\chi}^2_{i }\dd x }
			\Bigg(e(t)-\int_{\mathbb T^2}\Big(|\bar{v}_{q}|^2
			+2\bar{v}_{q}\cdot w^{(t)}_{q+1}\Big)\dd x\notag\\
			&\hspace{25mm}-2\sum_{n=1}^{L-1}\sum_i\bar{\chi}_i^2 c_{n,\xi}\delta_{q+1,n}
			-15\delta_{q+2}\Bigg).
			\label{limit-energy-gap-m6}
		\end{align}

		 It remains to verify  \eqref{e:initial}  at level  $q+1$, with  $h_q(t)$  given by  \eqref{limit-energy-gap-m6}.
		\begin{prop}[Energy estimate]\label{Pro-3.21}
			\label{p:energy} For all  $t\in[0,T]$, we have
			\[\bigg|e(t)-\int_{\mathbb T^2}|v_{q+1}|^2\dd x-15\delta_{q+2}\bigg|\leq\delta_{q+2}.\]
		\end{prop}
		\begin{proof}
			 Fix  $t\in[0,T]$. The energy error decomposes as
			\begin{align}\label{energy0}
				&e(t)-\int_{\mathbb{T}^2}|v_{q+1}|^2\dd x-15\delta_{q+2}\notag\\
				={}&\bigg[e(t)-\int_{\mathbb{T}^2}|\bar{v}_q|^2\dd x-15\delta_{q+2}\bigg]
				-\int_{\mathbb{T}^2}|w_{q+1}|^2\dd x
				-2\int_{\mathbb{T}^2}w_{q+1}\cdot\bar{v}_q\dd x\notag\\
				={}&\bigg[e(t)-\int_{\mathbb{T}^2}\Big(|\bar{v}_q|^2+2\bar{v}_q\cdot w^{(t)}_{q+1}\Big)\dd x-15\delta_{q+2}\bigg]\notag\\
				&-\int_{\mathbb{T}^2}|w^{(p_1)}_{q+1}|^2\dd x
					-\int_{\mathbb{T}^2}|w^{(p_2)}_{q+1}|^2\dd x-e_{\mathrm{low},0},
			\end{align}
		where we have used \eqref{bujiaobar}, which yields
				$w^{(p_1)}_{q+1}\cdot w^{(p_2)}_{q+1}=0$. Moreover,
			\begin{align}
				e_{\mathrm{low},0}:={}&2\int_{\mathbb{T}^2}\bigg[
				\big(w^{(p)}_{q+1}+w^{(c)}_{q+1}\big)\cdot\bar{v}_q
				+w^{(p)}_{q+1}\cdot w^{(t)}_{q+1}
				+w^{(c)}_{q+1}\cdot w_{q+1}
				-\frac12|w^{(c)}_{q+1}|^2
				+\frac12|w^{(t)}_{q+1}|^2\bigg]\dd x.
			\end{align}
			
			The perturbation $w^{(p_1)}_{q+1}$ corresponds to the level $n=0$. Since $g_{k,0}=h_k$ and $\int_{\mathbb T}h_k^2\,\dd y=1$, we have
			\begin{align*}
				g_{k,i,0}^2
				={}&1+\big(\mathbb P_{\ne0}h_k^2\big)
				\left(\mu_{q+1}t-\ell_q^{-1}
				\left\lceil\ell_q^{-\epsilon_0}\right\rceil\Phi_i^1\right),\\
				\bar g_{k,i,0}^2
				={}&1+\big(\mathbb P_{\ne0}h_k^2\big)
				\left(\mu_{q+1}t-\ell_q^{-1}
				\left\lceil\ell_q^{-\epsilon_0}\right\rceil\bar\Phi_i^1\right),\\
				\phi_k^2
				={}&1+\mathbb P_{\ne0}(\phi_k^2).
			\end{align*}
			Using these decompositions and adding and subtracting
				$g_{k,i,0}^2\tr(A_{k,i,0})$, we obtain
			\begin{align*}
				&\int_{\mathbb T^2}|w^{(p_1)}_{q+1}|^2\dd x\\
				={}&\sum_{i;k\in\Lambda}\int_{\mathbb T^2}
				g_{k,i,0}^2\tr(A_{k,i,0})\dd x\\
				&+\sum_{i;k\in\Lambda}\int_{\mathbb T^2}
				\Big(\bar g_{k,i,0}^2\tr(\bar A_{k,i,0})
				-g_{k,i,0}^2\tr(A_{k,i,0})\Big)\dd x\\
				&+\sum_{i;k\in\Lambda}\int_{\mathbb T^2}
				\bar g_{k,i,0}^2\tr(\bar A_{k,i,0})
				\big(\mathbb P_{\ne0}\phi_k^2\big)(\bar\Phi_i)\dd x.
			\end{align*}
			The level-zero part of the first term satisfies
			\begin{align*}
				&\sum_{i;k\in\Lambda}\int_{\mathbb T^2}
				g_{k,i,0}^2\tr(A_{k,i,0})\dd x\\
				={}&\sum_{i;k\in\Lambda}\int_{\mathbb T^2}
				\tr(A_{k,i,0})\dd x\\
				&+\sum_{i;k\in\Lambda}\int_{\mathbb T^2}
				\big(\mathbb P_{\ne0}h_k^2\big)
				\left(\mu_{q+1}t-\ell_q^{-1}
				\left\lceil\ell_q^{-\epsilon_0}\right\rceil\Phi_i^1\right)
				\tr(A_{k,i,0})\dd x\\
				={}&2h_q(t)\sum_i\int_{\mathbb T^2}\widetilde\chi_i^2\dd x\\
				&+\sum_{i;k\in\Lambda}\int_{\mathbb T^2}
				\big(\mathbb P_{\ne0}h_k^2\big)
				\left(\mu_{q+1}t-\ell_q^{-1}
				\left\lceil\ell_q^{-\epsilon_0}\right\rceil\Phi_i^1\right)
				\tr(A_{k,i,0})\dd x.
			\end{align*}
			Indeed, the geometric lemma and the trace-free property of $\RRR_q$ give
			\[
			\sum_{k\in\Lambda}A_{k,i,0}
			=\widetilde\chi_i^2\big(h_q(t)\operatorname{Id}-\RRR_q\big),
			\qquad
			\sum_{k\in\Lambda}\tr(A_{k,i,0})
			=2h_q(t)\widetilde\chi_i^2.
			\]
			Consequently,
			\begin{align}\label{energy2}
				\int_{\mathbb T^2}|w^{(p_1)}_{q+1}|^2\dd x
				=2h_q(t)\sum_i\int_{\mathbb T^2}
				\widetilde\chi_i^2\dd x+e_{\mathrm{low},1},
			\end{align}
			where
			\begin{align*}
				e_{\mathrm{low},1}:={}&\sum_{i;k\in\Lambda}\int_{\mathbb T^2}
				\Big(\bar g_{k,i,0}^2\tr(\bar A_{k,i,0})
				-g_{k,i,0}^2\tr(A_{k,i,0})\Big)\dd x\\
				&+\sum_{i;k\in\Lambda}\int_{\mathbb T^2}
				\bar g_{k,i,0}^2\tr(\bar A_{k,i,0})
				\big(\mathbb P_{\ne0}\phi_k^2\big)(\bar\Phi_i)\dd x\\
				&+\sum_{i;k\in\Lambda}\int_{\mathbb T^2}
				\big(\mathbb P_{\ne0}h_k^2\big)
				\left(\mu_{q+1}t-\ell_q^{-1}
				\left\lceil\ell_q^{-\epsilon_0}\right\rceil\Phi_i^1\right)
				\tr(A_{k,i,0})\dd x.
			\end{align*}

			For the higher levels $1\leq n\leq L-1$, the same decomposition gives
			\begin{align*}
				&\int_{\mathbb{T}^2}|w^{(p_2)}_{q+1}|^2\dd x\\
				={}&\sum_{n=1}^{L-1}\sum_{i;k\in\Lambda}\int_{\mathbb T^2}
				g_{k,i,n}^2\tr(A_{k,i,n})\dd x\\
				&+\sum_{n=1}^{L-1}\sum_{i;k\in\Lambda}\int_{\mathbb T^2}
				\Big(\bar g_{k,i,n}^2\tr(\bar A_{k,i,n})
				-g_{k,i,n}^2\tr(A_{k,i,n})\Big)\dd x\\
				&+\sum_{n=1}^{L-1}\sum_{i;k\in\Lambda}\int_{\mathbb T^2}
				\bar g_{k,i,n}^2\tr(\bar A_{k,i,n})
				\big(\mathbb P_{\ne0}\phi_k^2\big)(\bar\Phi_i)\dd x.
			\end{align*}
			Using the definition of $g_{k,i,n}$, the geometric lemma, the measure-preserving property of $\Phi_i$, and the fact that $\RRR_{q,n}$ is trace-free, we obtain
			\begin{align*}
				&\sum_{n=1}^{L-1}\sum_{i;k\in\Lambda}\int_{\mathbb T^2}
				g_{k,i,n}^2\tr(A_{k,i,n})\dd x\\
				={}&\sum_{n=1}^{L-1}\sum_{i;k\in\Lambda}\int_{\mathbb T^2}
				\Bigg[1+\big(\mathbb P_{\ne0}h_k^2\big)
				\left(\left\lceil\ell_q^{-\epsilon_0}\right\rceil^ny\right)\Bigg]
				\prod_{\gamma=0}^{n-1}h_\xi^2
				\left(\left\lceil\ell_q^{-\epsilon_0}\right\rceil^\gamma y\right)
				\tr(A_{k,i,n})\dd x\\
				={}&\sum_{n=1}^{L-1}\sum_i\bar\chi_i^2\delta_{q+1,n}
				\int_{\mathbb T^2}\prod_{\gamma=0}^{n-1}h_\xi^2
				\left(\left\lceil\ell_q^{-\epsilon_0}\right\rceil^\gamma y\right)
				\tr\left(\operatorname{Id}-\frac{\RRR_{q,n}}{\delta_{q+1,n}}\right)\dd x\\
				&+\sum_{n=1}^{L-1}\sum_{i;k\in\Lambda}\int_{\mathbb T^2}
				\big(\mathbb P_{\ne0}h_k^2\big)
				\left(\left\lceil\ell_q^{-\epsilon_0}\right\rceil^ny\right)
				\prod_{\gamma=0}^{n-1}h_\xi^2
				\left(\left\lceil\ell_q^{-\epsilon_0}\right\rceil^\gamma y\right)
				\tr(A_{k,i,n})\dd x\\
				={}&2\sum_{n=1}^{L-1}\sum_i\bar\chi_i^2c_{n,\xi}\delta_{q+1,n}\\
				&+\sum_{n=1}^{L-1}\sum_{i;k\in\Lambda}\int_{\mathbb T^2}
				\big(\mathbb P_{\ne0}h_k^2\big)
				\left(\left\lceil\ell_q^{-\epsilon_0}\right\rceil^ny\right)
				\prod_{\gamma=0}^{n-1}h_\xi^2
				\left(\left\lceil\ell_q^{-\epsilon_0}\right\rceil^\gamma y\right)
				\tr(A_{k,i,n})\dd x.
			\end{align*}
			Consequently,
			\begin{align*}
				\int_{\mathbb{T}^2}|w^{(p_2)}_{q+1}|^2\dd x
				=2\sum_{n=1}^{L-1}\sum_i\bar\chi_i^2c_{n,\xi}\delta_{q+1,n}+e_{\mathrm{low},2},
			\end{align*}
			where
			\begin{align*}
				e_{\mathrm{low},2}:={}&\sum_{n=1}^{L-1}\sum_{i;k\in\Lambda}\int_{\mathbb T^2}
				\Big(\bar g_{k,i,n}^2\tr(\bar A_{k,i,n})
				-g_{k,i,n}^2\tr(A_{k,i,n})\Big)\dd x\\
				&+\sum_{n=1}^{L-1}\sum_{i;k\in\Lambda}\int_{\mathbb T^2}
				\bar g_{k,i,n}^2\tr(\bar A_{k,i,n})
				\big(\mathbb P_{\ne0}\phi_k^2\big)(\bar\Phi_i)\dd x\\
				&+\sum_{n=1}^{L-1}\sum_{i;k\in\Lambda}\int_{\mathbb T^2}
				\big(\mathbb P_{\ne0}h_k^2\big)
				\left(\left\lceil\ell_q^{-\epsilon_0}\right\rceil^ny\right)
				\prod_{\gamma=0}^{n-1}h_\xi^2
				\left(\left\lceil\ell_q^{-\epsilon_0}\right\rceil^\gamma y\right)
				\tr(A_{k,i,n})\dd x.
			\end{align*}

			Recalling the definition of $h_q(t)$ in \eqref{limit-energy-gap-m6}, we conclude that
			\begin{align}
				&\int_{\mathbb{T}^2}|w^{(p_1)}_{q+1}|^2\dd x
				+\int_{\mathbb{T}^2}|w^{(p_2)}_{q+1}|^2\dd x\notag\\
				={}&\bigg[e(t)-\int_{\mathbb{T}^2}\Big(|\bar{v}_q|^2
				+2\bar{v}_q\cdot w^{(t)}_{q+1}\Big)\dd x-15\delta_{q+2}\bigg]
				+e_{\mathrm{low},1}+e_{\mathrm{low},2}.
				\label{e:energy-wpq+dpq}
			\end{align}
			
			Substituting \eqref{e:energy-wpq+dpq} into \eqref{energy0}, we obtain
			\begin{align}\label{energy4}
				e(t)-15\delta_{q+2}-\int_{\mathbb{T}^2}|v_{q+1}|^2\dd x
					=-e_{\mathrm{low},0}-e_{\mathrm{low},1}-e_{\mathrm{low},2},
					\quad\forall t\in[0,T].
			\end{align}
			It remains to estimate $e_{\mathrm{low},0}$, $e_{\mathrm{low},1}$, and $e_{\mathrm{low},2}$. Taking $n=0$ in \eqref{gA} and using
			\eqref{lambdaN}, we obtain
			\begin{align*}
				\Bigg|\sum_{i;k\in\Lambda}\int_{\mathbb{T}^2}
				\Big(\bar g_{k,i,0}^2\tr(\bar A_{k,i,0})
				-g_{k,i,0}^2\tr(A_{k,i,0})\Big)\dd x\Bigg|
				\lesssim\ell_q^{2\alpha}\delta_{q+2}.
			\end{align*}
			
			The bounds
			$\|A_{k,i,0}\|_N\lesssim\delta_{q+1}\ell_q^{-N-2\alpha}$ and
				$\|\bar A_{k,i,0}\|_N\lesssim\delta_{q+1}\ell_q^{-N}
				\ell_q^{-LN\epsilon_0}$, combined with the estimate
			\begin{align}\label{fenbu}
				\Big|\int_{\mathbb{T}^2} f \mathbb{P}_{\geq c} g \dd x\Big|
				=\Big|\int_{\mathbb{T}^2}|\nabla|^N f
				|\nabla|^{-N}\mathbb{P}_{\geq c} g\dd x\Big|
				\leq c^{-N}\|g\|_{L^2}\|f\|_{H^N},
			\end{align}
			yield
			\begin{align*}
				&\Bigg|\sum_{i;k\in\Lambda}\int_{\mathbb{T}^2}
				\big(\mathbb P_{\ne0}h_k^2\big)
				\left(\mu_{q+1}t-\ell_q^{-1}
				\left\lceil\ell_q^{-\epsilon_0}\right\rceil\Phi_i^1\right)
				\tr(A_{k,i,0})\dd x\Bigg|\\
				&\quad+\Bigg|\sum_{i;k\in\Lambda}\int_{\mathbb{T}^2}
				\bar g_{k,i,0}^2\tr(\bar A_{k,i,0})
				\big(\mathbb P_{\ne0}\phi_k^2\big)(\bar\Phi_i)\dd x\Bigg|
				\lesssim\ell_q^{2\alpha}\delta_{q+2}.
			\end{align*}
			Thus, we have
			\begin{equation}\label{e_low}
				|e_{\mathrm{low},1}|\lesssim\ell_q^{2\alpha}\delta_{q+2}.
			\end{equation}
			
			Arguing as in the estimate of $e_{\mathrm{low},1}$, we similarly obtain
				\begin{align}\label{small terms-elow}
					|e_{\mathrm{low},2}|+|e_{\mathrm{low},0}|
					\lesssim\ell_q^\alpha\delta_{q+2}.
			\end{align}
			
			Combining \eqref{energy4}, \eqref{e_low}, and \eqref{small terms-elow}, and increasing $a$ if necessary, we obtain
			\begin{align}\label{energy4-final}
				\Big|e(t)-15\delta_{q+2}-\int_{\mathbb{T}^2}|v_{q+1}|^2\dd x\Big|
				\leq{}&|e_{\mathrm{low},0}|+|e_{\mathrm{low},1}|+|e_{\mathrm{low},2}|\notag\\
				\lesssim{}&\ell_q^\alpha\delta_{q+2}
				\leq\delta_{q+2},\qquad\forall t\in[0,T].
			\end{align}
			
			This completes the proof of Proposition~\ref{p:energy}.
			
		\end{proof}

		\subsection{Proof of Proposition~\ref{p:main-prop}}\label{proof-main-prop}
		\begin{proof}
			It remains to verify that $(v_{q+1},P_{q+1},\RR_{q+1})$ satisfies \eqref{p_q} and \eqref{e:vq-C0}--\eqref{e:initial} with $q$ replaced by $q+1$.
			
			By \eqref{e:stability-vv_q-N} and \eqref{estimate-w}, after increasing $a$ so that the fixed implicit constants are absorbed by the margin in $M$, we obtain
				\[
				\|v_{q+1}\|_0
				\leq\|\bar v_q\|_0+\|w_{q+1}\|_0
				\leq M\Big(1+\sum_{i=1}^{q+1}\delta_i^{1/2}\Big).
				\]
				This proves \eqref{e:vq-C0} at level $q+1$.
			
			Likewise, \eqref{e:vv_q-bound} and \eqref{estimate-w} imply that, for $1 \leq N \leq N_0+2$,
			\begin{align*}
				\|v_{q+1}\|_N \leq \|\vv_q\|_N + \|w_{q+1}\|_N \leq M \delta_{q+1}^{1/2} \lambda_{q+1}^N.
			\end{align*}

			Furthermore, since $\beta<\frac{1}{3}$, Propositions~\ref{proptrans}--\ref{p:energy}, together with \eqref{lambdaN}, yield
			\begin{align*}
				&\|\RR_{q+1}\|_0 \leq \delta_{q+2} \lambda_{q+1}^{-4\alpha}, \\
				10\delta_{q+2} &\leq e(t) - \|v_{q+1}\|_{L^2(\TTT^2)}^2 \leq 20\delta_{q+2}.
			\end{align*}
			
			It remains to verify the increment estimate, which does not follow from the preceding $C^N$ bound alone. We write
				\[
				v_{q+1}-v_q=(\bar v_q-v_{\ell_q})+(v_{\ell_q}-v_q)+w_{q+1}.
				\]
				 \eqref{e:v_ell-vq-high}, \eqref{e:stability-vv_q}, and
					\eqref{estimate-w} give the required $C^0$ estimate. For one spatial derivative, \eqref{e:v_ell-vq-high}, Proposition~\ref{p:stability}, and interpolation give
				\[
				\|\bar v_q-v_q\|_1=o(\delta_{q+1}^{1/2}\lambda_{q+1}),
				\]
				whereas \eqref{estimate-w} yields
				$\|w_{q+1}\|_1\leq\frac{2M}{3}\delta_{q+1}^{1/2}\lambda_{q+1}$.
				After increasing $a$, the two low-frequency errors are absorbed by the remaining $\frac{M}{3}$ margin. Hence
				\[
				\|v_{q+1}-v_q\|_0+\lambda_{q+1}^{-1}\|v_{q+1}-v_q\|_1
				\leq M\delta_{q+1}^{1/2},
				\]
				which proves \eqref{e:velocity-diff}. Finally, subtracting the spatial mean of
					$P_{q+1}$ yields \eqref{p_q} without changing the equation.
			
			This completes the proof of Proposition~\ref{p:main-prop}.
		\end{proof}
		



	\subsection{Proof of Proposition~\ref{p:main-prop2}}\label{3.6}
		\begin{proof}The proof of Proposition~\ref{p:main-prop2} follows the argument in \cite{MNY}(Proposition 3.3), for completeness, we sketch the main steps below.
        Since $(v_q,\RR_q)=(\tv_q,\tRR_q)$ on $[0, 8\tau_q]$, it follows that
			$$
			(\vv_q, \RRR_q) = (\widetilde{\vv}_q, \widetilde{\RRR}_q) \quad \text{and} \quad \Phi_i = \widetilde{\Phi}_i, \quad t \in [0, 8\tau_q].
			$$
			This implies that
			$$
			A^{(0)}_{k,i,n} = \widetilde{A}^{(0)}_{k,i,n}, \quad t \in [0,   5\tau_q],
			$$
			so that
			$$
			{w^{(t)}_{q+1,0} = \widetilde{w}^{(t)}_{q+1,0}}, \quad \bar{\Phi}_{i,0} = \widetilde{\bar{\Phi}}_{i,0}, \quad
            h_q^{(0)}=\widetilde{h}_q^{(0)},\quad t \in [0,  5\tau_q].
			$$
            
			Since
			\begin{align}\label{wpq-wcq1}
				{w^{(p)}_{q+1,0} + w^{(c)}_{q+1,0}}
				= & {\lambda^{-1}_{q+1} \nabla^\perp \sum_{i; k \in \Lambda} \bar g_{k,i,0,0} \bar{a}_{k,i,0,0} \psi(\lambda_{q+1} k^\perp \cdot \bar{\Phi}_{i,0})} \notag \\
				& {+ \lambda^{-1}_{q+1} \nabla^\perp \sum_{n=1}^{L-1} \sum_{i; k \in \Lambda} \bar g_{k,i,n,0} \bar{a}_{k,i,n,0} \psi(\lambda_{q+1} k^\perp \cdot \bar{\Phi}_{i,0}),}
			\end{align}
			{where
				\[
				\bar a_{k,i,0,0}=\widetilde\chi_i\big(h_q^{(0)}\big)^{1/2}
				a_k\!\left(\nabla\bar\Phi_{i,0}
				\left(\operatorname{Id}-\frac{\RRR_{q,0}^{(0)}}{h_q^{(0)}}\right)
				\nabla\bar\Phi_{i,0}^{T}\right),
				\]
				while, for $1\le n\le L-1$,
				\[
				\bar a_{k,i,n,0}=\bar\chi_i\delta_{q+1,n}^{1/2}
				a_k\!\left(\nabla\bar\Phi_{i,0}
				\left(\operatorname{Id}-\frac{\RRR_{q,n}^{(0)}}{\delta_{q+1,n}}\right)
				\nabla\bar\Phi_{i,0}^{T}\right),
				\]
                we deduce that
$$({ {w}^{(p)}_{q+1,0} ,  {w}^{(c)}_{q+1,0}})=({\widetilde{w}^{(p)}_{q+1,0} , \widetilde{w}^{(c)}_{q+1,0}}),\quad t\in [0,5\tau_q].$$

Next, recalling the definitions of $h_q^{(m)}, ~ w^{(m)}_{q+1}$ and $\RR_{q+1}^{(m)}$ with $m=0$, we obtain, for   $t\in[0,5\tau_{q }]$,
			\begin{align}\label{xiandeng0}
			(h^{(0)}_q, w^{(0)}_{q+1}, \RR^{(0)}_{q+1}) = (\widetilde{h}^{(0)}_q, \widetilde{w}^{(0)}_{q+1}, \tRR^{(0)}_{q+1}) \Longrightarrow (v^{(0)}_{q+1}, P^{(0)}_{q+1}, \RR^{(0)}_{q+1}) = (\tv^{(0)}_{q+1}, \widetilde{P}^{(0)}_{q+1}, \tRR^{(0)}_{q+1}).
			\end{align}
            
            Finally, using the   Picard-type iteration as in Proposition \ref{repeatedly iteration proposition}, one can easily deduce that for $t\in[0,5\tau_{q }]$,
            \begin{align}\label{xiandeng1}
			(h^{(m)}_q, w^{(m)}_{q+1}, \RR^{(m)}_{q+1}) = (\widetilde{h}^{(m)}_q, \widetilde{w}^{(m)}_{q+1}, \tRR^{(m)}_{q+1}) \Longrightarrow (v^{(m)}_{q+1}, P^{(m)}_{q+1}, \RR^{(m)}_{q+1}) = (\tv^{(m)}_{q+1}, \widetilde{P}^{(m)}_{q+1}, \tRR^{(m)}_{q+1}),\quad m\geq1.
			\end{align}
			Passing to the limit as $m\to\infty$, we obtain
		$$
			(h_q, w_{q+1}, \RR_{q+1}) = (\widetilde{h}_q, \widetilde{w}_{q+1}, \tRR_{q+1}) \Longrightarrow (v_{q+1}, P_{q+1}, \RR_{q+1}) = (\tv_{q+1}, \widetilde{P}_{q+1}, \tRR_{q+1}),\quad t\in[0,8\tau_{q+1 }].
			$$
            
			This completes the proof of Proposition~\ref{p:main-prop2}.}
			\end{proof}
			
			\vspace{2mm}
			
			In the next section, we prove Theorem~\ref{main1}: the wild initial data associated with the dissipative weak solutions constructed in Theorem~\ref{main0} are dense, in the $C^{\gamma'}$ topology, in the divergence-free subspace of $C^\gamma(\mathbb T^2)$.
		\section{Density of ``wild'' initial data}\label{Proof of Theorem 2}
		\begin{proof}
			Fix $T<\infty$ and $0<\gamma'<\gamma<\frac13$, and let $v^{\mathrm{in}}\in C^\gamma(\mathbb T^2;\mathbb R^2)$ be divergence-free. Given $\epsilon>0$, the standard properties of mollification allow us to choose $\delta=\delta(\epsilon)>0$ sufficiently small that
			$$
			\|v^{\text{in}} * \psi_{\delta} - v^{\text{in}}\|_{C^{\gamma'}} \leq \epsilon.
			$$
			
			We next solve the following Cauchy problem for the two-dimensional Euler equations:
			\begin{equation}
				\left\{
				\begin{alignedat}{-1}
					&\partial_t v_{\delta} + \div (v_{\delta} \otimes v_{\delta}) + \nabla P_{\delta} = 0, \\
					&\nabla \cdot v_{\delta} = 0, \\
					&v_{\delta}\big|_{t=0} = v^{\text{in}} * \psi_{\delta}.
				\end{alignedat}
				\right. \label{2deuler}
			\end{equation}
			For $t\in[0,T]$, the resulting smooth solution satisfies
			$$
			\|v_{\delta}\|_{L^{\infty}_T C^N} \leq C_N e^{5e^{T\delta^{-2}}} \delta^{-N} \|v^{\text{in}}\|_{0}, \quad N \geq 0.
			$$
			
			Choose $\gamma<\beta<\frac13$, let $M$ satisfy Proposition \ref{p:main-prop} and additionally be greater than $2C_0 e^{10e^{T\delta^{-2}}}  \|v^{\text{in}}\|_{C^{\gamma}} $,
and then take the parameter $a$ sufficiently large, depending on
$M$, $\alpha$, $T$, $\epsilon$, and $v^{\mathrm{in}}$,
as prescribed in Section~\ref{Iteration proposition}. More precisely,
			$$
			a \geq \max\big\{10^{\alpha^{-1}},~10T^{-1},~M^{10},~
			(2M\epsilon^{-1})^{(b^2 (\beta-\gamma'))^{-1}},~C^3_{N_0+2}e^{30e^{T\delta^{-2}}}\delta^{-3(N_0+2)}\|v^{\text{in}}\|^3_{0}\big\}.
			$$
			We take $\lambda_q=\lambda_\ast\lceil a^{b^q}\rceil$ and choose all remaining parameters as in Section~\ref{Iteration proposition}.
			
			Set $(v_1,P_1,\RR_1):=(v_\delta,P_\delta,0)$. We may then choose a strictly decreasing smooth positive function $e(t)$ so that
			\begin{align}
				&\|v_1\|_{0} \leq M\big(1+ \delta^{1/2}_1\big), \label{e:vq-C01} \\
				&\|v_1\|_{N} \leq M \delta_{1}^{1/2} \lambda^N_1, \quad 1 \leq N \leq N_0+2, \label{e:vq-C11} \\
				&{\|\RR_1\|_{0} \leq \delta_{2} \lambda_1^{-4\alpha},} \label{e:RR_q-C01} \\
				&10\delta_{2} \leq e(t) - \|v_1(t)\|_{L^2(\TTT^2)}^2 \leq 20\delta_{2}. \label{e:initial1-start}
			\end{align}
			
			Applying Proposition~\ref{p:main-prop} with this initial triple on $[0,T]$, we obtain a weak solution $v_\delta^{(\mathrm{conv})}$ with strictly decreasing kinetic energy. Moreover, \eqref{Cbeta} with $\beta'=\gamma'<\beta$ gives
			$$
			\|(v^{(\text{conv})}_{\delta} - v_{\delta})(0,\cdot)\|_{C^{\gamma'}} \leq \|v^{(\text{conv})}_{\delta} - v_1\|_{C^{\gamma'}} \lesssim \sum_{q=2}^{\infty} \lambda^{\gamma'-\beta}_q \leq 2M \lambda^{\gamma'-\beta}_2 \leq \epsilon,
			$$
			where we used $a\geq(2M\epsilon^{-1})^{(b^2(\beta-\gamma'))^{-1}}$.

			Finally, since $v_\delta(0,x)=v^{\text{in}}*\psi_\delta$, we take $v_\delta^{(\mathrm{conv})}(0,x)$ as the wild initial datum and obtain
			$$
			\|v^{(\text{conv})}_{\delta}(0,\cdot) - v^{\text{in}}\|_{C^{\gamma'}} \leq {\|v^{\text{in}} * \psi_{\delta} - v^{\text{in}}\|_{C^{\gamma'}}} + \|v^{(\text{conv})}_{\delta}(0,\cdot) - v^{\text{in}} * \psi_{\delta}\|_{C^{\gamma'}} \leq 2\epsilon.
			$$
          Because of $\gamma<\beta$, interpolation in
\eqref{e:velocity-diff} also gives
\[
\|v_\delta^{(\mathrm{conv})}-v_\delta\|_{C_tC_x^\gamma}
\lesssim
M\sum_{q=2}^{\infty}\lambda_q^{\gamma-\beta}
<\infty.
\]
Thus $v_\delta^{(\mathrm{conv})}\in C_tC_x^\gamma$,
and its divergence-free initial trace belongs to
$C^\gamma(\mathbb T^2)$.
            
			By the freedom in choosing $e(t)$, we may select infinitely many
distinct admissible strictly decreasing profiles with uniform
smooth bounds that agree on $[0,10\tau_1]$, so
Proposition~\ref{p:main-prop2}, applied with the same parameters
and auxiliary choices, yields infinitely many distinct solutions
with the same initial datum.
Together with the $C^\gamma$ regularity and divergence-free
property established above, this shows that the common initial
datum is wild in the sense of Definition~\ref{WID}.
			
			This completes the proof of Theorem~\ref{main1}.
		\end{proof}
		
		\begin{appendices}
			\begingroup
			\section{\texorpdfstring{$L_{t}^\infty$}{L-infinity} contraction estimate for the Picard iteration}
			\label{app:picard-C0}
			
			Fix $q\geq1$. We prove that the map
			\begin{align}\label{energy-gap1}
				\Gamma[h]
				:=& \frac{1}{2\sum_i\int_{\mathbb T^2}\widetilde\chi_i^2\,\dd x}
				\Bigg(
				e(t)-\int_{\mathbb T^2}\big(|\bar v_q|^2+2\bar v_q\cdot w_{q+1}^{(t)}[h] \big)\,\dd x\notag\\
				&-2\sum_{n=1}^{L-1}\sum_i\bar\chi_i^2c_{n,\xi}\delta_{q+1,n}
				-15\delta_{q+2}
				\Bigg)
			\end{align}
			is a contraction. In particular, \eqref{hq-m-bounds} gives
			\[
			\delta_{q+1}\leq h(t) \leq100\delta_{q+1}.
			\]
			\begin{prop}[Contraction of the map in \eqref{energy-gap1}]
				\label{app:picard-C0-contraction}
				If $a$ is sufficiently large that $\lambda_q^{-\alpha}\ll1$, then
				\begin{align}
					&\| \Gamma[h]-\Gamma[\widetilde h]\|_{L^{\infty}_t}
					\leq \lambda^{-\alpha}_q\|h-\widetilde h\|_{L_{t}^\infty},\qquad t\in [0,T].         \label{app:full-mapcontraction}
				\end{align}
				It suffices to prove
				\begin{align}
					\left\|\int_{\Ttwo}\bar v_q\cdot
					\bigl(w_{q+1}^{(t)}[h]-w_{q+1}^{(t)}[\widetilde h]\bigr)\,\dd x\right\|_{L_{t}^\infty}
					\leq\lambda^{-\alpha}_q\|h-\widetilde h\|_{L_{t}^\infty}.                                     \label{app:picard-C0-final}
				\end{align}
			\end{prop}
			
			\begin{proof}
				Upon subtracting the two right-hand sides of \eqref{energy-gap1}, all terms cancel except for the Newton interaction. Since
				\[
				w_{q+1}^{(t)}[h]-w_{q+1}^{(t)}[\widetilde h]
				= \sum_{i=0}^{{i_{max}}-1}\sum_{n=0}^{L-1}\bar\eta_i\bigl(w_{i,n+1}[h]-w_{i,n+1}[\widetilde h]\bigr)
				\]
				and $\operatorname{div}\mathcal Ru=u$ for mean-free fields, periodic integration by parts and the bounded overlap of the cutoffs give
				\begin{align}
					&\left\|\int_{\Ttwo}\bar v_q\cdot
					\bigl(w_{q+1}^{(t)}[h]-w_{q+1}^{(t)}[\widetilde h]\bigr)\,\dd x\right\|_{L_{t}^\infty}
					\notag\\[-1mm]
					&\quad\leq L\|\nabla\bar v_q\|_{L_t^\infty L_x^2}
					\sup_i\sup_{\substack{\,0\leq n\leq L-1\\t\in[t_{i}-2\tau_{q},t_{i+1}+2\tau_{q}]}}
					\left\|\mathcal R\bigl(w_{i,n+1}[h]-w_{i,n+1}[\widetilde h]\bigr)(t)\right\|_{L_x^2}.       \label{app:target-reduction}
				\end{align}
				It therefore remains to prove that
				\begin{align}
					\sup_{\substack{\,0\leq n\leq L-1\\t\in[t_{i}-2\tau_{q},t_{i+1}+2\tau_{q}]}}
					\left\|\mathcal R\bigl(w_{i,n+1}[h]-w_{i,n+1}[\widetilde h]\bigr)(t)\right\|_{L_x^2}
					\leq C_{0}^{L}\tau_q\|h-\widetilde h\|_{L^{\infty}_t},\quad 0\leq i\leq i_{max}-1.                                  \label{app:needed-Rw-bound}
				\end{align}
				Here the sufficiently large constant $C_0$ is independent of $q$, $a$, $h$, $\widetilde h$, and $n$.
				
				Set $z=w_{i,n+1}[h]-w_{i,n+1}[\widetilde h]$, $\pi=p_{i,n+1}[h]-p_{i,n+1}[\widetilde h]$, and
				$F=\sum_{k\in\Lambda}f_{k,i,n}(A_{k,i,n}[h]-A_{k,i,n}[\widetilde h])$.
				Subtracting the two Newton equations yields
				\[
				\partial_tz+\bar v_q\cdot\nabla z+z\cdot\nabla\bar v_q+\nabla\pi
				=\mathbb P_H\operatorname{div}F,
				\qquad z(t_i)=0,
				\qquad \operatorname{div}z=0.
				\]
				To establish the required $H^{-1}$ stability estimate, let $\psi=(-\Delta)^{-1}z$ and
					$Y=\|\nabla\psi\|_{L_{x}^{2}}=\|z\|_{\dot H^{-1}_{x}}$. Periodic integration by parts, the identities $\operatorname{div}\bar v_q=\operatorname{div}\psi=0$, and the self-adjointness of $\mathbb P_H$ give
				\[
				\frac12\frac{\dd}{\dd t}Y^2
				\leq C\|\nabla\bar v_q\|_{L^{\infty}_{x}}Y^2+\|F\|_{L^{2}_{x}}Y,
				\qquad
				\frac{\dd}{\dd t}Y \leq C\|\nabla\bar v_q\|_{L^{\infty}_{x}}Y+\|F\|_{L^{2}_{x}}.
				\]
				Since $Y(t_i)=0$, Gronwall's inequality and $\|\mathcal Rz\|_{L_{x}^{2}}\lesssim\|z\|_{\dot H^{-1}_{x}}=Y$ yield, for $t\geq t_i$,
				\begin{align}
					\|\mathcal Rz(t)\|_{L_{x}^{2}}\lesssim Y(t)
					&\leq C \int_{t_i}^t
					\exp\!\left(C\int_s^t\|\nabla\bar v_q(r)\|_{L_{x}^{\infty}}\,\dd r\right)\|F(s)\|_{L_{x}^{2}}\,\dd s\notag\\
					&\leq C\tau_q\sup_{t}\sum\limits_{k\in\Lambda} \|A_{k,i,n}[h]-A_{k,i,n}[\widetilde h]\|_{L_{x}^{2}},                 \label{app:Hminus1-step}
				\end{align}
				where we used $|t-t_i|\leq4\tau_q$ and $\tau_q\|\nabla\bar v_q\|_{L_{x}^{\infty}}\lesssim\lambda_q^{-5\alpha}\ll1$.
				
				We now prove \eqref{app:needed-Rw-bound} by induction. More precisely, we claim that
				\begin{align}\label{n-step}
					\sup_{t}\left\|\mathcal R\bigl(w_{i,n+1}[h]-w_{i,n+1}[\widetilde h]\bigr)\right\|_{L_{x}^{2}}
					\leq C_0^{n+1}\tau_q\|h-\widetilde h\|_{L_{t}^\infty},\quad 0\leq n\leq L-1,~ 0\leq i\leq i_{max}-1.
				\end{align}
				For $n=0$,
				\[
				A_{k,i,0}[h]=\widetilde\chi_i^2h\,
				a_k^2\!\left(\nabla\Phi_i\left(\mathrm{Id}-\frac{\RRR_q}{h}\right)\nabla\Phi_i^{\mathsf T}\right)
				\nabla\Phi_i^{-1}(k\otimes k)\nabla\Phi_i^{-\mathsf T}.
				\]
				The arguments of the geometric coefficients remain in a fixed compact set, and the flow matrices are uniformly bounded. Hence the mean-value theorem gives
				\begin{align}
					\|A_{k,i,0}[h]-A_{k,i,0}[\widetilde h]\|_{L_x^2}
					&\lesssim    \|h-\widetilde h\|_{L_{t}^\infty}
					+\widetilde h\left\|\frac{\RRR_q}{h}-\frac{\RRR_q}{\widetilde h}\right\|_{L_x^2} \notag\\
					&\lesssim\|h-\widetilde h\|_{L_{t}^\infty}
					+\delta_{q+1}\|\RRR_q\|_{L_x^2}\frac{\|h-\widetilde h\|_{L_{t}^\infty}}{\delta_{q+1}^2}
					\leq C\|h-\widetilde h\|_{L_{t}^\infty}.                                          \label{app:A0-L2}
				\end{align}
				Substituting \eqref{app:A0-L2} into \eqref{app:Hminus1-step}, we obtain
				\begin{align*}
					\sup_{t}\left\|\mathcal R\bigl(w_{i,1}[h]-w_{i,1}[\widetilde h]\bigr)\right\|_{L_{x}^{2}}\leq C^2\tau_q\|h-\widetilde h\|_{L_{t}^\infty}
					\leq C_0\tau_q\|h-\widetilde h\|_{L_{t}^\infty}.
				\end{align*}
				Thus \eqref{n-step} holds for $n=0$, provided that the fixed constant $C_0\gg C^3$ is chosen sufficiently large.

				\smallskip
				
				Assume that \eqref{n-step} holds for some $n\geq0$. We must prove the corresponding estimate at level $n+1$, namely,
				\begin{align*}
					\sup_{t}\left\|\mathcal R\bigl(w_{i,n+2}[h]-w_{i,n+2}[\widetilde h]\bigr)\right\|_{L_{x}^{2}}
					&\leq C^{n+2}_0\tau_q \|h-\widetilde h\|_{L_{t}^\infty}.
				\end{align*}

				By \eqref{app:Hminus1-step},
				\begin{align*}
					\sup_{t}\left\|\mathcal R\bigl(w_{i,n+2}[h]-w_{i,n+2}[\widetilde h]\bigr)\right\|_{L_{x}^{2}}
					&\leq C\tau_q\sup_{t}
					\|A_{k,i,n+1}[h]-A_{k,i,n+1}[\widetilde h]\|_{L_{x}^{2}}.
				\end{align*}
				For $n\geq0$,
				\begin{align*}
					A_{k,i,n+1}[h]
					={}&\bar\chi_i^2\delta_{q+1,n+1}
					a_k^2\!\left(
					\nabla\Phi_i\left(\mathrm{Id}-\frac{\RRR_{q,n+1}[h]}{\delta_{q+1,n+1}}\right)
					\nabla\Phi_i^{\mathsf T}\right)
					\nabla\Phi_i^{-1}(k\otimes k)\nabla\Phi_i^{-\mathsf T},
				\end{align*}
				and the analogous formula holds with $h$ replaced by $\widetilde h$. Since
				$\bar\chi_i$, $\delta_{q+1,n+1}$, and $\Phi_i$ are independent of the
					candidate corrector, subtraction yields
				\begin{align*}
					&A_{k,i,n+1}[h]-A_{k,i,n+1}[\widetilde h]
					\\
					={}&\bar\chi_i^2\delta_{q+1,n+1}
					\Bigg[
					a_k^2\!\left(
					\nabla\Phi_i\left(\mathrm{Id}-\frac{\RRR_{q,n+1}[h]}{\delta_{q+1,n+1}}\right)
					\nabla\Phi_i^{\mathsf T}\right)
					\\[-1mm]
					&\hspace{27mm}-
					a_k^2\!\left(
					\nabla\Phi_i\left(\mathrm{Id}-\frac{\RRR_{q,n+1}[\widetilde h]}{\delta_{q+1,n+1}}\right)
					\nabla\Phi_i^{\mathsf T}\right)
					\Bigg]
					\nabla\Phi_i^{-1}(k\otimes k)\nabla\Phi_i^{-\mathsf T}.
				\end{align*}
				Using the integral form of the mean-value theorem and the uniform bounds for
					$D(a_k^2)$ on the preceding compact set, we obtain
				\begin{align*}
					&\|A_{k,i,n+1}[h]-A_{k,i,n+1}[\widetilde h]\|_{L_x^2}
					\\
					&\quad\leq C
					\delta_{q+1,n+1}\delta_{q+1,n+1}^{-1}
					\|\nabla\Phi_i\|_{L_{x}^{\infty}}^2\|\nabla\Phi_i^{-1}\|_{L_{x}^{\infty}}^2
					\|\RRR_{q,n+1}[h]-\RRR_{q,n+1}[\widetilde h]\|_{L_x^2}
					\\
					&\quad\leq C
					\|\RRR_{q,n+1}[h]-\RRR_{q,n+1}[\widetilde h]\|_{L_x^2}	\\
					&\quad
					\leq C^{2}\tau_q^{-1}\sup_{t}
					\left\|\mathcal R\bigl(w_{i,n+1}[h]-w_{i,n+1}[\widetilde h]\bigr)(t)\right\|_{L_x^2},
				\end{align*}
				where we used
					$\RRR_{q,n+1}=\sum_i\partial_t\bar\eta_i\mathcal Rw_{i,n+1}$, the bounded overlap of the cutoffs, and
					$\|\partial_t\bar\eta_i\|_{L_{t}^{\infty}}\leq C\tau_q^{-1}$.
				
				Consequently,
				\begin{align*}
					\sup_{t}\left\|\mathcal R\bigl(w_{i,n+2}[h]-w_{i,n+2}[\widetilde h]\bigr)\right\|_{L_{x}^{2}}
					&\leq C\tau_q\sup_{t}
					\|A_{k,i,n+1}[h]-A_{k,i,n+1}[\widetilde h]\|_{L_{x}^{2}}\notag\\
					&\leq C^{3}\sup_{t}
					\left\|\mathcal R\bigl(w_{i,n+1}[h]-w_{i,n+1}[\widetilde h]\bigr)(t)\right\|_{L_{x}^{2}}\notag\\
					&\leq C^{3}C_0^{n+1}\tau_q\|h-\widetilde h\|_{L_{t}^\infty}\notag\\
					&\leq C_0^{n+2}\tau_q\|h-\widetilde h\|_{L_{t}^\infty}.
				\end{align*}
				This proves the estimate at level $n+1$ and closes the induction. Iterating up to $n=L-1$ establishes \eqref{app:needed-Rw-bound}.
				
				Finally, \eqref{app:target-reduction},
				$\|\nabla\bar v_q\|_{L_t^\infty L_x^2}\leq\delta_q^{1/2}\lambda_q$, and
				$\tau_q\delta_q^{1/2}\lambda_q=\lambda_q^{-5\alpha}$ give
				\[
				\left\|\int_{\Ttwo}\bar v_q\cdot
				\bigl(w_{q+1}^{(t)}[h]-w_{q+1}^{(t)}[\widetilde h]\bigr)\,\dd x\right\|_{L_{t}^\infty}
				\leq LC_{0}^{L}\lambda_q^{-5\alpha}\|h-\widetilde h\|_{L_{t}^\infty}
				\leq \lambda_q^{- \alpha}\|h-\widetilde h\|_{L_{t}^\infty}
				\]
				for sufficiently large $a$. This proves \eqref{app:picard-C0-final}.
			\end{proof}
			\endgroup
			
			{
				\section{Analytic tools for the convex integration scheme}}\label{Preliminaries}
			
			Let $\alpha\in(0,1)$ and $N\in\mathbb N\cup\{-1\}$. We use the following notation for the H\"older and Besov norms:
			\begin{align*}
				\|f\|_N \triangleq \|f\|_{C^N}, \quad \|f\|_{N+\alpha} \triangleq \|f\|_{B^{N+\alpha}_{\infty,\infty}}.
			\end{align*}
			
			\begin{defn}[Mollifiers]  \label{Moll}
				Choose a spatial mollifier $\psi\in C_c^\infty(\mathbb R^d)$ satisfying
					\[
					\int_{\mathbb R^d}\psi(x)\,\dd x=1,
					\qquad
					\int_{\mathbb R^d}x^\kappa\psi(x)\,\dd x=0,
					\quad1\le|\kappa|\le N_0+2,
					\]
					where $\kappa$ is a multi-index. We also choose a temporal mollifier $\varphi\in C_c^\infty(\mathbb R)$ satisfying
					\[
					\int_{\mathbb R}\varphi(t)\,\dd t=1,
					\qquad
					\int_{\mathbb R}t^n\varphi(t)\,\dd t=0,
					\quad1\le n\le N_0+2.
					\]
					For each $\epsilon>0$, define
				\begin{align}
					{\varphi_\epsilon(t)=\frac1\epsilon\varphi\!\left(\frac t\epsilon\right),
						\qquad
						\psi_\epsilon(x)=\frac1{\epsilon^d}\psi\!\left(\frac x\epsilon\right).} \label{e:defn-mollifier-x}
				\end{align}
			\end{defn}
			
			\begin{defn}[Nonhomogeneous Besov spaces $B^s_{p,q}(\mathbb{T}^2)$ \cite{Ba}]  \label{def.Be}
				Let $s\in\mathbb R$ and $1\leq p,q\leq\infty$. The nonhomogeneous Besov space $B^s_{p,q}(\mathbb T^2)$ consists of all periodic distributions $u\in D'(\mathbb T^2)$ such that
				\begin{equation*}
					{\|u\|_{B^s_{p,q}(\mathbb{T}^2)} \overset{\text{def}}{=} \Big\|\big(2^{js} \|\Delta_j u\|_{L^p(\mathbb{T}^2)}\big)_{j\geq-1}\Big\|_{\ell^q(\{-1,0,1,\ldots\})} < \infty,}
				\end{equation*}
				where  $\Delta_j$ is a localization operator.
			\end{defn}
			
			\begin{defn}[The operator $\mathcal{R}$]  \label{def.R}
				For a smooth vector field $u\in C^\infty(\mathbb T^2;\mathbb R^2)$, the operator $\mathcal R$ introduced in \cite{BDIS15} is defined by
				\begin{align*}
					(\mathcal{R} u)_{ij} = \Delta^{-1} (\partial_i u_j + \partial_j u_i - \div u \delta_{ij}).
				\end{align*}
				This is a matrix-valued right inverse of the divergence on mean-free vector fields:
				\[ \div \mathcal{R} u = u - \fint_{\mathbb{T}^2} u. \]
				Moreover, $\mathcal Ru$ is symmetric and trace-free.
			\end{defn}
			
			We conclude with the standard geometric lemma. See, for instance, \cite{ADS}.
			
			\begin{lem}[Geometric Lemma]  \label{first S}
				Let $r_0=\frac{1}{100}$, and $B_{r_0}(\mathrm{Id})$ denote the ball of radius $r_0$ centered at $\mathrm{Id}$ in the space of $2\times2$ symmetric matrices. There exist two pairwise disjoint sets $\Lambda_i\subset\mathbb S\cap\mathbb Q^2$, $i=1,2$, whose elements $k$ determine orthonormal bases $(k,k^\perp)$, and smooth positive functions $a_{k}$ defined on a neighborhood of  $B_{r_0}(\mathrm{Id})$ such that, for all $R\in B_{r_0}(\mathrm{Id})$,
				\[ R = \sum_{k \in \Lambda_i} a_k^2(R) k \otimes k, \quad i=1,2. \]
				Moreover, $\Lambda_i=-\Lambda_i$ and $a_k=a_{-k}$. There exist constants $C_0=C_0(N_0)$, for $0\leq N\leq N_0+3$ and $x\in\overline{B}_{r_0}(\mathrm{Id})$,
				\begin{equation}
					\|a_k\|_{C^N(\bar{B}_{r_{0}}(\mathrm{Id}))} \leq C_0.
				\end{equation}
			\end{lem}
			
		\end{appendices}
		
		\section*{Acknowledgments}
		The authors would like to thank Jiahong Wu for helpful discussions related to this work. We are grateful to Philip Isett for his valuable comments and suggestions, which helped improve the manuscript. Lili Du was supported by the National Natural Science Foundation of China (Grant No.12125102, 12671245, 12631009). Xinliang Li was supported by the National Natural Science Foundation of China (Grant No.~12401297) and the STU Scientific Research Initiation Grant. Weikui Ye was supported by the National Natural Science Foundation of China (Grant No.~12401277).
		\section*{AI Declaration}
		The research manuscript and preparation of the submission were completed
		 without AI assistance, and the initial manuscript was first posted
		on arXiv on 18 June 2025 (arXiv:2506.15396). 
		\section*{Conflicts of interest}
		The authors declare that they have no conflicts of interest.
		\section*{Data Availability}
		Data sharing is not applicable because no data were generated or analyzed in this study.

\end{document}